\documentclass[11pt, reqno]{amsart}

\usepackage{amsmath, amssymb, comment, amsthm, amssymb}
\usepackage{amscd}
\usepackage{graphicx}
\usepackage{ragged2e}
\usepackage{epstopdf}
\usepackage{setspace}
\usepackage{abstract}
\usepackage{subfigure}
\usepackage{color}

\usepackage{enumitem}
\usepackage{tabularx}
\usepackage[percent]{overpic}

\usepackage[pdftex, linktocpage=true]{hyperref}
\usepackage{todonotes}

\DeclareFontFamily{U}{mathx}{\hyphenchar\font45}
\DeclareFontShape{U}{mathx}{m}{n}{
      <5> <6> <7> <8> <9> <10>
      <10.95> <12> <14.4> <17.28> <20.74> <24.88>
      mathx10
      }{}
\DeclareSymbolFont{mathx}{U}{mathx}{m}{n}
\DeclareFontSubstitution{U}{mathx}{m}{n}
\DeclareMathAccent{\widecheck} {0}{mathx}{"71}

\usepackage[format=plain,labelfont=bf,up,width=.9\textwidth]{caption}

\theoremstyle{plain}
\newtheorem{thm}{Theorem}[section]
\newtheorem{lemma}[thm]{Lemma}
\newtheorem{prop}[thm]{Proposition}

\theoremstyle{definition}
\newtheorem{defn}[thm]{Definition}
\newtheorem{remarka}[thm]{Remark}

\newtheorem{cor}{Corollary}[thm]

\theoremstyle{remark}
\newtheorem{remark}[thm]{Remark}

\newcommand{\R}{\mathbb{R}}
\newcommand{\N}{\mathbb{N}}
\newcommand{\C}{\mathbb{C}}
\newcommand{\Z}{\mathbb{Z}}
\newcommand{\D}{\mathbb{D}}

\newcommand{\bigO}{\mathcal{O}}
\newcommand{\fp}{\mathbf{p}}
\newcommand{\fq}{\mathbf{q}}
\newcommand{\s}{\mathbb{S}}
\newcommand{\Pet}{\mathcal{P}}

\newcommand{\ds}{\displaystyle}
\newcommand{\con}{\tau}
\newcommand{\cU}{c_{U''}}
\newcommand{\rD}{3}

\def\proof{\par\medskip\noindent {\bf Proof.\ \ }}

\def\qed{\hfill $\square$\\ }
\def\qedof{\hfill $\square$\ }

\def \pvec #1#2{\begin{pmatrix}#1\\ \noalign{\vskip  -0 pt} #2\end{pmatrix}}
\def \hvec #1#2{\left(#1,\\ #2\right)\!}

\def \diag #1#2#3#4#5#6#7#8{\begin{CD} #1 @>#5>> #2\\ @V#6VV  @VV#7V\\ #3 @>#8>> #4 \end{CD}}
\def\He{{H\'enon }}

\def\limind{{\underrightarrow{\lim}}}

\title[A structure theorem for semi-parabolic H\'enon maps]{A structure theorem for semi-parabolic\\ H\'enon maps}
\author{Remus Radu and Raluca Tanase}
\date{\today}

\begin{document}
\maketitle

\vspace{.5cm}
\begin{abstract}
\noindent {\sc abstract.}
Consider the parameter space $\mathcal{P}_{\lambda}\subset \mathbb{C}^{2}$ of complex H\'enon maps 
\[
H_{c,a}(x,y)=(x^{2}+c+ay,ax),\ \ a\neq 0
\]
which have a semi-parabolic fixed point with one eigenvalue $\lambda=e^{2\pi i p/q}$. We give a characterization of those H\'enon maps from the curve $\mathcal{P}_{\lambda}$ that are small perturbations of a quadratic polynomial $p$ with a parabolic fixed point of multiplier $\lambda$. 
We prove that there is an open disk of parameters in $\mathcal{P}_{\lambda}$ for which the semi-parabolic H\'enon map has connected Julia set $J$ and is structurally stable on $J$ and $J^{+}$. The Julia set $J^{+}$  has a nice local description: inside a bidisk $\mathbb{D}_{r}\times \mathbb{D}_{r}$  it is a trivial fiber bundle over $J_{p}$, the Julia set of the polynomial $p$, with fibers biholomorphic to $\mathbb{D}_{r}$. The Julia set $J$ is homeomorphic to a quotiented solenoid. 
\end{abstract}
\vspace{.5cm}

\tableofcontents


\newpage

\section{Introduction}\label{sec:intro}
A \He map is a polynomial automorphism of $\C^{2}$ and can be written as
\[
H_{c,a}\hvec{x}{y}=\hvec{x^{2}+c + ay}{ax}\!, \ \ \mbox{for}\ a\neq 0,
\]
where $a$ and $c$ are complex parameters. In this parametrization, the \He map has constant Jacobian $-a^{2}$.
In order to study the dynamics of polynomial automorphisms of $\C^{2}$ we need to understand their behavior under forward and backward iterations. The dynamical objects that we need to analyze are the sets $K^{\pm}$ (the set of points with bounded forward/backward orbits) and their topological boundaries $J^{\pm}=\partial K^{\pm}$. The set $J=J^{+}\cap J^{-}$ is the analogue of the one-dimensional Julia set for polynomials.

We say that the \He map is hyperbolic if it is hyperbolic on its Julia set $J$. If $H_{c,a}$ is hyperbolic and $|a|<1$ then the interior of $K^{+}$ consists of the basins of attraction of finitely many attractive periodic points \cite{BS1}. Each basin of attraction is a Fatou-Bieberbach domain (a proper subset of $\C^{2}$, biholomorphic to $\C^{2}$). The common boundary of the basins is the set $J^{+}$ \cite{BS1}. The set $J^{+}$ is where the most interesting chaotic behavior takes place.
For hyperbolic \He maps, periodic points are dense in $J$ and the map is structurally stable on $J$ \cite{BS1}. When $x\mapsto x^{2}+c$ is a hyperbolic polynomial, the \He map $H_{c,a}$ is also hyperbolic for small values of $a$.
Hyperbolic \He maps that come from perturbations of hyperbolic polynomials are very well understood, by work of Hubbard and Oberste-Vorth \cite{HOV1}, \cite{HOV2} and Forn{\ae}ss and Sibony \cite{FS}.
However, there is very little known about \He maps which are not hyperbolic.

In this paper, we study \He maps with a semi-parabolic fixe point (or cycle). A fixed point of $H_{c,a}$ is called semi-parabolic if the derivative of $H_{c,a}$ at the fixed point has two eigenvalues $\lambda=e^{2\pi i p/q}$ and $\mu$, with $|\mu|<1$.  For clarity and simplicity of exposition, we will call a \He map semi-parabolic if it has a semi-parabolic fixed point.

Unlike hyperbolic \He maps, which exhibit structural stability, semi-parabolic \He maps are not expected to be structurally stable. The general assumption is that bifurcations will occur as we perturb from a semi-parabolic \He map. Bedford, Smillie, and Ueda show in \cite{BSU} some of the complications that can arise by describing the phenomenon of ``semi-parabolic implosion'' in $\C^{2}$ (discontinuity of $J$ and $J^{+}$ on the parameters).
We prove that there are classes of semi-parabolic \He maps that are structurally stable on the sets $J$ and $J^{+}$ inside a parametric region of codimension one in $\C^{2}$. We give a complete characterization of the dynamics of these \He maps. In particular, we show that $J$ is homeomorphic to a solenoid with identifications, hence it is connected.

This parametric region of structural stability will be obtained by considering appropriate perturbations in $\C^{2}$ of a polynomial with a parabolic fixed point of multiplier $\lambda=e^{2 \pi i p/q}$. 
The set of parameters $(c,a)\in\C^{2}$ for which the H\'enon map $H_{c,a}$ has a fixed point with one eigenvalue $\lambda$ is an algebraic curve $\mathcal{P}_{\lambda}$ in $\C^{2}$. The parameter $c=c(a)$ is a function of $a$ and the parametric line $a=0$ intersects the curve $\mathcal{P}_{\lambda}$ at a point $c_{0}$. The polynomial $p(x)=x^{2}+c_{0}$ has a parabolic fixed point of multiplier $\lambda$. Let $J_{p}$ be the Julia set of the parabolic polynomial $p$.

In this article, we will describe the class of semi-parabolic \He maps $H_{c,a}$ where $(c,a)$ lies in a small disk around $0$ inside the curve $\mathcal{P}_{\lambda}$. Consider $r>\rD$ and let $V$ denote the bidisk $\D_{r}\times\D_{r}$ throughout this section. 
We can now state the following theorem. 

\begin{thm}[\textbf{Structure Theorem}]\label{thm:Parabolics-Jp}
Let $p(x)=x^{2}+c_{0}$ be a polynomial with a parabolic fixed point of multiplier $\lambda=e^{2 \pi i p/q}$. There exists $\delta>0$  such that for all parameters $(c,a)\in \mathcal{P}_{\lambda}$ with $0<|a|<\delta$ there exists a homeomorphism
\begin{equation*}
    \Phi: J_{p}\times \D_{r} \rightarrow J^{+}\cap V
\end{equation*}
which makes the diagram
\begin{equation*}
\diag{J_{p}\times \D_{r}}{J^{+}\cap  V}{J_{p}\times \D_{r}}{J^{+}\cap V}
{\Phi}{\psi}{H_{c,a}}{\Phi}
\end{equation*}
commute, where
\begin{equation*}
	\psi(\zeta,z)=\left(p(\zeta),a\zeta - \frac{a^{2}z}{p'(\zeta)}\right).
\end{equation*}
\end{thm}

The map $\psi$ depends on $a$, but we will show in Lemmas \ref{lemma:h-conj1} and \ref{lemma:h-conj2} that all maps $\psi$ are conjugate to each other, for sufficiently small $0<|a|<\delta$. Thus it does not matter which one we use and we can assume that the model map is $\psi(\zeta,z)=\left(p(\zeta),\epsilon\zeta - \frac{\epsilon^{2}z}{p'(\zeta)}\right)$, for some $\epsilon>0$ independent of $a$. The function $\psi$ is a solenoidal map in the sense of \cite{HOV1}; it behaves like angle-doubling in the first coordinate, and contracts strongly in the second coordinate.

Theorem \ref{thm:Parabolics-Jp} shows that $J^{+}\cap V$ is a trivial fiber bundle over $J_{p}$, the Julia set of the parabolic polynomial $p(x)=x^{2}+c_{0}$, with fibers biholomorphic to $\D_{r}$. The set $J^{+}$ is laminated by Riemann surfaces isomorphic to $\C$. In fact, the current $\mu^{+}$ supported on $J^{+}$ defined by Bedford and Smillie in \cite{BS1} is laminar.

\begin{thm}[\textbf{Model for $J$}]\label{cor:J} The Julia set $J$ for the \He map is homeomorphic to a quotiented solenoid
\[
J \simeq \bigcap_{n\geq 0} \psi^{\circ n}(J_{p}\times \D_{r}),
\]
hence connected. Moreover $J=J^{*}$, where $J^{*}$ is the closure of the saddle periodic points.
\end{thm}

\begin{figure}[htb]
\begin{center}
\mbox{\subfigure{
\begin{overpic}[scale=0.33]{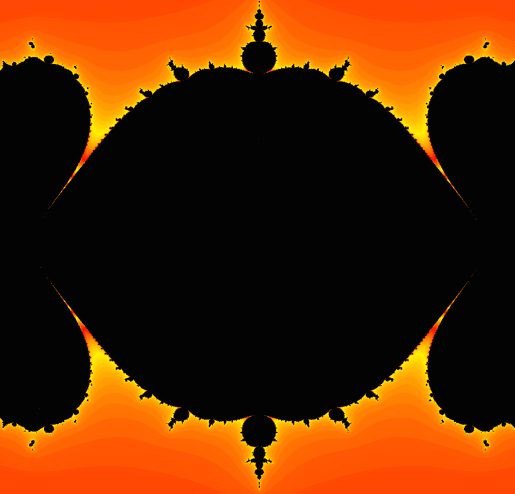}
 \put (50,48.5) {\color{white} {\tiny $\bullet$} {\large $0$} }
\end{overpic}
}
\quad
\subfigure{
\begin{overpic}[scale=0.33]{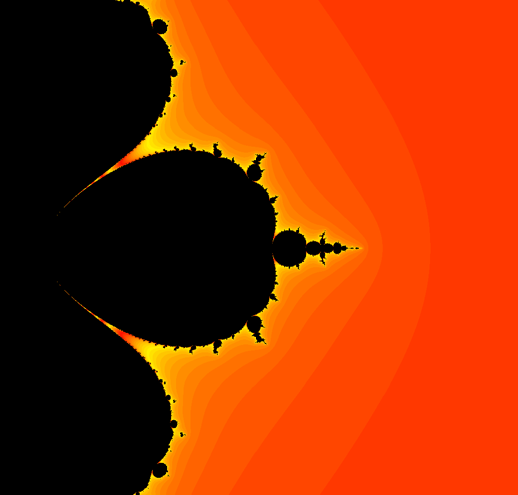}
 \put (37,47) {\color{white} {\tiny $\bullet$} {\large $0$} }
\end{overpic}
}
}
\end{center}
\caption{A parameter plot inside the curve $\mathcal{P}_{-1}$. In both pictures the large region in the center contains the disk $|a|<\delta$. The black region represents (a rough approximation of) the set of parameters $(c,a)\in \mathcal{P}_{-1}$ for which $J$ is connected.  The picture on the left is a double cover of the picture on the right. Both pictures were generated using FractalStream. {\sc left:} The \He map is written as $H_{c,a}(x,y)=(x^{2}+c+ay,ax)$. {\sc right:} The \He map is written in the standard form $H_{c,a}(x,y)=(x^{2}+c-ay,x)$.}
\label{fig:param}
\end{figure}

We describe $J$ as the quotient of a topological solenoid; it is easy to pass from the topological model to a combinatorial description and see that $J$ is equivalent to a dyadic solenoid as in \cite{HOV2} (as a projective limit of $J_{p}$ under the polynomial $p$).

In Corollary \ref{cor:J*} we establish that $J=J^{*}$, the closure of the saddle periodic points. This was not known for this particular class of \He maps and it is still an open question whether $J=J^{*}$ in general. It was shown to be true if $J$ is hyperbolic \cite{BS1}.

Let $f:X\rightarrow X$ be an open, injective map from a space $X$ to itself. Define the inductive limit as $\limind(X,f) :=X\times \N\big{/}\sim$, where the equivalence relation is defined by $(x,n)\sim (f(x),n+1)$. In this setting, the inductive limit is an increasing union of sets homeomorphic to $X$, so locally it looks like $X$. The limit space comes with a natural bijective map $\widecheck{f}:\limind(X,f)\rightarrow \limind(X,f)$ given by $(x,n)\mapsto (f(x),n)$.

Passing to the inductive limit as done by Hubbard and Oberste-Vorth in \cite{HOV2} we get a global model for $J^{+}$.
\begin{thm}[\textbf{Model for $J^{+}$}]\label{cor:globalJ+}
The map $\Phi$ extends naturally to a homeomorphism $\widecheck{\Phi}$ and the following diagram
\[
\diag{\limind(J_{p}\times \D_{r},\psi)}{J^{+}}{\limind(J_{p}\times
\D_{r},\psi)}{J^{+}} {\widecheck{\Phi}}{\widecheck{\psi}}{H_{c,a}}{\widecheck{\Phi}}
\]
commutes.
\end{thm}

As a consequence of the previous theorems, we get that the family of semi-parabolic \He maps $P_\lambda \ni (c,a) \to  H_{c,a}$   is a structurally stable family  on $J$ and $J^{+}$ for  $|a| < \delta$. By structural stability on $J$ and $J^{+}$ we understand the following:

\begin{thm}[\textbf{Stability}]\label{cor:stability}
If $(c_{1},a_{1})$ and $(c_{2},a_{2})$ belong to $\mathcal{P_{\lambda}}$ and if $0<|a_{i}|<\delta$ then $J_{c_{1},a_{1}}$ is homeomorphic to $J_{c_{2},a_{2}}$ and $(H_{c_{1},a_{1}}, J_{c_{1},a_{1}})$ is conjugate to $(H_{c_{2},a_{2}}, J_{c_{2},a_{2}})$. The same is true for $J^{+}$ instead of $J$.
\end{thm}

Let $\lambda=1$ and consider perturbations of the parabolic polynomial $p(x)=x^{2}+1/4$ (the root of the main cardioid of the Mandelbrot set) inside $\mathcal{P}_{1}$. The Julia set $J_{p}$ of this polynomial is a Jordan curve and therefore Theorem \ref{cor:J} implies that the Julia set $J_{c,a}$ of the \He map is homeomorphic to a solenoid. Together with Theorem \ref{cor:stability} this gives a positive answer to some questions of Bedford (Questions 1 and 2 in \cite{B}). Moreover,  the set $J^{+}_{c,a}$ is homeomorphic to a $3$-sphere with a dyadic solenoid removed, for all $(c,a)\in \mathcal{P}_{1}$ and $a$ sufficiently small \cite{R}.

In order to have nontrivial identifications in the description of the Julia set $J$ from Theorem \ref{cor:J} we need to consider  $\lambda=e^{2\pi i p/q}$, different from 1.  To do so, we have generalized a theorem of Ueda \cite{U} and Hakim \cite{Ha} regarding the local normal form around the semi-parabolic fixed point from the case $\lambda=1$ to the case $\lambda=e^{2\pi i p/q}$. This is given in Section \ref{sec: NormalForm}. In Section \ref{sec: Petals} we define big attractive petals for semi-parabolic germs of $(\C^{2},0)$. Both sections are of independent interest.  In Theorem \ref{thm:UniformNF} we show how to control the size of the normalizing neighborhood for our family of semi-parabolic \He maps.

\begin{remarka}\label{remark:J-} We were able to characterize $J$ without using $J^{-}$, by carefully describing the set $J^{+}$ inside the polydisk $\D_{r}\times \D_{r}$. We know from \cite{BS8} that the \He map is hyperbolic on $J$ if and only if there is a neighborhood $\mathcal{N}$ of $J$ and Riemann surface laminations $\mathcal{L}^{\pm}$ of $\mathcal{N}\cap J^{\pm}$ such that $\mathcal{L}^{+}$ and $\mathcal{L}^{-}$ intersect transversely at all points of $J$. In our setting, we have shown that $J^{+}$ is laminar and $J=J^{*}$, but the \He map is semi-parabolic (hence not hyperbolic), so $J^{-}$ is non-laminar or $J^{+}$ and $J^{-}$ may have points of non-transverse intersection. From Theorem \ref{cor:J}, $J$ is connected and by Theorem 1.5 in \cite{Du} it follows that the set $J^{-}-K^{+}$ supports a unique Riemann surface lamination which is uniquely ergodic. In fact, it seems reasonable that $J^{-}$ is non-laminar precisely at the semi-parabolic fixed point.
\end{remarka}

\begin{remarka} Let $(H_{a})_{a\in\D_{\delta}}$ be the family of complex \He maps with a semi-parabolic fixed point with one eigenvalue $\lambda=e^{2\pi i p/q}$ from Theorem \ref{thm:Parabolics-Jp}. It follows from Bedford, Lyubich and Smillie \cite{BLS} that $H_{a}$ admits an invariant measure $\mu_{a}$ which is the unique measure of maximal entropy $\log(2)$. The measure $\mu_{a}$ has two non-zero Lyapunov exponents $\lambda^{-}_{a}<0<\lambda^{+}_{a}$.
Let $J_{a}$ denote the Julia set of $H_{a}$. We have the following dichotomy from \cite{BS5}:  $\lambda^{-}_{a}=\log(2)$ if and only if $J_{a}$ is connected. We have shown in Theorem \ref{cor:J}, that the Julia set $J_{a}$ is connected for each $a\in\D_{\delta}$. Thus $\lambda^{+}_{a}=\log(2)$ and $\lambda^{-}_{a}=2\log|a|-\log(2)$ for this family of semi-parabolic \He maps.
\end{remarka}

In Theorem \ref{thm:Parabolics-Jp} we give a characterization of semi-parabolic \He maps $H_{c,a}$ that are perturbations of a parabolic polynomial $p(x)=x^{2}+c_{0}$. This generalizes the theorem of Hubbard and Oberste-Vorth \cite{HOV2}, which describes H\'enon maps that are perturbations of a hyperbolic polynomial, to the semi-parabolic setting.
The technique of our proof is quite new and is inspired by the proof of Douady and Hubbard \cite{DH}, Section X, (see also \cite{H}) that the Julia set of a parabolic polynomial is locally connected. They show that  the (inverse) B\"ottcher isomorphism extends continuously to the boundary for quadratic polynomials with a parabolic cycle, thus showing local connectivity of the Julia set. We create a two-dimensional analogue to show connectivity of the Julia set for the semi-parabolic \He map.

The key to proving Theorem \ref{thm:Parabolics-Jp} is to build a metric  on a neighborhood of $J^{+}\cap V$ for which the \He map is expanding in the horizontal direction.  
We will consider the infimum of a pull-back of an Euclidean metric in a small tubular neighborhood of the local stable manifold of the semi-parabolic fixed point and 
a product of a Poincar\'e metric with an Euclidean metric outside. The details are given in Section \ref{sec:metric}.
The expanding factor depends on the distance to the local stable manifold of the semi-parabolic fixed point, so it is strictly bigger than one, but there is no constant of uniform expansion.

We prove the result in Theorem \ref{thm:Parabolics-Jp} as a Browder fixed point theorem \cite{Br}. We will recover the set $J^{+}$ inside a bidisk as the image of the unique fixed point of a weakly contracting graph-transform operator in an appropriate function space $\mathcal{F}$. The space $\mathcal{F}$ and the contraction are described in Sections \ref{sec:distance} and \ref{subsec:contraction}.
In order to establish the conjugacy of the semi-parabolic \He map to the model map $\psi$ we have used some heavy-duty topology: a  theorem of Hamstrom \cite{Ham}, which states that if $S$ is a compact surface with nonempty boundary then the components of the group of homeomorphisms which are the identity on the boundary are contractible. This is described in detail in Section \ref{subsec: conjugacy}.
The approach outlined above was used in \cite{R} to reprove the theorem of Hubbard and Oberste-Vorth  about hyperbolic \He maps \cite{HOV2} as an application of the  Banach fixed point theorem.

This article is built on previous work done by the authors in \cite{R} and \cite{T}. We will further use the techniques developed in this paper in \cite{RT} to study perturbations of semi-parabolic \He maps. We show in \cite{RT} that the family of semi-parabolic \He maps $H_a$, where $a$ belongs to an open disk of parameters $|a|<\delta$ from $\mathcal{P}_{\lambda}$, lies in the boundary of a hyperbolic component of the \He connectedness locus. 

\begin{remarka}\label{remark:nonwandering} For the family of semi-parabolic \He maps with small enough Jacobian (suppose $|a|<\delta$ as in Theorem \ref{thm:Parabolics-Jp}) there are no wandering components of $int(K^{+})$. The proof is given in \cite{R} and is similar to the hyperbolic case from \cite{BS2}. This essentially follows from the fact that the \He map expands in horizontal cones on a neighborhood of  $J^{+}\cap V$ as shown in Sections \ref{sec:metric} and \ref{sec:cones}.
\end{remarka}

\noindent {\it Acknowledgements.} We thank John Hubbard for his entire support and guidance with this project, for explaining us the one-dimensional technique used in proving that the Julia set of a parabolic polynomial is locally connected and  for his help in designing a two-dimensional  technique to understand semi-parabolic \He maps. We would also like to thank John Smillie and Eric Bedford for many useful discussions and suggestions.

\section{Preliminaries}\label{sec: Preliminaries}

For a polynomial $p$ of degree $d\geq 2$, the {\it filled Julia set}  of $p$ is
\[
    K_p = \{z\in \C\ |\ |p^{\circ n}(z)|\ \mbox{bounded as}\ n\rightarrow \infty \}.
\]
The set  $J_p=\partial K_p$ is the {\it Julia set} of $p$. If $K_p$ is connected (or equivalently $J_{p}$ is connected) then there exists a unique analytic isomorphism
\begin{equation}\label{eq:Phip}
    \Phi_{p}:\C-\overline{\D} \rightarrow \C-K_p
\end{equation}
such that $\Phi_{p}(z^d)=p(\Phi_{p}(z))$ and $\Phi_{p}(z)/z\rightarrow 1$ as
$z\rightarrow \infty$. Furthermore, if $J_p$ is locally connected then $\Phi_{p}$ extends to the boundary $\s^{1}$ and defines a continuous, surjective map $\gamma:\s^{1} \rightarrow J_p$ \cite{M}. The Julia set of a hyperbolic or parabolic polynomial is locally connected \cite{DH}. 
The boundary map $\gamma$ is called the Carath\'eodory loop. We will use this map in an essential way in Section \ref{subsec: conjugacy}. 

Fix $\lambda=e^{2\pi i p/q}$ a root of unity. The set of parameters $(c,a)\in\C^{2}$ for which the H\'enon map $H_{c,a}$ has a fixed point with one eigenvalue $\lambda$ is a curve of equation
\begin{equation}\label{eq:P-lambda}
\mathcal{P}_{\lambda}:=\left\{(c,a)\in\C^{2}\ | \ c=c(a):=(1-a^{2})\left(\frac{\lambda}{2}-\frac{a^{2}}{2\lambda}\right)-\left(\frac{\lambda}{2}-\frac{a^{2}}{2\lambda}\right)^{2}\right\}\!.
\end{equation}
To see this, let $(x,y)$ be a fixed point of the \He map such that the derivative $DH_{c,a}$ at ${(x,y)}$ has an eigenvalue $\lambda$.
Then  $\lambda$ is a root of the characteristic polynomial $\lambda^{2}-2x\lambda-a^{2}=0$.
The parameters $(c,a)$ must verify the equations $x^{2}+c+ay = x$, $y=ax$ and $x = \frac{\lambda}{2}-\frac{a^{2}}{2\lambda}$.
The solution set is the curve $\mathcal{P}_{\lambda}$.

The parameter $c$ is a quadratic function of the Jacobian $-a^{2}$, so we will refer to the curves $\mathcal{P}_{\lambda}$ as parabolas. For $(c,a)$ in $\mathcal{P}_{\lambda}$,  $H_{c,a}$ has a fixed point $\fq_{a}$ such that $DH_{c,a}(\fq_{a})$ has one eigenvalue $\lambda$ and one eigenvalue $\mu=-\frac{a^{2}}{\lambda}$.  When $|a|<1$, the eigenvalue $\mu$ is  smaller than one in absolute value, and we call $\fq_{a}$ a {\it semi-parabolic fixed point}, and $H_{c,a}$ a {\it semi-parabolic \He map}.
The fixed point has an explicit equation
\begin{equation}\label{eq:fqa}
\fq_{a}:=\left(\frac{\lambda}{2}-\frac{a^{2}}{2\lambda}, a\left(\frac{\lambda}{2}-\frac{a^{2}}{2\lambda}\right)\right)\!.
\end{equation}
We will use this notation throughout this paper. We will see that for $\delta$ small enough and $(c,a)\in \mathcal{P}_{\lambda}$ with $0<|a|<\delta$, the semi-parabolic fixed point $\fq_{a}$ has multiplicity $q+1$ as a solution of the equation $H_{a}^{\circ q}(x,y)=(x,y)$. In analogy with the one dimensional dynamics, we  say that the semi-parabolic multiplicity of $\fq_{a}$ in this case is $1$.

The semi-parabolic fixed point $\fq_{a}$ has a strong stable manifold $W^{s}(\fq_{a})$ biholomorphic to $\C$ corresponding to the eigenvalue $\mu$ whose absolute value is strictly less than $1$,
\begin{equation}\label{eq:Ws-qa}
W^{s}(\fq_{a}):=\{\fp\in\C^{2}\ | \ ||H^{\circ m}(\fp)-\fq_{a}||<C|\mu|^{m}\ \mbox{for}\ m\geq 0 \},
\end{equation}
where $C>0$ is a constant  \cite{U}. This is the set of points for which $H^{\circ m}(\fp)\rightarrow \fq_{a}$ exponentially as $m\rightarrow \infty$.
Bedford, Smillie and Ueda show that $W^{s}(\fq_{a})$ is dense in $J^{+}$ \cite{BSU}. The basin of attraction of the semi-parabolic fixed point $\fq_{a}$ belongs to $int(K^{+})$ and is a Fatou-Bieberbach domain  \cite{Ha}, \cite{U}. The rate of convergence to $\fq_{a}$ is parabolic.

The parametric line $a=0$ intersects the curve $\mathcal{P}_{\lambda}$ at the point $c_{0}=\frac{\lambda}{2}-\frac{\lambda^{2}}{4}$. It is easy to see that the polynomial $p(x)=x^{2}+c_{0}$ has a parabolic fixed point $q_{0}=\frac{\lambda}{2}$ of multiplier $\lambda$. 
This is the polynomial from which we are perturbing in $\C^{2}$. 
Let $J_{p}$ and $K_{p}$ denote the Julia set, respectively the filled-in Julia set of the polynomial
$p$.

For $(c,a)\in \mathcal{P}_{\lambda}$ the equation for $c$ from \ref{eq:P-lambda} can be rewritten as
\begin{equation}\label{eq:w}
c =  \frac{\lambda}{2}-\frac{\lambda^{2}}{4} +a^{2}w,\ \ \ \mbox{where}\ \ w :=  \frac{2\lambda-2\lambda^{2}-1}{4\lambda}+a^{2}\frac{2\lambda-1}{4\lambda^{2}}.
\end{equation}
Thus we can also write the semi-parabolic \He map $H_{c,a}(x,y)=(x^{2}+c+ay,ax)$ as
\begin{equation}\label{eq:Hw}
H_{a}\hvec{x}{y} = \hvec{p(x)+a^{2}w+ay}{ax},
\end{equation}
with inverse
\begin{equation}\label{eq:Hw-1}
H_{a}^{-1}\left(x, y\right) = \frac{1}{a}\hvec{y}{x-p\left(y/a\right)-a^{2}w},
\end{equation}
where $p$ is the parabolic polynomial $p(x)=x^{2}+c_{0}$. This emphasizes the dependency on the parabolic polynomial $p$. The constant $w$ depends only on $a$ and $\lambda$ and clearly $|w|<2$. 

Following \cite{HOV1}, we choose a constant $r$ greater than the largest root of the quadratic equation $|x|^{2}-(|a|+2)|x|-|c_{0}|-|a|^{2}|w|=0$. 
Then the dynamical space $\C^2$ can be divided
into three regions: the bidisk
$ \D_{r}\times\D_{r}=\{(x,y)\in \C^2\ |\ |x|\leq r, |y|\leq r\},$
\[
V^+=\{(x,y)\ |\ |x|\geq \max(|y|,r) \}\ \
\mbox{and}\ \ V^-=\{(x,y)\ |\ |y|\geq \max(|x|, r)\}.
\]
The sets $J$ and $K$ are contained in $\D_{r}\times\D_{r}$. The escaping sets are $U^{+}=\C^{2}-K^{+}$ and $U^{-}=\C^{2}-K^{-}$ and they can be
described in terms of $V^+$ and $V^-$ as follows \cite{HOV1}:
\[
U^+=\bigcup_{k\geq 0}
H^{-\circ k}(V^+) \ \ \mbox{and}\  \  U^-=\bigcup_{k\geq 0} H^{\circ k}(V^-).
\]
Therefore the Julia set $J^{+}$ is the common boundary of $K^{+}$ and $U^{+}$. We know that $|c_{0}|<2$ from \cite{DH} so any constant $r>3$ works for our purpose. 

In order to understand these objects better it is useful to look first at the case when $a=0$. In this case the dynamics of the \He map reduces to the dynamics of the polynomial $p$. The relevant sets under forward dynamics can then be easily described: $J^{+}=J_{p}\times \C$, $K^{+}=K_{p}\times \C$ and $U^{+}=(\C-K_{p})\times \C$.

\section{Normal form of semi-parabolic \He maps}\label{sec: NormalForm}

Hakim in \cite{Ha} and Ueda in \cite{U} have studied normal forms for germs of semi-attractive transformations $H$ of $(\C^{n},0)$ for which $DH_{(0)}$ has one eigenvalue $\lambda=1$, and the other eigenvalues $\mu_{2},\ldots \mu_{n}$ have absolute values $|\mu_{j}|<1, j=2,\ldots,n$.

The following results are similar to Proposition 2.1, 2.2, and 2.3 from \cite{Ha} and to Section 6 from \cite{U}. We have adapted the propositions in \cite{Ha} to semi-parabolic germs of transformations of $(\C^{2},0)$ with eigenvalues $\lambda=e^{2\pi i p/q}$ and $|\mu|<1$. As a consequence we get that $0$ is a fixed point with multiplicity $\nu q+1$ for some constant $\nu$ which we call the {\it (semi) parabolic multiplicity} of the fixed point, like in one-dimensional dynamics.

\begin{prop}\label{thm:normalform1}
Let $H$ be a semi-parabolic germ of transformation of $(\C^{2},0)$, with eigenvalues $\lambda$ and $\mu$, with $\lambda=e^{2\pi i p/q}$ and $|\mu|<1$. There exist local coordinates $(x,y)$  in which $H$ has the form $H(x,y)=(x_{1},y_{1})$, with
\begin{equation}\label{eq:NF1}
\left\{\begin{array}{l}
    x_{1}= a_{1}(y)x+a_{2}(y)x^{2}+\ldots \\
    y_{1}=\mu y + xh(x,y)
\end{array}\right.
\end{equation}
where $a_{j}(\cdot)$ and $h(\cdot,\cdot)$ are germs of holomorphic functions from $(\C,0)$ to $\C$, respectively from $(\C^{2},0)$ to $\C$, with $a_{1}(0)=\lambda$ and $h(0,0)=0$.
\end{prop}
\proof The proof is the same as in \cite{Ha} and \cite{U} and is based on the straightening of the local strong stable manifold of the fixed point. \qed

\begin{prop}\label{thm:normalform2}
Let $H$ be a semi-parabolic germ of transformation of $(\C^{2},0)$, with eigenvalues $\lambda$ and $\mu$, with $\lambda=e^{2\pi i p/q}$ and $|\mu|<1$. For every integer $m$ there exist local coordinates $(x,y)$  in which $H$ has the form $H(x,y)=(x_{1},y_{1})$, with
\begin{equation}\label{eq:NF2}
\left\{\begin{array}{l}
    x_{1}= \lambda x+a_{2}x^{2}+\ldots + a_{m}x^{m} + a_{m+1}(y)x^{m+1}+\ldots \\
    y_{1}= \mu y + xh(x,y)
\end{array}\right.
\end{equation}
where $a_{2},\ldots, a_{m}$ constants.
\end{prop}
\proof The proof is the same as in Proposition 2.2 from \cite{Ha} (proved also in Section 6 of \cite{U}). We will refer to this proof when we discuss the domain of convergence of the functions $u(\cdot)$ and $v(\cdot)$ defined below. We know from Proposition \ref{thm:normalform1} that there exist local coordinates $(x,y)$ in which $H$ has the form
\begin{equation*}
\left\{\begin{array}{l}
    x_{1}= a_{1}(y)x+a_{2}(y)x^{2}+\ldots \\
    y_{1}=\mu y + xh(x,y).
\end{array}\right.
\end{equation*}
The germs $a_{i}(\cdot)$ and $h(\cdot,\cdot)$ germs of holomorphic functions from $(\C,0)$ to $\C$, respectively from $(\C^{2},0)$ to $\C$, with $a_{1}(0)=\lambda$ and $h(0,0)=0$.

\noindent $(1)$ Reduction to $a_{1}(y)=\lambda$. Consider as in \cite{Ha} and \cite{U} a coordinate transformation
\begin{equation*}
\left\{\begin{array}{l}
    X=u(y)x \\
    Y=y
\end{array}\right. \ \ \ \mbox{with inverse}\ \ \
\left\{\begin{array}{l}
    x=X/u(Y) \\
    y=Y
\end{array}\right.
\end{equation*}
where $u$ is a germ of analytic functions from $(\C,0)$ to $\C$ with $u(0)=\lambda$.  We need to find $u$ such that
\begin{eqnarray*}
X_{1}&=& u(y_{1})x_{1}=u(\mu y +x h(x,y))\left(a_{1}(y)x+a_{2}(y)x^{2}+\ldots\right)\\
    &=& u(\mu Y + X/u(Y)h(X/u(Y) ,Y) )\left(a_{1}(Y)X/u(Y) +a_{2}(Y)(X/u(Y))^{2}+\ldots\right)\\
    &=& \frac{u(\mu Y)a_{1}(Y)}{u(Y)}X + \bigO(X^{2}) = \lambda X + \bigO(X^{2}).
\end{eqnarray*}
Thus $u$ satisfies the equation $u(Y) = u(\mu Y)\frac{a_{1}(Y)}{\lambda}$. We successively substitute $\mu Y$ instead of $Y$ in this equation and obtain the unique solution
\begin{equation}\label{eq:u(Y)}
u(Y) = \prod\limits_{n=0}^{\infty} \frac{a_{1}(\mu^{n}Y)}{\lambda}.
\end{equation}
This series converges in a neighborhood of $0$ since $\mu<1$ and $a_{1}(0) = \lambda$.

\noindent $(2)$ Reduction to $a_{2}(y),\ldots, a_{m}(y)$ constants. We proceed by induction on $m$. The base case $m=1$ was discussed above. Suppose that $m\geq 2$ and that  there exist local coordinates $(x,y)$ in which $H$ has the form
\begin{equation*}
\left\{\begin{array}{l}
    x_{1}= \lambda x+a_{2}x^{2}+\ldots + a_{m-1}x^{m-1} + a_{m}(y)x^{m}+\ldots \\
    y_{1}= \mu y + xh(x,y),
\end{array}\right.
\end{equation*}
with $a_{2}, \ldots , a_{m-1}$  constant. We would like to find local coordinates so that $a_{m}(y)$ is also constant.  Consider the transformation
\begin{equation*}
\left\{\begin{array}{l}
    X=x+v(y)x^{m} \\
    Y=y
\end{array}\right. \ \ \ \mbox{with inverse}\ \ \
\left\{\begin{array}{l}
    x=X-v(Y)X^{m} + \ldots \\
    y=Y
\end{array}\right.
\end{equation*}
where $v$ is a germ of analytic functions from $(\C,0)$ to $\C$ with $v(0)=0$.  Using the coordinates given by this transformation we get
\begin{eqnarray*}
X_{1}&=& x_{1} + v(y_{1})x_{1}^{m}\\
    &=& \lambda x+a_{2}x^{2}+\ldots + a_{m-1}x^{m-1} + \left(a_{m}(y)+v(\mu y)\right)x^{m}+\bigO(x^{m+1})\\
    &=& X - v(Y)X^{m} + a_{2}X^{2} + \ldots + a_{m-1}X^{m-1}+\left(a_{m}(Y) + v(\mu Y)\right)X^{m} + \bigO(X^{m+1})\\
    &=& X + a_{2}X^{2} + \ldots + a_{m-1}X^{m-1}+\left( a_{m}(Y) + v(\mu Y) - v(Y)\right)X^{m} + \bigO(X^{m+1}).
\end{eqnarray*}
We need $v$ such that the coefficient of $X^{m}$ is constant, i.e. $ a_{m}(Y) + v(\mu Y) - v(Y)=a_{m}(0)$ is constant. This gives the equation $v(Y)-v(\mu Y)=a_{m}(Y)-a_{m}(0)$. We successively substitute $\mu Y$ instead of $Y$ in this equation and obtain
\begin{equation}\label{eq:v(Y)}
v(Y) = \sum\limits_{n=0}^{\infty}( a_{m}\left( \mu^{n}Y)-a_{m}(0) \right).
\end{equation}
The series clearly converges in a neighborhood of $0$ since $\mu<1$.
\qed

\begin{prop}\label{thm:normalform3}
Let $H$ be a semi-parabolic germ of transformation of $(\C^{2},0)$, with eigenvalues $\lambda$ and $\mu$, with $\lambda=e^{2\pi i p/q}$ and $|\mu|<1$. There exist local coordinates $(x,y)$  in which $H$ has the form $H(x,y)=(x_{1},y_{1})$, with
\begin{equation}\label{eq:NF3}
\left\{\begin{array}{l}
    x_{1}= \lambda (x+x^{\nu q+1} + Cx^{2\nu q+1}+ a_{2\nu q+2}(y)x^{2\nu q +2}+\ldots )\\
    y_{1}= \mu y + xh(x,y)
\end{array}\right.
\end{equation}
and $C$ a constant. Moreover the multiplicity of the fixed point is $\nu q +1$.
\end{prop}
\proof
Suppose that the map has the form from Equation \ref{eq:NF2}, where $m$ is big enough, and fixed
\begin{equation*}
\left\{\begin{array}{l}
    x_{1}= \lambda x+a_{k}x^{k}+\ldots + a_{m}x^{m} + a_{m+1}(y)x^{m+1}+\ldots \\
    y_{1}= \mu y + xh(x,y).
\end{array}\right.
\end{equation*}
Consider the coordinate transformation
\begin{equation*}
\left\{\begin{array}{l}
    X=x+bx^{k} \\
    Y=y
\end{array}\right. \ \ \ \mbox{with inverse}\ \ \
\left\{\begin{array}{l}
    x=X-bX^{k}+\ldots \\
    y=Y
\end{array}\right.
\end{equation*}
In the new coordinate system, we get
\begin{eqnarray*}
X_{1}&=&x_{1}+bx_{1}^{k} = ( \lambda x+a_{k}x^{k}+\ldots )+ b(\lambda x+a_{k}x^{k}+\ldots )^k\\
    &=& \lambda x+a_{k}x^{k}+\ldots + b\lambda^{k}x^{k} +\ldots \\
    &=& \lambda x +(a_{k}+b\lambda^{k})x^{k} + \ldots\\
    &=&\lambda ( X -  b X^{k} + \ldots) +(a_{k}+b\lambda^{k})(X-bX^{k} + \ldots)^{k} + \ldots\\
    &=&\lambda X  +  (a_{k}+b(\lambda^{k}-\lambda)) X^{k}+ \ldots
\end{eqnarray*}
If $k$ is not congruent to $1$ modulo $q$ (i.e. $\lambda^{k}\neq \lambda$), then we can set
\[b = \frac{a_{k}}{\lambda-\lambda^{k}}
\]
and eliminate the term $a_{k}x^{k}$. This proves that by successive coordinate transformations of the form $X=x+bx^{k}, Y=y$ we can eliminate terms with powers that are not congruent to $1$ modulo $q$, so the first term that cannot be eliminated in this way will have a power of the form $\nu q + 1$ for some $\nu$.

Thus the map takes the form
\begin{equation}\label{eq:NF3a}
\left\{\begin{array}{l}
    x_{1}= \lambda (x+a_{\nu q+1}x^{\nu q+1} + \ldots + a_{m}x^{m} + a_{m+1}(y)x^{m+1}+\ldots) \\
    y_{1}= \mu y + xh(x,y).
\end{array}\right.
\end{equation}
Here we will assume that $m$ was chosen so that $m>2\nu q+1$. Thus the coefficients up to order $m$ are still constants.  We can further reduce Equation \ref{eq:NF3a} to $a_{\nu q+1}=1$ by considering a transformation of the form $X=Ax, Y=y$, where $A$ is a constant such that $A^{\nu q}=a_{\nu q+1}$. Consider therefore the transformation $H$ written as
\begin{equation}\label{eq:NF3b}
\left\{\begin{array}{l}
    x_{1}= \lambda (x+x^{\nu q+1} + \ldots + a_{m}x^{m} + a_{m+1}(y)x^{m+1}+\ldots) \\
    y_{1}= \mu y + xh(x,y).
\end{array}\right.
\end{equation}

We have previously showed that we can eliminate any term of degree $k$ between $1$ and $m$,
which is not congruent to $1$ modulo $q$. One can also eliminate all terms of degree $jq+1$, where $\nu<j<2\nu$. Assume $a_{jq+1}$ is the first such coefficient different from $0$ as below
\begin{equation}\label{eq:NF3c}
\left\{\begin{array}{l}
    x_{1}= \lambda (x+x^{\nu q+1} + a_{jq+1}x^{jq+1} + \ldots )\\
    y_{1}= \mu y + xh(x,y).
\end{array}\right.
\end{equation}
Consider the transformation
\begin{equation*}
\phi\pvec{x}{y}=\pvec{X}{Y}\ \ \mbox{where} \ \
\left\{\begin{array}{l}
    X=x+bx^{(j-\nu)q+1} \\
    Y=y
\end{array}\right. \ \ \mbox{and}\ \ \
b=\frac{a_{jq+1}}{(2\nu-j)q}.
\end{equation*}

Define $G:=\phi\circ H\circ \phi^{-1}$ and suppose that $G(X,Y)=(X_{1},Y_{1})$ with
\begin{equation}\label{eq:NF3d}
\left\{\begin{array}{l}
    X_{1}= \lambda (X+X^{\nu q+1}  + AX^{k}+\ldots )\\
    Y_{1}= \mu Y + Xh(X,Y)
\end{array}\right.
\end{equation}
and $k\leq jq$. We show that $A=0$ by comparing the terms of the power series of $G\circ \phi$ and $\phi\circ H$. We will only need to analyze the $x$-coordinate. The first coordinate of $\phi\circ  H$ is
\begin{eqnarray*}
X_{1}&=&x_{1}+bx_{1}^{(j-\nu)q+1}\\
    &=& \lambda(x+x^{\nu q+1} + a_{jq+1}x^{jq+1} + \ldots)+ b\lambda( x+x^{\nu q+1} + a_{jq+1}x^{jq+1} + \ldots )^{(j-\nu)q+1}\\
    &=& \lambda \left(x+bx^{(j-\nu)q+1}+x^{\nu q+1}+(a_{jq+1}+b((j-\nu)q+1))x^{jq+1}+\bigO_{y}(x^{jq+2})\right).
\end{eqnarray*}
The first coordinate of $G\circ  \phi$ is
\begin{eqnarray*}
X_{1}&=& \lambda (X+X^{\nu q+1}  + AX^{k}+\ldots )\\
    &=&\lambda \left( (x+bx^{(j-\nu)q+1})+(x+bx^{(j-\nu)q+1} )^{\nu q+1}  + A(x+bx^{(j-\nu)q+1} )^{k}+\ldots \right)\\
    &=& \lambda \left(x+bx^{(j-\nu)q+1}+x^{\nu q+1}+b(\nu q +1)x^{jq+1}+ Ax^{k}+\bigO_{y}(x^{k+1})\right).
\end{eqnarray*}
We have that $a_{jq+1}+b((j-\nu)q+1)=b(\nu q +1)$ by the choice for $b$. The two power series are equal so the coefficient of $x^{k}$ vanishes, so $A=0$. Thus in the first coordinate of $\phi\circ H\circ \phi^{-1}$ the coefficient of $x^{jq+1}$ is zero and the coordinate transformation did not introduce additional terms of lower powers.

Using similar transformations we can eliminate all terms between $\nu q+1$ and $2\nu q +1$
and write $H(x,y)=(x_{1},y_{1})$ with
\begin{equation*}\label{eq:NF3e}
\left\{\begin{array}{l}
    x_{1}= \lambda (x+x^{\nu q+1} + Cx^{2\nu q+1}+ \bigO_{y}(x^{2\nu q+2}))\\
    y_{1}= \mu y + xh(x,y)
\end{array}\right.
\end{equation*}
for some constant $C$.

It is easy to prove that in the last coordinate system, $H^{\circ q}$ takes the form
\begin{equation}\label{eq:NF4}
\left\{\begin{array}{l}
    x_{1}= x+qx^{\nu q+1} + \widetilde{C}x^{2\nu q +1}+\bigO_{y}(x^{2\nu q+2}) \\
    y_{1}= \mu^{q} y + x\widetilde{h}(x,y).
\end{array}\right.
\end{equation}
The partial derivative $\frac{\partial y_{1}}{\partial y}(0,0)=\mu^{q}<1$, hence by the Implicit Function Theorem, the equation $\mu^{q} y + xh(x,y)=y$ has a unique solution $y=\varphi(x)$ in a neighborhood of $0$, where $\varphi$ is a holomorphic function.  From the first equation it then follows that $x=0$ is a fixed point of $H^{\circ q}$ of multiplicity $\nu q+1$.
\qed

The normalizing form as proven in the previous theorem holds locally around the semi-parabolic fixed point. The disadvantage of the ``local'' statement is that it does not allow us to control the size of the neighborhood of the fixed point where we can put on normalizing coordinates. However, in Section \ref{subsec: U-normalizing} we show how to control the size of this neighborhood.
We  consider a class of semi-parabolic \He maps which are perturbations of a polynomial with a parabolic fixed point or cycle, and show how to extend this theorem in order to get uniform bounds (with respect to the parameters) on the size of the normalizing neighborhood.

\section{Attracting and repelling sectors}\label{sec: Petals}

Set $m:=\nu q$ and let
\[
\Delta_{R}= \left\{ x\in \C\ \big{|}\ \left(Re(x^{m})+\frac{1}{2R}\right)^{2} + \left(|Im(x^{m})|-\frac{1}{2R}\right)^{2} < \frac{1}{2R^{2}}\right\}.
\]

There are $m$ connected components of $\Delta_{R}$, which we denote $\Delta_{R,j}$, for $1\leq j\leq m$. Define
\[
\Pet_{att}= \left\{ (x,y)\in \C^{2}\ |\ x\in \Delta_{R}, |y|<r \right\}
\]
and let
\[
\Pet_{att,j}= \left\{ (x,y)\in \C^{2}\ |\ x\in \Delta_{R,j}, |y|<r \right\}
\]
be the connected components of $\Pet_{att}$. These are called {\it (big)  attractive petals} for the \He map, similar to the one-dimensional case.

\begin{prop}\label{prop:sectors} For $R$ large enough and $r$ small enough
\[
H(\overline{\Pet_{att,j}})\subset \Pet_{att,j+\nu p}\cup \{0\}\times \D_{r}\ \ \ \ \mbox{for} \ 1 \leq j \leq \nu q.
\]
In particular $H(\overline{\Pet_{att}})\subset \Pet_{att}\cup \{0\}\times \D_{r}$ and all points of $\Pet_{att}$ are attracted to the origin under iterations by $H$.
\end{prop}
\proof The analysis is similar to \cite{Ha}, but one should have in mind the formalism from the one-dimensional case (see \cite{DH}, \cite{BH}) to resolve the ambiguity about which branch of $x^{1/m}$ we are talking about.

Assume that $R$ is large enough and $r$ is small enough so that $H$ is well defined and has the expansion from Proposition \ref{thm:normalform3}.
Define the region $U_{R_{1}}$
\[
U_{R_{1}}:=\left\{ X\in \C \ |\ R_{1}-Re(X)<|Im(X)|\right\}
\]
where $R_{1}=R/m$ and set $W_{R_{1},r}:=U_{R_{1}}\times \D_{r}\subset\C^{2}$.

Consider the \He map $H$ written as
\begin{equation*}
\left\{\begin{array}{l}
    x_{1}= \lambda (x+x^{m+1} + Cx^{2m+1}+ a_{2m+2}(y)x^{2m +2}+\ldots )\\
    y_{1}= \mu y + xh(x,y).
\end{array}\right. 
\end{equation*}
Suppose $(x,y)\in \Pet_{att,j}$ and consider the transformation
\begin{equation*}
\left\{\begin{array}{l}
    \ds X=-\frac{1}{mx^{m}} \\
    Y=y.
\end{array}\right. 
\end{equation*}
It maps each $\Pet_{att,j}$ to $W_{R_{1},r}$ (it maps points $(0,y)$  to $(\infty,y)$). Let $\widehat{H}(X,Y)=(X_{1},Y_{1})$ be the map in these coordinates
\begin{eqnarray*}
X_{1}&=&-\frac{1}{mx_{1}^{m}}= - \frac{1}{m\left(\lambda (x+x^{m+1} + Cx^{2m+1}+ a_{2m+2}(y)x^{2m +2}+\ldots )\right)^{m}}\\
    &=& \frac{X}{\left(1+x^{m} + Cx^{2m}+ a_{2m+2}(y)x^{2m +1}+\ldots \right)^{m}}\\
    &=&X\left(1-m(x^{m} + Cx^{2m}+\ldots) + \frac{m(m+1)}{2}x^{2m}+\ldots \right)\\
    &=&X+1+\frac{A}{X}+\bigO_{Y}\left(\frac{1}{|X|^{1+1/m}}\right)
\end{eqnarray*}
where $ A := \frac{1}{m}\left(\frac{m+1}{2}-C\right)$ is a constant. The notation $\bigO_{Y}\left(\frac{1}{|X|^{\alpha}}\right)$ represents a holomorphic function of $(X,Y)$ in $W_{R_{1},r}$ which is bounded by $\frac{K}{|X|^{\alpha}}$ for some constant $K$.  Similarly
\begin{eqnarray*}
Y_{1} &=&\mu y + xh(x,y) = \mu Y + \bigO_{Y}\left(\frac{1}{|X|^{1/m}}\right).
\end{eqnarray*}
Note that $|X|>\frac{R_{1}}{\sqrt{2}}$ for all $X\in U_{R_{1}}$. There exists constants $K'$ and $K''$ such that
\begin{eqnarray*}
|X_{1}-X-1|&\leq& \frac{K'}{|X|}<\frac{K_{1}}{R_{1}}\ \ \ \mbox{where}\ \ K_{1}:=K'\sqrt{2}\\
|Y_{1}-\mu Y|&\leq& \frac{K''}{|X|^{1/m}}<\frac{K_{2}}{R_{1}^{1/m}} \ \ \ \mbox{where}\ \ K_{2}:=K''\sqrt{2}^{1/m}.
\end{eqnarray*}
Choose $R_{1}$ large enough and $r$ small enough so that
\begin{equation}\label{eq:K1K2}
\left\{\begin{array}{ll}
	\frac{K_{1}}{R_{1}}<\frac{1}{2}  \\
	\vspace{-0.37cm}\\
	\frac{K_{2}}{R_{1}^{1/m}}<(1-|\mu|)r \ .
\end{array}\right.
\end{equation}
The first condition gives $|X_{1}-X-1|<\frac{1}{2}$, which implies $Re(X_{1})> Re(X)+\frac{1}{2}$ and $|Im(X_{1})|>|Im(X)|-\frac{1}{2}$. Thus $R_{1}-Re(X_{1})<|Im(X_{1})|$.
The second condition gives
\[
	|Y_{1}|\leq |Y_{1}-\mu Y|+|\mu| |Y|<\frac{K_{2}}{R_{1}^{1/m}}+|\mu| r <r.
\]
Hence $\widehat{H}(W_{R_{1},r})\subset W_{R_{1},r}$.

We need to show that points in $W_{R_{1},r}$ are attracted by $(\infty,0)$ under iterations by $\widehat{H}$. Let $(X,Y)\in W_{R_{1},r}$ and set $(X_{n},Y_{n})=\widehat{H}^{\circ n}(X,Y)$. Assume without loss of generality that $Re(X)> \rho$, where $\rho>0$ is a constant to be defined later. We can make this assumption since $Re(X_{k})> Re(X)+\frac{k}{2}$ for every positive integer $k$. We take the first integer $k_{0}$ such that $Re(X_{k_{0}})> \rho$ and let $X:=X_{k_{0}}$ and $Y:=Y_{k_{0}}$.

Clearly
\begin{equation}\label{eq:ReXn}
Re(X_{n})> \rho + \frac{n}{2}
\end{equation}
for every $n\geq 0$. This follows immediately by induction since
\[
Re(X_{n+1})> Re(X_{n})+1/2>\rho+(n+1)/2.
\]

We now show by induction that
\begin{equation*}
|Y_{n}|<2NrR_{1}^{1/m}\left(\frac{1}{\rho+\frac{n}{2}}\right)^{1/m},\ \  n\geq 0
\end{equation*}
where $N$ is an integer number such that $NR_{1}^{1/m}>\rho^{1/m}$. When $n=0$, $|Y|<r$ and
\[
r<2NrR_{1}^{1/m}\frac{1}{\rho^{1/m}} \Leftrightarrow \rho^{1/m}<2NR_{1}^{1/m}.
\]

We now proceed to the induction step. First note that $|X_{n}|\geq Re(X_{n})>\rho+\frac{n}{2}$ and $K''<K_{2}$. We get
\begin{eqnarray*}
|Y_{n+1}|&\leq& |Y_{n+1}-\mu Y_{n}|+|\mu| |Y_{n}|<\frac{K''}{|X_{n}|^{1/m}}+|\mu||Y_{n}| \\
&<& (K_{2}+ |\mu| 2NrR_{1}^{1/m}) \left(\frac{1}{\rho+\frac{n}{2}}\right)^{1/m}<\left(1+(2N-1)|\mu|\right)rR_{1}^{1/m}\left(\frac{1}{\rho+\frac{n}{2}}\right)^{1/m}
\end{eqnarray*}
and we want to show that
\[
|Y_{n+1}|<2NrR_{1}^{1/m}\left(\frac{1}{\rho+\frac{n+1}{2}}\right)^{1/m}.
\]
This inequality is satisfied if
\begin{equation}\label{eq:mu1}
\left(\frac{\rho+\frac{1}{2}}{\rho}\right)^{1/m}=\left(1+\frac{1}{2\rho}\right)^{1/m}<\frac{2}{1+|\mu|}\leq\frac{2N}{1+(2N-1)|\mu|}.
\end{equation}
But $|\mu|<1$, so $2/(1+|\mu|)>1$. This allows us to choose a number $\rho$ large enough so that Equation \ref{eq:mu1} is satisfied. Then choose an integer $N$ such that $NR_{1}^{1/m}>\rho^{1/m}$.

It follows that $(X_{n},Y_{n})\rightarrow (\infty,0)$ as $n\rightarrow \infty$.
\qed

\begin{figure}[htb]
\begin{center}
\includegraphics[scale=0.375]{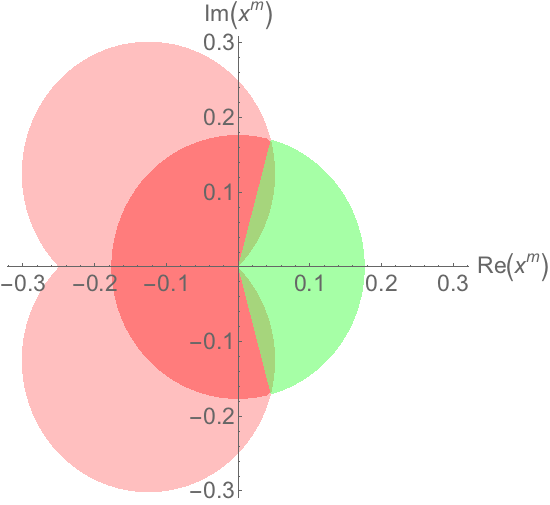}
\end{center}
\caption{The image of $\Pet_{att}$ under the map $x\mapsto x^{m}$ at height $y=0$ is shown in light red. Similarly the attracting sector $\Delta^{+}$ is shown in red and the repelling sector $\Delta^{-}$ in green. The angle opening of the green region is $5\pi/6$.}
\label{fig:Sectors}
\end{figure}

Let $\epsilon_{0}=\tan(\pi/12)$.  
Define the {\it attractive sectors}
\begin{equation}\label{eq:W+}
	\Delta^{+} := \left\{ x\in \C\ |\ Re(x^{m})\leq \epsilon_{0}|Im(x^{m})| \ \mbox{and}\ |x^{m}|< \frac{1}{\sqrt{2}R} \right\}
\end{equation}
and {\it repelling sectors}
\begin{equation}\label{eq:W-}
\Delta^{-}:= \left\{ x\in \C\ |\ Re(x^{m})> \epsilon_{0}|Im(x^{m})|\ \mbox{and} \ |x^{m}|< \frac{1}{\sqrt{2}R} \right\}.
\end{equation}
Let $W^{+}=\Delta^{+}\times\D_{r}\subset \Pet_{att}\cup\{0\}\times\D_{r}$ and $W^{-}=\Delta^{-}\times\D_{r}$. 
We call $W^{-}$ repelling because as we will see, the \He map expands horizontally when the Jacobian is small enough.
There are $m$ components of $W^{\pm}$ which we denote $W^{\pm}_{j}$ for $1\leq j\leq m$. These are the preimages of the red/green regions in Figure \ref{fig:Sectors} under $x\mapsto x^{m}$. The choice of $\epsilon_{0}$ means that the angle of the image of $W^{-}$ under $x\mapsto x^{m}$ is $5\pi/6$.

Furthermore, since $Re(x^{m})> \epsilon_{0}|Im(x^{m})|$ on $W^{-}$, we have
\begin{equation}\label{eq:eps1}
	Re(x^{m})> \epsilon_{1}|x^{m}|,\ \ \ \mbox{where}\ \  \epsilon_{1} := \frac{\epsilon_{0}}{\sqrt{1+\epsilon_{0}^{2}}}>\frac{1}{4}.
\end{equation}
\begin{prop}\label{prop:int(K+)} The interior of the union $\bigcup_{n\geq 0}H^{-\circ n}(W^{+})$ is the basin of attraction of the semi-parabolic fixed point.
\end{prop}
\proof The proof follows immediately by analyzing the situation at infinity using Equation \ref{eq:ReXn} as in the proof of Proposition \ref{prop:sectors}. \qed

\section{The parametrizing map of the stable manifold}\label{subsec:WS}

This is a self-contained section where we study the degeneracy of the parametrization of the stable manifold of the semi-parabolic fixed point $\fq_{a}$ as $a\rightarrow 0$. Consider the \He map $H$ and its inverse $H^{-1}$ written as in Equation \ref{eq:Hw}, respectively \ref{eq:Hw-1}.  

Fix $\lambda=e^{2\pi i p/q}$. Suppose $H$ has a semi-parabolic fixed point at $\fq_{a}$ such that $DH(\fq_{a})$ has eigenvalues $\lambda$ and $\mu$, with $|\mu|<1$.  We have $\lambda\mu=-a^{2}$ so $\mu = -\frac{a^{2}}{\lambda}$ and $|\mu|=|a|^{2}$. Set for simplicity
$
	q_{a}:= \frac{\lambda}{2}-\frac{a^{2}}{2\lambda}=\frac{\lambda+\mu}{2}.
$
With this notation, the equation \ref{eq:fqa} of the fixed point $\fq_{a}$ reduces to
\[
\fq_{a} := \pvec{q_{a}}{aq_{a}}=\frac{\lambda+\mu}{2}\pvec{1}{a}.
\]
Let $v = \pvec{-a/\lambda}{1}=\pvec{\mu/a}{1}$ be an eigenvector for the eigenvalue $\mu$.
The semi-parabolic fixed point $\fq_{a}$ has a stable manifold $W^{s}(\fq_{a})\subset \C^{2}$ as in Equation \ref{eq:Ws-qa}. The stable manifold is biholomorphic to $\C$ and has a natural parametrization given by the following proposition.
\begin{prop}\label{prop:Ws1}
The stable manifold $W^{s}(\fq_{a})$ has a parametrization $F_{a}:\C \rightarrow W^{s}(\fq_{a})$ given by
\begin{equation}\label{eq:funFa}
F_{a}(z) = \lim\limits_{m\rightarrow \infty} H^{-\circ m}(\fq_{a}+\mu^{m}vz).
\end{equation}
 $F_{a}$ is an injective immersion of $\C$ onto $W^{s}(\fq_{a})$ with the property that $F_{a}(\mu z) = H(F_{a}(z))$.
\end{prop}
\proof
The proof is similar to the proof of Theorem 1 from \cite{H1}. Consider the inverse map $H^{-1}$ instead of $H$. Then $\fq_{a}$ is a fixed point of $H^{-1}$ and $DH^{-1}(\fq_{a})$ has eigenvalues $\overline{\lambda}$ and $\mu'=1/\mu$, where $|\mu'|>1$. The fixed point $\fq_{a}$ has now an unstable manifold $W^{u}(\fq_{a})$ which has a natural parametrization given by $F_{a}$ as shown in \cite{H1}.
\qed

\begin{prop}\label{prop:Ws} The parametrizing function $F_{a}\rightarrow F_{0}$ as $a\rightarrow 0$, where
\[ F_{0}(z) := \fq_{0}+\hvec{0}{z} = \hvec{\lambda/2}{z}.
\]
\end{prop}
\proof Define a sequence of points
\[
    \hvec{x_{i}}{y_{i}}=H^{-1}\hvec{x_{i-1}}{y_{i-1}} \ \ \ \mbox{for}\ 1\leq i\leq m,
\]
where
\[
    \hvec{x_{0}}{y_{0}} = \fq_{a}+\mu^{m}vz = \hvec{q_{a}+\frac{\mu^{m+1}}{a}z}{aq_{a}+\mu^{m}z}.
\]
At the first step we have
\[
    \hvec{x_{1}}{y_{1}} = H^{-1}\hvec{x_{0}}{y_{0}} =\frac{1}{a} \hvec{y_{0}}{x_{0}-p(y_{0}/a)-a^{2}w},
\]
so
\[
	x_{1}=\frac{y_{0}}{a}=q_{a}+\frac{\mu^{m}}{a}z.
\]
From the fixed point equation $H(\fq_{a})=\fq_{a}$ we get that $p(q_{a}) + a^{2} q_{a}+a^{2}w=q_{a}$ so $q_{a}-p(q_{a})-a^{2}w=a^{2} q_{a}$. Moreover, the matrix $DH(\fq_{a})$ has eigenvalues $\lambda$ and $\mu$ so $\lambda + \mu =tr(DH(\fq_{a}))$, which gives $p'(q_{a})=\lambda + \mu$. 
Since $y_{0}=aq_{a}+\mu^{m}z$ and $\mu = -\frac{a^{2}}{\lambda}$, we can write
\[
p(y_{0}/a)=p\left(q_{a}+\frac{\mu^{m}}{a}z\right) = p(q_{a})+p'(q_{a})\frac{\mu^{m}}{a}z+\bigO(\mu^{2m-1}).
\]
Note that this is a finite sum. Thus the equation for $y_{1}$ has the following form
\begin{eqnarray*}
    y_{1}=\frac{x_{0}-p(y_{0}/a)-a^{2}w}{a} &=& \frac{q_{a}-p(q_{a})-a^{2}w}{a} + \frac{ \mu^{m+1}z-p'(q_{a})\mu^{m}z}{a^{2}}+\bigO(a\mu^{2m-2})\\
    &=&aq_{a} + \mu^{m-1}z + \bigO(a\mu^{2m-2}).
\end{eqnarray*}
By induction we can show that for $1\leq i\leq m$ we have
\begin{eqnarray*}
    x_{i} &=& q_{a}+\frac{\mu^{m-(i-1)}}{a}z+\bigO(\mu^{2(m-(i-1))})\\
    y_{i} &=& aq_{a}+\mu^{m-i}z+\bigO(a\mu^{2(m-i)}).
\end{eqnarray*}
For $i=m$ these reduce to
\begin{eqnarray*}
    x_{m} &=& q_{a} -a\overline{\lambda}z+\bigO(a^{4})\\
    y_{m} &=& aq_{a}+z+\bigO(a).
\end{eqnarray*}
Thus $x_{m}\rightarrow q_{0}$ and $y_{m}\rightarrow z$ as $a\rightarrow 0$.  
Therefore $F_{a}$ converges to $F_{0}$ uniformly on compact subsets of $\C$ as $a\rightarrow 0$. 
\qed

\section{Choosing uniform normalizing coordinates}\label{subsec: U-normalizing}

In Section \ref{sec: NormalForm}, we described the normal form of germs of transformations of $\C^{2}$ with a semi-parabolic fixed point. In this section, we will study the normal form for the family of \He maps with a semi-parabolic fixed point, which are small perturbations of the parabolic polynomial $p(x)$ inside the parabola $\mathcal{P}_{\lambda}$. Since $\mathcal{P_{\lambda}}$ is parametrized by $a$, we write the \He map as $H_{a}(x,y)=\hvec{x^{2}+c+ay}{ax}$, where $c$ is chosen as in Equation \ref{eq:P-lambda}, so that $(c,a)\in \mathcal{P_{\lambda}}$.

We will show how to extend the results from Section \ref{sec: NormalForm} in order to get uniform bounds (with respect to the parameter $a$) on the size of the normalizing neighborhood. We will prove that the coordinate transformation $\phi_{a}$ that puts the \He map $H_{a}$ in the normal form is holomorphic with respect to $a$. Then we will use the theory developed in Section \ref{sec: Petals} to exhibit local attractive and repelling sectors for the \He map $H_{a}$. The attractive sectors will belong to the interior of $K^{+}$. In the repelling sectors we will show that the derivative of the \He map is weakly expanding in the ``horizontal'' direction and strongly contracting in the ``vertical'' direction.  

We first look at the normal form from \cite{DH} and \cite{H} for the polynomial $p(x)=x^{2}+c_{0}$
which has a parabolic fixed point $q_{0}=\frac{\lambda}{2}$, of multiplier $\lambda=e^{2\pi i p/q}$.
Denote, for the clarity of exposition, $\rho :=(\sqrt{2}R)^{-1/q}$ in the definition of the set $\Delta^{\pm}$ from Equations \ref{eq:W+} and \ref{eq:W-}.

\begin{lemma}\label{lemma:Upolynom}
There exists a neighborhood $V_{0}$ of $q_{0}$ and an isomorphism $\phi: V_{0}\rightarrow \D_{\rho}$ such that $\widetilde{p}(x)=\phi \circ p\circ \phi^{-1} (x)$ where $\widetilde{p}(x)=\lambda x\left(1+x^{q}+C x^{2 q}+\bigO(x^{2q+1})\right)$.
Furthermore, there exists $\rho$ small enough such that in the region
\[
\Delta^{-}= \left\{ |x|<\rho\ |\ Re(x^{q})> \epsilon_{0}|Im(x^{q})| \right\}
\]
the map $\widetilde{p}$ satisfies $|\widetilde{p}\,'(x)|>1+\epsilon_{1}|x|^{q}$. The compact region
\[
\Delta^{+}= \left\{|x|<\rho\ |\ Re(x^{q})\leq \epsilon_{0}|Im(x^{q})|\right\}\
\]
satisfies $p(\Delta^{+})\subset  int(K_{p})\cup \{0\}$.
\end{lemma}
\proof After a global coordinate change that brings the parabolic fixed point at the origin, we can  write the polynomial as $p(x)=\lambda x + x^{2}$. Since $p$ is a quadratic polynomial, the fixed point $q_{0}$ can only have parabolic multiplicity $1$, hence its multiplicity as a solution of the equation $p^{\circ q}(z)-z=0$ is $q+1$. 
The local normal form of $p$ around $q_{0}$ is obtained by successive elimination of the terms of degree less than $2q+1$ which are not congruent to $1$ mod $q$, using the same coordinate transformations as in Theorem \ref{thm:normalform3}.

The derivative $\widetilde{p}\, '(x)$ is weakly expanding in  $\Delta^{-}$. To show this, let $m$ be chosen so that  $\big{|}\widetilde{p}\, '(x)-\lambda(1+(q+1)x^{q})\big{|}<m|x|^{2q}$ on $\D_{\rho}$.
Since $|\lambda|=1$ and $Re(x^{q})>\epsilon_{1}|x|^{q}$ from Equation \ref{eq:eps1}, we can estimate $|\widetilde{p}\, '(x)|$ on $\Delta^{-}$ as follows:
\begin{eqnarray*}
|\widetilde{p}\, '(x)|&=&\big{|}1+(q +1)x^{ q}+\bigO(x^{2 q})\big{|}\geq\big{|}1+(q +1)x^{ q}\big{|}-m|x|^{2q}\\
&\geq& 1 + (q+1)\epsilon_{1}|x|^{q} - m|x|^{2q} > 1 + \epsilon_{1}|x|^{q},
\end{eqnarray*}
for $x$  small enough so that $|x|^{q}<q\epsilon_{1}/m$.
It follows that $|\widetilde{p}\, '(x)|>1+\epsilon_{1}|x|^{q}$ for $x\in \Delta^{-}$ and $|x|$ sufficiently small.
\qed

Choose $\rho'>0$ such that the disk $\D_{2\rho'}(q_{0})$ of radius $2\rho'$ centered at $q_{0}$ is contained in the neighborhood $V_{0}$. We make this choice for technical reasons.

\begin{thm}\label{thm:UniformNF}
Let $r>\rD$ be a fixed constant.
There exists $\delta>0$ such that for any $(c,a)\in \mathcal{P}_{\lambda}$ with $|a|<\delta$  we can find a coordinate transformation $\phi_{a}$ from a tubular neighborhood $B=\D_{\rho'}(q_{0})\times\D_{r}$
of the local stable manifold of the semi-parabolic fixed point $\fq_{a}$
\[
\phi_{a}:B\rightarrow \D_{\rho}\times \D_{r+\bigO(|a|)}
\]
in which the \He map has the form $\widetilde{H}_{a}(x,y)=(x_{1},y_{1})$, with
\begin{equation}\label{eq:UniformNF3}
\left\{\begin{array}{l}
    x_{1}= \lambda (x+x^{q+1} + Cx^{2q+1}+ a_{2q+2}(y)x^{2q +2}+\ldots )\\
    y_{1}= \mu y + xh(x,y)
\end{array}\right.
\end{equation}
and $C$ is a constant (depending on $a$) and $xh(x,y)=\bigO(a)$. Moreover,
the maps $\phi_{a}$ are holomorphic in $a$ and
\[
\lim\limits_{a\rightarrow 0}\phi_{a}=\phi_{0}(x,y)=(\phi(x),y),
\] where
 $\phi:\D_{\rho'}\rightarrow
\D_{\rho}$
  is the change of coordinates for the polynomial $p(x)=x^2+c_0$ with a parabolic fixed point at $\displaystyle q_{0}$,
  \[
  \phi \circ p\circ \phi^{-1} (x)=\lambda x(1+x^{q}+C x^{2q}+\bigO(x^{2q+1})).
  \]
\end{thm}
\proof We will follow the same steps as in Section \ref{sec: NormalForm}. The following two propositions are part of the proof.

The degenerate map $H_{0}(x,y)=(p(x),0)$ has a semi-parabolic fixed point $\fq_{0}=(\frac{\lambda}{2},0)$ of multiplicity $q+1$ and the stable manifold $W^{s}(\fq_{0})$ is just a vertical line passing through $\fq_{0}$. The multiplicity of the semi-parabolic fixed point is constant in a neighborhood of $a=0$ in $\mathcal{P}_{\lambda}$. When $a\neq 0$, $W^{s}(\fq_{a})$ is an analytic submanifold biholomorhic to $\C$. By \cite{H1}, $W^{s}(\fq_{a})$ depends analytically on $a$ in a neighborhood of $a=0$ inside $\mathcal{P}_{\lambda}$.
\begin{defn}
Denote by $S_{r}$ the horizontal strip $S_{r}:=\{(x,y)\in \C^{2}\ |\ |y|<r\}$ and by $W^{s}_{loc}(\fq_{a})$ the connected component of $W^{s}(\fq_{a})\cap S_{r}$ that contains the fixed point $\fq_{a}$.
\end{defn}

Let us choose $\delta>0$, such that for all $(c,a)\in \mathcal{P}_{\lambda}$ with $|a|<\delta$  the \He map $H_{a}$ has a fixed point $\fq_{a}$ of multiplicity $q+1$ and such that the local stable manifold $W^{s}_{loc}(\fq_{a})$ is ''vertical-like''. Rigorously, we require that the horizontal distance between $W^{s}_{loc}(\fq_{a})$ and the vertical line that contains $\fq_{a}$ is less than $\frac{\rho'}{4}$, and that $W^{s}_{loc}(\fq_{a})$ has no horizontal foldings.

The parametrizing map $F_{a}:\C\rightarrow W^{s}(\fq_{a})$ defined in \ref{eq:funFa} is analytic in the parameter $a$. By Proposition \ref{prop:Ws}, 
it degenerates to a translation in the horizontal direction when $a=0$, given by 
$F_{0}(y)=\fq_{0}+(0,y)$.
By Proposition \ref{prop:Ws} we know that
\[
F_{a}(y)=F_{0}(y)+\bigO(a),
\]
so $F_{a}$ will map the disk $\{y\in \C\ |\ |y|<r\}$ onto a holomorphic disk inside $W^{s}(\fq_{a})$ around $\fq_{a}$ of size approximately $r+\bigO(a)$. For $a$ small, fix therefore $3<r'\leq r$ such that $W^{s}_{loc}(\fq_{a})\cap S_{r'}\subset F_{a}(S_{r})$. In principle $r'=r+\bigO(|a|)$, but  the vertical size is not a delicate issue, so we can think of $r'$ as simply being $r$.

\begin{prop}\label{thm:UniformNF1}
Choose $\delta>0$ as before. For all $(c,a)\in \mathcal{P}_{\lambda}$ with $|a|<\delta$ there exists a coordinate transformation $\phi_{a}^{1}:S_{r'}\rightarrow S_{r}$, such that in the new coordinates, the \He map $H_{a}$ has the form $H_{a}(x,y)=(x_{1},y_{1})$, with
\begin{equation}\label{eq:UniformNF1}
\left\{\begin{array}{l}
    x_{1}= a_{1}(y)x+a_{2}(y)x^{2}+\ldots \\
    y_{1}=\mu y + xh(x,y)
\end{array}\right.,
\end{equation}
where $a_{j}(\cdot)$ and $h(\cdot,\cdot)$ are holomorphic functions from $\{y\in\C, |y|<r'\}$ to $\C$, respectively from $\{(x,y)\in\C^{2}, |y|<r'\}$ to $\C$, with $a_{1}(0)=\lambda$ and $h(0,0)=0$.
\end{prop}
\proof
Suppose $F_{a}(y)=(f(y), g(y))$, and let $\psi_{a}:S_{r}\rightarrow \C^2$ be the map
\[
\psi_{a}(x,y)=(x+f(y), g(y)).
\]
It is easy to see that $\psi_{a}$ is an invertible function. The Jacobian matrix is given by
\[
D\psi_{a}|_{(x,y)}=\left(\begin{array}{cc}1 & f'(y) \\0 & g'(y)\end{array}\right).
\]
The local stable manifold $W^{s}_{loc}(\fq_{a})$ is vertical-like. In particular it has no horizontal foldings, hence $g'(y)\neq 0$ for $|y|<r$. This means that $\psi_{a}$ is invertible in the strip $S_{r}$. Define $\phi_{a}^{1}(x,y):=\psi_{a}^{-1}(x,y)$.

The fact that $\phi_{a}^{1}(x,y)$ is holomorphic in $a$ follows immediately, since we know that $F_{a}(y)$ depends holomorphically on $a$. From Proposition \ref{prop:Ws} we obtain that
\[
\phi_{a}^{1}(x,y)=\phi_{0}^{1}(x,y)+\bigO(a).
\]
The transformation $\phi_{0}^{1}$ is straightforward to compute
\[\phi_{0}^{1}(x,y)=(x,y)-\fq_{0}=\left(x-\lambda/2,y\right).
\]
In the new coordinate system the \He map $H_{0}$ becomes $H_{0}(x,y)=(x_{1},y_{1})$, where
\begin{equation*}
\left\{\begin{array}{l}
    x_{1}= \lambda x+x^{2} \\
    y_{1}= 0
\end{array}\right.
\end{equation*}
Therefore, when $a\neq 0$, it is easy to control the size of the coefficients in Equation \ref{eq:UniformNF1} in terms of $a$, as follows:
$a_{1}(y)=\lambda+\bigO(a)$, $a_{2}(y)=1+\bigO(a)$, $a_{i}(y)=\bigO(a)$ for $i>2$ and $h(x,y)=\bigO(a)$.
\qed

\begin{prop}\label{thm:UniformNF2}
There exists $\delta>0$ such that for all parameters $(c,a)\in \mathcal{P}_{\lambda}$ with $|a|<\delta$ there exists a coordinate transformation $\phi^{2}_{a}:\D_{1/2}\times\D_{r}\rightarrow S_{r}$ in which $H_{a}$ has the form $H_{a}(x,y)=(x_{1},y_{1})$, with
\begin{equation}\label{eq:UniformNF2}
\left\{\begin{array}{l}
    x_{1}= \lambda x+a_{2}x^{2}+\ldots + a_{2q+1}x^{2q+1} + a_{2q+2}(y)x^{2q+2}+\ldots \\
    y_{1}= \mu y + xh(x,y)
\end{array}\right.
\end{equation}
where $a_{2}$ is close to $1$ and the coefficients $a_{3},\ldots, a_{2q+1}$ are constants close to $0$.
\end{prop}
\proof
Suppose $H_{a}$ is written as in Equation \ref{eq:UniformNF1}. The proof of this proposition is the same as that of Theorem \ref{thm:normalform2} with $m=2q+1$. Notice that $a_{1}(y)=\lambda+\bigO(a)$ for $|y|<r$, so one can perform the same change of coordinates as in the proof of Theorem \ref{thm:normalform2}
\[
T_{1}: (x,y)\rightarrow (u(y)x,y),\ \mbox{where}\ u(y)=\prod\limits_{n\geq 0} a_{1}(\mu^{n}y)
\]
in order to set $a_{1}(y)=\lambda$.  Since $a_{1}(y)$ is close to $\lambda$ when $|y|<r$, it follows that the product is convergent when $|y|<r$. We get that $u(y)\neq 0$, hence $T_{1}(x,y)$ is invertible.

The coordinate changes that make $a_{j}(y)$ constant for $2\leq j \leq 2q+1$ are of the form
\[
T_{j}:(x,y)\rightarrow (x+v(y)x^{j},y),\ \mbox{where}\ v(y)=\sum\limits_{n\geq 0} a_{j}(\mu^{n}y)-a_{j}(0).
\]
Clearly the sum is convergent when $|y|<r$. We get $v(y)=\bigO(a)$. The transformation $T_{j}$ is invertible because $x$ is bounded ($1/2$ would be a reasonable bound for $x$, but any bound less than 1 would do), so for $a$ small $1+v(y)jx^{j-1}$ does not vanish.

The coordinate changes that are done in order to make the first $2q+1$ coefficients constants are identity on the second coordinate. Denote by $\phi^{2}_{a}(x,y)$ their composition.
Notice also that in Equation \ref{eq:UniformNF1}, $H_{0}(x,y)=(\lambda x+x^{2},0)$ already has constant coefficients, so $\phi^{2}_{0}$ is just the identity map. It is easy to check that 
\[
\phi^{2}_{a}(x,y)=(x+\bigO(a),y)
\]
and $h(x,y)=\bigO(a)$. We also have that $a_{2}=1+\bigO(a)$, $a_{i}=\bigO(a)$ for $2<i\leq 2q+1$ and $a_{i}(y)=\bigO(a)$ for $i>2q+1$.
\qed

We are now able to finish the proof of Theorem \ref{thm:UniformNF}. The coordinate changes done in Proposition \ref{thm:UniformNF1} did not require any bounds on $x$. The coordinate transformations done in Proposition \ref{thm:UniformNF2} required only a mild assumption on $x$ (such as $|x|<1/2$). We will now use the coordinate changes for the polynomial $p$ to put the \He map in the normal form given in Equation \ref{eq:UniformNF3}. We will thus require a bound on $x$ comparable to the size of the normalizing neighborhood $V_{0}$ from Lemma  \ref{lemma:Upolynom}.

Assume that $H_{a}$ is written in the form \ref{eq:UniformNF2}.
We use the same transformations as in Theorem \ref{thm:normalform3} in order to eliminate the terms $x^{i}$, where $1<i<q+1$ and $q+1<i<2q+1$. Let $\phi^{3}_{a}:\D_{2\rho'-\bigO(|a|)}\times \D_{r}\rightarrow \D_{\rho}\times \D_{r}$ denote the coordinate change from Theorem \ref{thm:normalform3}.
When $a=0$, $H_{0}(x,y)=(\lambda x+x^{2},0)$ and
\[
\phi^{3}_{0}(x,y)=(\phi(x),y),
\]
where $\phi(x)$ is the coordinate transformation used in Lemma \ref{lemma:Upolynom} to put $p(x)=\lambda x+x^{2}$ in the normal form $\widetilde{p}(x)= \lambda (x+x^{q+1} + Cx^{2q+1}+\bigO(x^{2q+2}))$.

Define $\phi_{a}(x,y)=\phi_{a}^{3}\circ \phi_{a}^{2}\circ \phi_{a}^{1}(x,y)$.
Recall that $\phi_{0}^{1}$ is a horizontal translation by $\lambda/2$ and $ \phi_{0}^{2}$ is the identity map, so when $a=0$ the composition of the three transformations yields exactly the coordinate transformation used in Lemma \ref{lemma:Upolynom} to put $p(x)=x^{2}+\frac{\lambda}{2}-\frac{\lambda^{2}}{4}$ in the normal form $\widetilde{p}(x)= \lambda (x+x^{q+1} + Cx^{2q+1}+\bigO(x^{2q+2}))$.
\qedof \textbf{of Theorem \ref{thm:UniformNF}}

\medskip
In the normalizing coordinates we define attractive and repelling sectors for the \He map and study the behavior of $\widetilde{H}_{a}$ for $|a|<\delta$.

\begin{lemma}[\textbf{Attractive/Repelling sectors}]\label{lemma:HenonSectors}
Let $W^{\pm}$ be defined as in Equations \ref{eq:W-} and \ref{eq:W+}. There exists $\rho>0$ and $\delta>0$ such that for all $|a|<\delta$ the derivative $D\widetilde{H}_{a}$ expands horizontally by a factor of $(1+\frac{\epsilon_{1}}{2}|x|^{q})$ in the region
\[
W^{-}=\Delta^{-}\times\D_{r}=\left\{ |x| \leq \rho\ |\   \mbox{Re}(x^{q})> \epsilon_{0}|Im(x^{q})| \right\} \times \D_{r}
\]
The compact region
\[
W^{+}=\Delta^{+}\times\D_{r}=\left\{ |x| \leq \rho\ |\  \mbox{Re}(x^{q})\leq  \epsilon_{0}|Im(x^{q})| \right\} \times \D_{r}
\]
satisfies $\widetilde{H}_{a}(W^{+})\subset int(K^{+})\cup \{0\}\times \D_{r}$.
\end{lemma}
\proof 
By construction, $W^{+}\subset\Pet_{att}\cup\{0\}\times \D_{r}$ and all points in $\Pet_{att}$ are attracted to the origin under forward iterations by Proposition \ref{prop:sectors}. Hence $W^{+}\subset int(K^{+})\cup\{0\}\times \D_{r}$. The horizontal expansion in $W^{-}$ follows from Proposition \ref{prop:hcone-in-petal} below. 
\qed

The multiplicity of the semi-parabolic fixed point is $q+1$, so there are exactly $q$ connected components of $W^{-}$ and $q$ components of $int(W^{+})$.

\begin{figure}[htb]
\begin{center}
\includegraphics[scale=0.95, bb = 250 480 370 735]{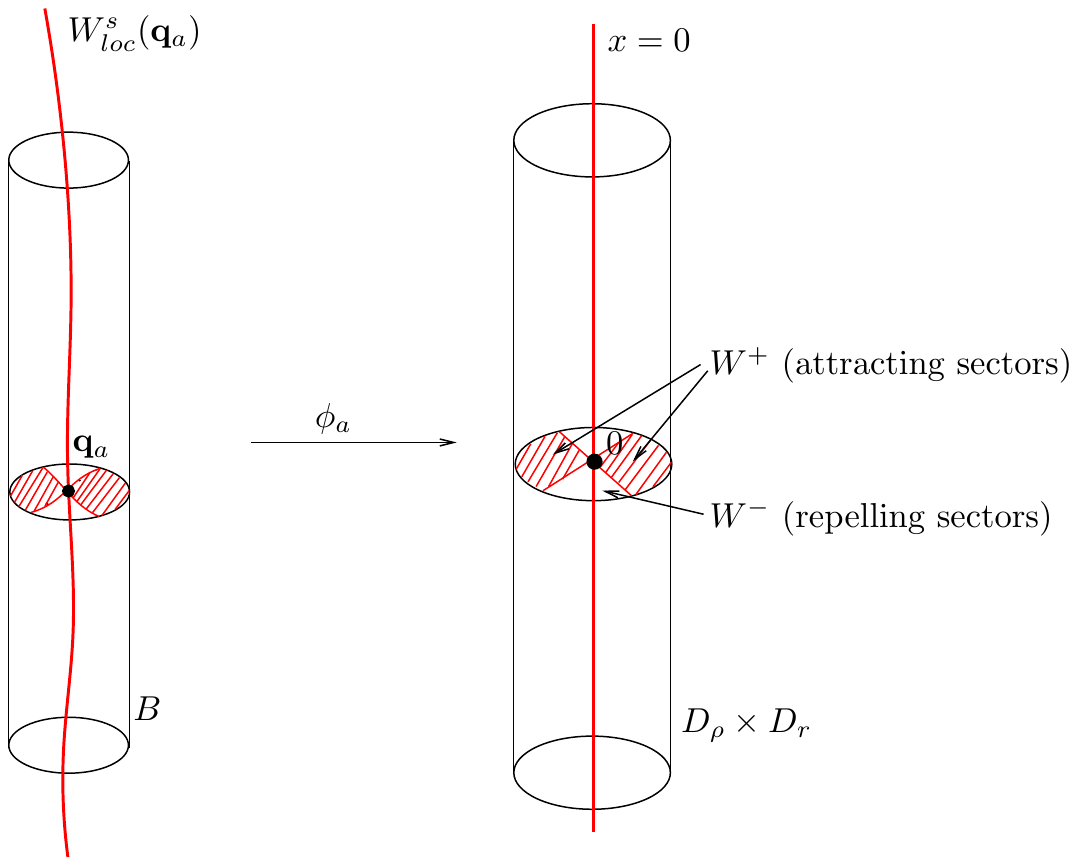}
\end{center}
\caption{The transformation $\phi_{a}$ and the sectors $W^{\pm}$ for $q=2$.}
\end{figure}

Let $|a|<\delta$ as before and consider the \He map $\widetilde{H}:\D_{\rho}\times \D_{r}\rightarrow \C^{2}$ written in normal coordinates as in Theorem \ref{thm:UniformNF},
\[
\widetilde{H}_{a}\hvec{x}{y} = \hvec{\lambda
(x+x^{q+1}+g_{a}(x,y))}{\mu y + xh_{a}(x,y)},
\] 
where  $g_{a}(x,y)=g_{0}(x)+\bigO(a)$ and $h_{a}(x,y)=\bigO(a)$ and 
\begin{eqnarray*}
    g_{a}(x,y)&=& C_ax^{2q+1}+ a_{2q+2}(y)x^{2q +2}+\ldots \\
    h_{a}(x,y)&=& b_{1} (y) +\ldots + b_{k}(y)x^{k}+\ldots.
\end{eqnarray*}
When $a=0$, $\widetilde{H}_{0}(x,y)=(\widetilde{p}(x),0)$, where $\widetilde{p}(x)=\lambda(x+x^{q+1}+g(x))$ and 
\begin{eqnarray*}
g(x) &=& C_0x^{2q+1}+a_{2q+2}x^{2q+2}+\ldots.
\end{eqnarray*}
\noindent The function $g_{0}(x,y)=g(x)$ is just a function
of the variable $x$, hence $\partial_{y} g_{0}(x,y)\equiv 0$. For $|a|<\delta$ we
can assume that there exists a constant $M_{a}$ with $0<M_{a}<1$
such that 
\begin{equation}\label{eq:Ma}
\big{|}\partial_{y} g_{a}(x,y)\big{|}<M_{a}|x|^{2q+2}.
\end{equation}
As usual, $\partial_{x}$ and $\partial_{y}$ denote the partial derivatives with respect to the variable $x$, and respectively $y$.  
When $a=0$ we also know that $xh_0(x,y)\equiv 0$.
Moreover by the construction of the normalizing coordinates we have
$xh_a(x,y)=\bigO(a)$. There exists a constant $N_{a}$, depending on $a$,  with $0<N_a<1$ such that when $|a|<\delta$ the following bounds hold
\begin{eqnarray}\label{eq:Na}
\big{|}\partial_{x} (xh_a)(x,y)\big{|} < N_a \ \
\mbox{and}\ \ \big{|}\partial_{y} (xh_a)(x,y)\big{|} < N_a.
\end{eqnarray}
Let $\partial_{x} g_a(x,y)=x^{2q}t_a(x,y)$ and denote by $m$ the supremum of $|t_a(x,y)|$ on the set $W^{-}$, where the supremum is taken after all $|a|<\delta$. Hence for any $a$ with $|a|<\delta$ and any $(x,y)$ taken from the repelling sectors $W^{-}=\Delta^{-}\times \D_{r}$ of the \He map we have
\[
\big{|}\partial_{x} g_a(x,y)\big{|}< m|x|^{2q}.
\]
By eventually reducing
$\rho>0$, we can assume as in Equation \ref{eq:eps1} that
\begin{equation}\label{eq:eps2}
|1+(q+1)x^{q}|-m|x|^{2q}>1+\epsilon_{1}|x|^{q}>1,\ \mbox{for all}\ x\in \Delta^{-}.
\end{equation}

\begin{defn}\label{def:hcone-in-petal}
Let $(x,y)$ be a point in the repelling sectors $W^{-}$ of the \He
map. Define the horizontal cone at $(x,y)$ to be
\begin{equation*}
\mathcal{C}^{h}_{(x,y)}=\{(\xi,\eta)\in T_{(x,y)} W^{-},\ |\xi|>|\eta|\}.
\end{equation*}
\end{defn}
We will show that the horizontal cones are invariant
under $D\widetilde{H}$ and that $D\widetilde{H}$ is expanding inside the horizontal cones.

\begin{prop}[\textbf{Horizontal cones}]\label{prop:hcone-in-petal}
Let $(x,y)$ and $(x',y')$ be two points from
$W^{-}$ such that $\widetilde{H}(x,y)=(x',y')$. Then
\[
D\widetilde{H}_{(x,y)}\left(\mathcal{C}^{h}_{(x,y)}\right)\subset Int\
\mathcal{C}^{h}_{(x',y')}
\]
and $ \big{\|} D\widetilde{H}_{(x,y)}(\xi,\eta)\big{\|} \geq (1+\frac{\epsilon_{1}}{2}|x|^{q})
\|(\xi,\eta)\| \ $ for $\ (\xi,\eta)\in \mathcal{C}^{h}_{(x,y)}$.
\end{prop}
\proof
Pick $(\xi, \eta)\in \mathcal{C}^{h}_{(x,y)}$ 
and let
$D\widetilde{H}_{(x,y)}(\xi,\eta)=(\xi',\eta')$. The derivative of $\widetilde{H}$ is
\[
D\widetilde{H}_{(x,y)}=\left(
                 \begin{array}{cc}
                   \lambda(1+(q+1)x^q+\partial_{x} g_a(x,y)) & \lambda \partial_{y}g_a(x,y) \\
                   \partial_{x} (xh_a)(x,y) & \mu+\partial_{y} (xh_a)(x,y) \\
                 \end{array}
               \right).
\]
Consider now the Euclidean metric on the set
$\D_{\rho}\times \D_{r}$ and estimate
\begin{eqnarray}\label{eq:Ho1}
|\eta'|&\leq&  N_{a}|\xi|+\left(|\mu|+N_{a}\right)|\eta| < (2N_{a}+|\mu|)|\xi| \\
|\xi'| &\geq& \left(|1+(q+1)x^q|-m|x|^{2q}\right)|\xi|-
M_{a}|x|^{2q+2}|\eta| \nonumber \\
&>&   \left(|1+(q+1)x^q|-m|x|^{2q}-M_{a}|x|^{2q+2}\right)|\xi| \label{eq:ineq}.
\end{eqnarray}
We then obtain
\[
|\eta'|<\frac{B_{2}}{B_{1}}|\xi'|,
\]
where $B_{1}$ and $B_{2}$ are defined in the obvious way
\begin{eqnarray*}
B_{2} &:=& 2N_{a}+|\mu| \\
B_{1} &:=&
|1+(q+1)x^q|-m|x|^{2q}-M_{a}|x|^{2q+2}.
\end{eqnarray*}
The bounds $N_{a}$, $M_{a}$ and $|\mu|=|a|^{2}$ tend to $0$ as
$a\rightarrow0$, so we can make $B_{2}$ as small as we want, for example we assume $B_{2}<\frac{1}{2}$. The
points $(x,y)$ and $(x',y')$ are chosen from the repelling sectors, so we can assume that
\begin{equation}\label{eq: expansion-h}
B_{1}\geq 1+\frac{\epsilon_1}{2}|x|^q.
\end{equation}
In conclusion we get $|\eta'|<B_{2}|\xi'|$, so
$(\xi',\eta')\in Int\ \mathcal{C}^{h}_{(x',y')}$.
In fact, when $\eta=0$, we have that $|\eta'|<N_{a}|\xi'|$, which will be useful in Lemma \ref{lemma:Step2}.

The same computation \ref{eq:ineq} also shows that $|\xi'|>B_{1}|\xi|$ so $D\widetilde{H}$ expands the horizontal length of vectors, i.e.
\begin{equation} \label{prop: expansion-hcones-in-petal}
\|(\xi',\eta')\|=\max\left\{|\xi'|,|\eta'|\right\}=|\xi'| > B_{1}|\xi| 
\geq |\xi|=\max\left\{|\xi|,|\eta|\right\}=\|(\xi,\eta)\|. 
\end{equation}
\qed

\section{Construction of a neighborhood $V$ for $J^+$}\label{subsec: V-parabolic}
We will build a neighborhood of $J^{+}$ for a semi-parabolic \He map $H_{a}$ inside a polydisk $\D_{r}\times\D_{r}$, inspired by the construction of a neighborhood of the Julia set of a parabolic polynomial $p$ on which $p$ is strictly, but not strongly expanding, as in \cite{DH}. Inside a tubular neighborhood $B$ of the local stable manifold $W^{s}(\fq_{a})$, we want to forget about the dynamics of the polynomial $p$ and construct a neighborhood of $J^{+}\cap B$ that is meaningful for the dynamics of the  \He map.

Let $q_{0}$ be the parabolic fixed point of the polynomial $p$ and let $q_{1}$ be the other preimage of $q_{0}$ under $p$. Suppose $|a|<\delta$ and consider $B=\D_{\rho'}(q_{0})\times \D_{r}$ as defined in Theorem \ref{thm:UniformNF}. Let $W^{\pm}$ be the attractive/repelling sectors from Lemma \ref{lemma:HenonSectors}. Define $W^{\pm}_{B}=\phi_{a}^{-1}(W^{\pm})\cap\D_{r}\times\D_{r}$. For our purpose, it is more convenient to view these sectors in $B$ rather than in the normalized coordinates. 
The set $H^{-1}(B)\cap \D_{r}\times\D_{r}$ consists of two connected components. We denote by $B'$ the connected component which contains $q_{1}$. 
Define $W^{\pm}_{B'}$ to be the preimage of the attractive/repelling sectors $W^{\pm}_{B}$ in $B'$.

\begin{figure}[htb]
\begin{center}
\includegraphics[scale=1.53, bb = 230 620 370 730]{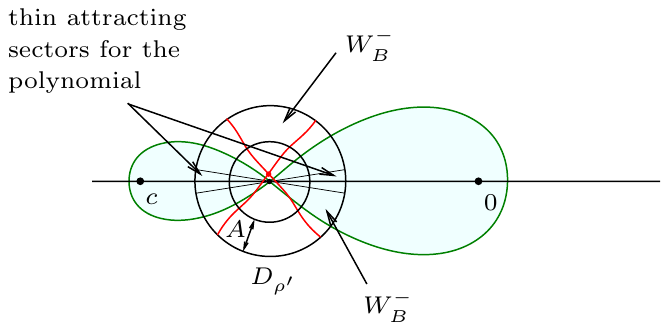}
\end{center}
\caption{Here $q=2$. This is a cross section around the parabolic fixed point of the polynomial $p(x)=x^{2}+c_{0}$. The red lines are the boundaries of the attractive sectors for the \He map. The thin attractive sectors for the polynomial and their preimages are shown in green. }
\end{figure}

Choose $\rho''>0$ as large as possible so that $p^{\circ 2}(\D_{\rho''}(q_{0}))\subset\D_{\rho'}(q_{0})$. Clearly this choice depends only on the parabolic polynomial $p$ and the radius $\rho'$. Consider the annulus $A:=A(q_{0}; \rho', \rho'')$  between the disk of radius $\rho''$ and the disk of radius $\rho'$ centered at $q_{0}$. Let $n$ be the first iterate of $p$ such that $p^{\circ(n+1)}(0)\in \D_{\rho'}(q_{0})$ and implicitly $p^{\circ(n+1)}(0)\in A$. We now construct  attractive sectors $S_{att}$ associated with the parabolic polynomial $p$ in $\D_{\rho'}(q_{0})$, thin enough along the attractive axes of the polynomial so that
\begin{equation}\label{eq:Delta-}
\left( p^{-\circ(n+1)}(S_{att})\cap A\right)\times \D_{r}\subset W^{+}_{B}.
\end{equation}

Let $\partial^{in}(W^{+}_{B})$ be the part of the boundary of $W^{+}_{B}$ that lies strictly inside $A\times \D_{r}$.  Similarly, let $\partial^{in}(p^{-\circ (n+1)}(S_{att})\cap A)$ be the part of the boundary of $p^{-\circ(n+1)}(S_{att})\cap A$ that is strictly inside the annulus $A$. We will further require that $S_{att}$  be thin enough so that the distance between the two boundaries $\partial^{in}(W^{+}_{B})$ and $\partial^{in}(p^{-\circ(n+1)}(S_{att})\cap A)\times \D_{r}$ is at least $\eta_{0}>0$. The constant $\eta_{0}$ depends only on the local dynamics of the polynomial $p$ and it can be taken to be a fraction of the distance between $\partial^{in}(\Delta^{+}\cap A)$ and an attractive axes that passes through $\Delta^{+}$, where $\Delta^{+}$ is defined in Lemma \ref{lemma:Upolynom}.

Let $\Omega= p^{-\circ n}(S_{att})$. In the definition of the set $\Omega$, we only consider the preimages of $S_{att}$ that contain the parabolic fixed point in the boundary, so they are local preimages. The set $\Omega$ has $q$ connected components and contains the critical value, but not the critical point of the polynomial $p$. Let us now define a set $U$ as the complement of $\Omega$ inside an equipotential of the Green's function of $p$, i.e.
\begin{equation}\label{eq:nbdU}
U: = \C - \Omega-\{ z\in \C-K_{p}\ |\ |\Phi_{p}^{-1}(z)|\geq R\}
\end{equation}
for some large enough $R>2$. Define $U':=p^{-1}(U)$ and $\Omega':=p^{-1}(\Omega)$. The constant $R$ is chosen so that the outer boundary of $U'$ (which is an equipotential of the polynomial $p$) is in the escaping set $U^{+}$ (more precisely in the set $V^{+}$).

\begin{figure}[htb]
\begin{center}
\includegraphics[scale=1.05, bb = 225 535 370 700]{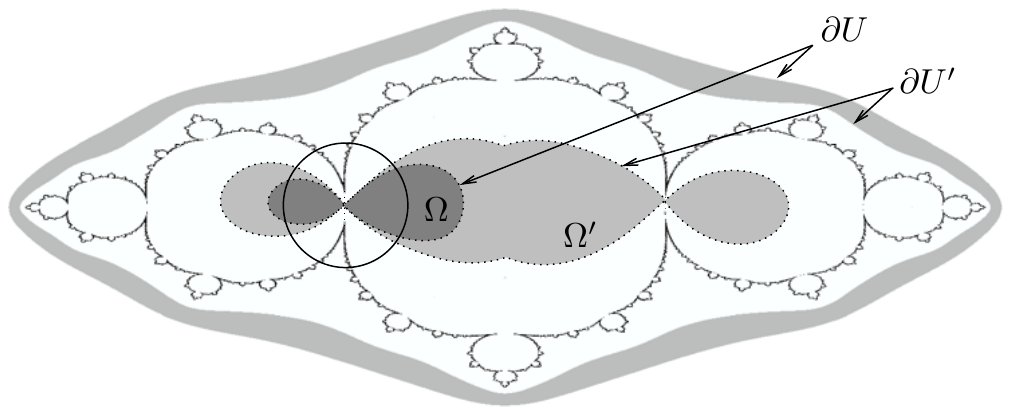}
\end{center}
\caption{The polynomial $p(x)=x^{2}-\frac{3}{4}$ has a parabolic fixed point at $-\frac{1}{2}$ and locally connected Julia set $J_{p}$. The corresponding neighborhoods $U$ and $U'$ are also shown, but $U'$ is not compactly contained in $U$, as in the hyperbolic case. Their boundaries touch at the parabolic fixed point.}
\end{figure}

We endow $U'$ with the Poincar\'e metric of $U$. The set $U'$ is contained in $U$ and $p:U'\rightarrow U$ is a covering map, hence expanding. 
However $U'$ is not relatively compact in $U$, so there is no constant of uniform expansion.

Define
\begin{equation}\label{eq:nbdV}
V := \left(U'\times \D_{r}-(B\cup B')\right)\cup \left(W^{-}_{B}\cup W^{-}_{B'}\right)\!,
\end{equation}
where $B$ and $B'$ are the two tubular neighborhoods defined above. 

\begin{figure}[htb]
\begin{center}
\includegraphics[scale=1, bb = 230 370 370 700]{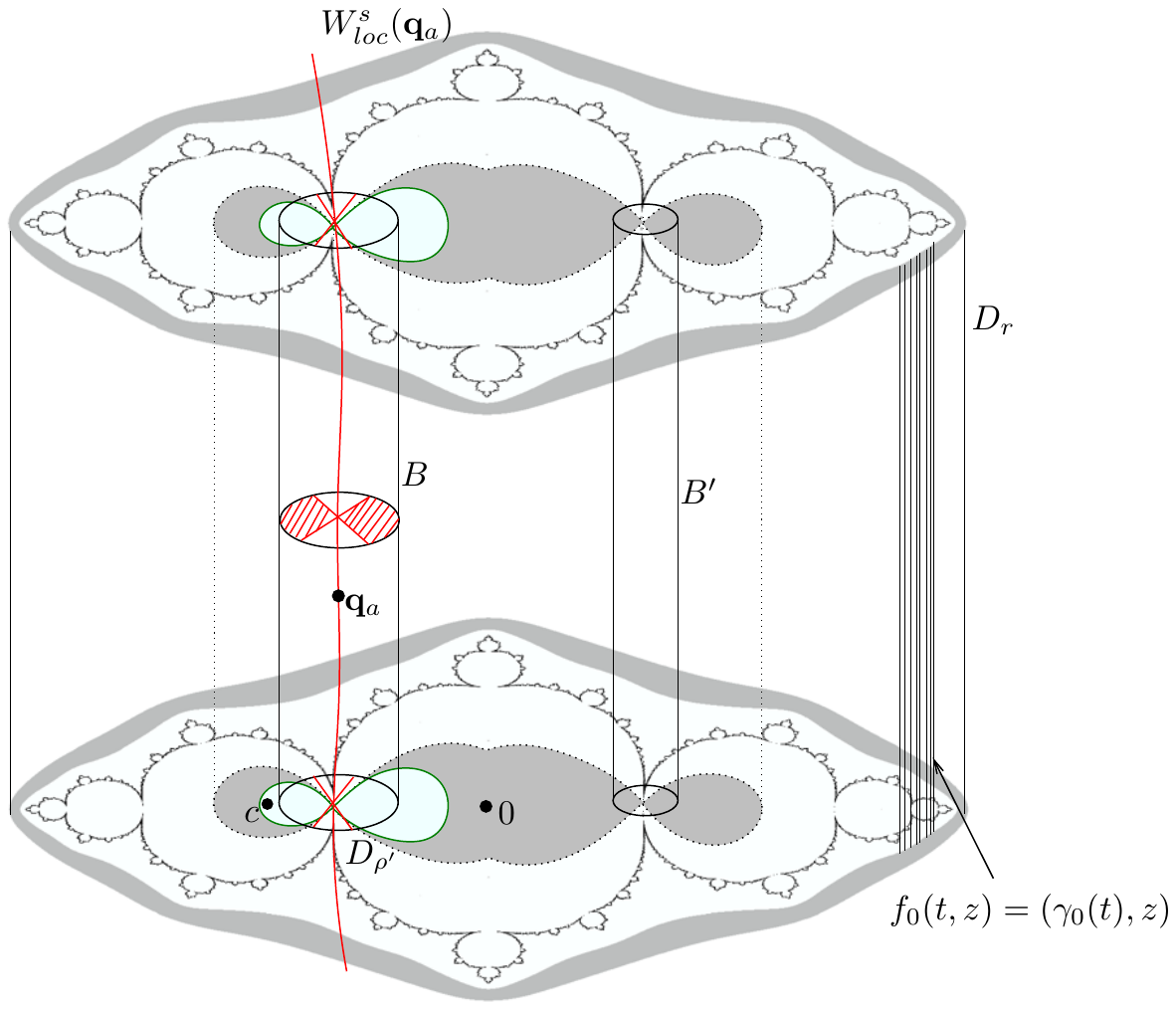}
\end{center}
\caption{A neighborhood $V$ of the set $J^{+}$ in $\D_{r}\times\D_{r}$. The map $\gamma_{0}$ used in the definition of the fiber $f_{0}(t,z)$ is the equipotential that gives the outer boundary of the set $U'$ (see also Equation \ref{eq:gamma}).}
\end{figure}

The vertical size of the neighborhood $V$ is $r>\rD$ where $r$ is chosen so that $ \D_{r}\subset U'\cup \Omega'$. We also require that  $\overline{H(V)}$ does not intersect the horizontal boundary of $V$, that is $|ax|<r$ for any $x\in U'$. 

The horizontal size of the neighborhood $V$ is given by an equipotential of the parabolic polynomial, contained entirely in the escaping set $U^{+}$. 

Let $a$ be small enough so that the following two conditions hold:
\begin{itemize}
  \item $r|a|<|p(x)-c_{0}|$ for any $x$ in $U'$. This is possible because we removed a disc around the critical value $c_{0}$ of the polynomial $p$, hence $\inf_{x\in U'}|p(x)-c_{0}|>0$.
  \item $2r|a|<d(\partial U'-\D_{\rho'}(q_{0}), \partial U)$. This assures that for all $x$ in $U'- \D_{\rho'}(q_{0})$ the disk of radius $2r|a|$ around $x$ belongs to $U$. In other words, the $2r|a|-$neighborhood of the set $U'-\D_{\rho'}(q_{0})$ is compactly contained in $U$.
\end{itemize}

Furthermore, we choose $a$ small enough so that in the construction of the set $V$ we make sure to remove points only from the interior of $K^{+}$ and not $J^{+}$. We only need to check for points that are outside the tubes $B$ and $B'$. This is guaranteed by the following lemma.

\begin{lemma}\label{lemma:omega} The removed set $\Omega'\times \D_{r} - (B\cup B')$ belongs to the interior of $K^{+}$.
\end{lemma} 
\proof  From the construction of the sets $U'$ and $\Omega'$, after at most $n+1$ iterations of the polynomial $p$,  points from $\Omega'$ are mapped to $\D_{\rho'}(q_{0})$. In fact, points from $\Omega-\D_{\rho'}(q_{0})$ are mapped to the region $p^{-\circ (n+1)}(S_{att})\cap A$ inside the annulus $A$ after at most $n+1$ iterates.

The \He map is given by $H(x,y)=(p(x)+a^{2}w +ay, ax)$. The $y$ component does not pose any problems as $|ax|<r$ for any $x\in \Omega'$. Let $(x,y)\in (\Omega-\D_{\rho'}(q_{0}))\times \D_{r}$. Suppose $k$ is the first iterate for which $p^{\circ k}(x)\in \D_{\rho'}(q_{0})$. Then $p^{\circ k}(x)$ $\in p^{-\circ (n+1-k)}(S_{att})\cap A$. After $k$ iterates, the distance between the $x$ coordinate of $H^{\circ k}(x,y)$ and $p^{k}(x)$ is at most $|a|\eta_{1}$, where $\eta_{1}$ is a constant which depends only on the parabolic polynomial $p$ and the integer $n$. Notice that $n$ is fixed and depends only on the polynomial $p$ and $\rho'$, which is also fixed.  However, if after $n+1$ iterates $p^{\circ(n+1)}(x)$ is too close to the outer boundary of $A$, then we take one more iterate and it is still in $A$, by construction of this annulus.  Based on the construction of $S_{att}$ we know that for $|a|<\eta_{0}/\eta_{1}$, $H^{\circ k}(x,y)\in W_{B}^{+}\cap A\times\D_{r}$, which is in the interior of $K^{+}$. 

This shows that $\Omega\times \D_{r} - B$ belongs to the interior of $K^{+}$, but a similar technical argument can be made for  $\Omega'\times \D_{r} - (B\cup B')$.
\qed

Let $\overline{V}$ denote the set $V$ together with $W^{s}_{loc}(\fq_{a})$ and $H^{-1}(W^{s}_{loc}(\fq_{a})) \cap B'$. In all other cases $\overline{X}$ denotes the closure of the set $X$.

\begin{lemma}\label{lemma: neighPolydisk} $J^{+}\cap \left(\D_{r}\times\D_{r}\right)=J^{+}\cap \overline{V}$.
\end{lemma}
\proof
The outer boundary of the set $V$ is an equipotential of the polynomial cross $\D_{r}$, which belongs to $U^{+}$. From the tubular neighborhood $B$ of the local stable manifold we removed only the attractive sectors $W^{+}_{B}$, which are contained inside the interior of $K^{+}$ union the local stable manifold $W^{s}_{loc}(\fq_{a})$.  From $B'$ we only removed the attractive sectors $W^{+}_{B'}$, which are contained inside the interior of $K^{+}$ union a preimage of the local stable manifold $H^{-1}(W^{s}_{loc}(\fq_{a})) \cap B'$. Outside of $B\cup B'$, we removed the set $\Omega'\times \D_{r} - (B\cup B')$ which belongs to the interior of $K^{+}$, as shown in Lemma \ref{lemma:omega}.
Therefore
\[
J^{+}\cap \left(\D_{r}\times\D_{r}\right)\  =\  \left(J^{+}\cap V\right) \cup W^{s}_{loc}(\fq_{a}) \cup \left(H^{-1}(W^{s}_{loc}(\fq_{a}))\cap B' \right)\ =\  J^{+}\cap \overline{V}.
\]
In this sense we say that $\overline{V}$ is a neighborhood of $J^{+}$ inside the bidisk $\D_{r}\times \D_{r}$.
\qed

\begin{cor}\label{lemma: J+invariance} $H(J^{+}\cap \overline{V}) \subset J^{+}\cap \overline{V}$ and $H^{-1}(J^{+}\cap \overline{V})\cap (\D_{r}\times \D_{r})\subset J^{+}\cap\overline{V}$.
\end{cor}

\begin{lemma}\label{lemma: neighJ}
    $J^{+}\cap \overline{V}=\bigcap_{n\geq 0}H^{-\circ n}(\overline{V}\cap \overline{U^{+}})$.
\end{lemma}
\proof
Let $q\in \bigcap_{n\geq 0}H^{-\circ n}(\overline{V}\cap \overline{U^{+}})$, where $\overline{U^{+}}=U^{+}\cup J^{+}$. Since all forward iterates of $q$ remain in the bounded set $\overline{V}$, $q$ cannot belong to $U^{+}$.
Hence $q\in J^{+}$. Suppose now that $q\in J^{+}\cap\overline{V}$. By construction of the neighborhood $V$, $H(J^{+}\cap\overline{V})\subset J^{+}\cap\overline{V}$, so all forward iterates of $q$ remain in $\overline{V}$. Hence $q\in \bigcap_{n\geq 0}H^{-\circ n}(\overline{V})$.
\qed

From Proposition \ref{lemma: neighPolydisk} we immediately get that the Julia set $J = \bigcap_{n\geq 0}H^{\circ n}(J^{+}\cap\overline{V})$.

\section{Infinitesimal metrics on $V$}\label{sec:metric}

On the set $V$ we will define two infinitesimal metrics with respect to which the derivative of the \He map is weekly expanding horizontally and strongly contracting vertically. 

To formalize our definitions, recall that $q_{0}$ is the parabolic fixed point of the quadratic polynomial $p$,
$B=\D_{\rho'}(q_0)\times \D_r$ and
$B''=\D_{\rho''}(q_0)\times \D_r$, where $0<\rho''<\rho'$. We have chosen $a$ small
enough so that the local stable manifold $W_{loc}^s(\fq_a)$ of the
semi-parabolic fixed point $\fq_a$ is contained in $B''$ and that Equation \ref{eq:Delta-} is satisfied.
In addition, the set
$U'$ is compactly contained in $U$ outside the disk
$\D_{\rho''}(q_{0})$.

\begin{defn}[\textbf{Euclidean metric}] In the repelling sectors $W^-_B$ of the tubular neighborhood $B$ of
the local stable manifold of the semi-parabolic fixed
point, we have a natural metric defined as a pull-back of the Euclidean metric from the
normalizing coordinates by $\phi_{a}:W^-_B\rightarrow W^-\subset
\D_{\rho}\times\D_r$,
\begin{equation}
\mu_{B}((x,y),(\xi,\eta)):= \max\left\{
|\widetilde{\xi}|,|\widetilde{\eta}|\right\},
\end{equation}
where
 $\phi_{a}$
is the
coordinate transformation from Lemma \ref{thm:UniformNF}, $(\widetilde{\xi},\widetilde{\eta})=D\phi_{a}\big{|}_{(x,y)}\left(\xi,\eta\right)$ and $|\widetilde{\xi}|$ and $|\widetilde{\eta}|$ represent the
length of $\widetilde{\xi}$ and $\widetilde{\eta}$ with respect to
the Euclidean metric.
\end{defn}

\begin{remarka} By
construction, the coordinate transformation $\phi_{a}$ takes
horizontal curves to horizontal curves. Therefore, if we choose a point $(x,y)\in B$ and a tangent vector
$(\xi,0)\in T_{(x,y)}B$, then 
$D\phi_{a}\big{|}_{(x,y)}(\xi,0)=(\widetilde{\xi}, 0)$ and $\mu_{B}\left((x,y),(\xi,0)\right)=|\widetilde{\xi}|$.
\end{remarka}

\begin{defn}[\textbf{Poincar\'e metric}]\label{def:Poincaremetric}
The set $V$, outside of a small neighborhood
$B''$ of the local stable manifold
$W^s_{loc}(\fq_{a})$ of the semi-parabolic fixed point
$\fq_{a}$, is contained in the product space $U\times \D_{r}$. 
On $V-B''$ we will use a product metric $\mu_{U}\times \mu_E$ of the Poincar\'{e} metric $\mu_{U}$ of the set $U$ and the Euclidean metric $\mu_E$ on the vertical disk $\D_r$.  Tangent vectors $(\xi,\eta)$ from $T_{(x,y)}V-B''$ will be
measured with respect to the metric
\begin{equation}
\mu_P((x,y),(\xi,\eta)) :=\max(\mu_{U}(x,\xi),|\eta|),
\end{equation}
where $|\eta|$ is the absolute value of the complex number $\eta$.
\end{defn}

\begin{defn}[\textbf{Combining the metrics}] Let $B'=(H^{-1}(B)-B)\cap V$ be one of the preimages of $B$ in $V$ as in Section \ref{subsec: V-parabolic}. Choose a number $M$ such that
\begin{equation}\label{eq:supM} M\geq\sup\limits_{\substack{(x,y)\in B'\\ (\xi,\eta)\in
\mathcal{C}^{h,P}_{(x,y)}}} \frac{2\cdot\mu_{P}\left((x,y),(\xi,\eta)\right)}{
\mu_{B}\left(H(x,y),DH_{(x,y)}(\xi,\eta)\right)}.
\end{equation}
Define as in \cite{DH} $\mu:=\inf\{\mu_{P},M\mu_{B}\}$, where the
infimum is taken pointwise between the metrics on $V$.
\end{defn}

\begin{remarka}
Note that the supremum from \ref{eq:supM} is a finite number. Since $(\xi,\eta)\in
\mathcal{C}^{h,P}_{(x,y)}$, the numerator $\mu_{P}\left((x,y),(\xi,\eta)\right)$ is equal to $\mu_{U}(x,\xi)$ which is bounded above for all $(x,y)$ in $B'$, because $B'$ is far away from the boundary of the set $U$, so the Poincar\'e metric $\mu_{U}$ is finite on $B'$. The denominator is bounded away from zero, because the tangent vector $DH_{(x,y)}(\xi,\eta)$ belongs to the horizontal cone $\mathcal{C}^{h,B}(H(x,y))$ when the vector $(\xi,\eta)$ belongs to the horizontal cone 
$\mathcal{C}^{h,P}_{(x,y)}$.
\end{remarka}

\begin{remarka}\label{remark: metricBPC}
The constant $M$ is chosen so that the \He map expands in horizontal cones with respect to the combined metrics, as we will see in  Theorem \ref{thm: mu-expansion}. By eventually increasing $M$ we can assume that inside horizontal cones on
$\partial \D_{\rho'}(q_0)\times
\D_{r}$, the infimum of the two metrics is attained
by the Poincar\'e metric $\mu_{P}$. As we approach the point $q_{0}$ which belongs to the boundary of $U$, the Poincar\'e metric explodes in horizontal cones, whereas the pull-back Euclidean metric $\mu_{B}$ is finite, so the infimum of the two metrics will be realized by the metric $\mu_{B}$. By eventually reducing $a$, we can 
assume that inside horizontal cones on $\partial \D_{\rho''}(\fq_0)\times
 \D_{r}$, the infimum is attained by the pull-back
metric $\mu_B$.  
\end{remarka}

\begin{thm}[\textbf{$\mu$-Expansion}] \label{thm: mu-expansion} Consider $(x,y),(x_1,y_1)\in
V$ with $H(x,y)=(x_1,y_1)$. Let $(\xi,\eta)$ and $(\xi_1,\eta_1)$ be two
tangent vectors such that  $DH_{(x,y)}(\xi,\eta)=(\xi_1,\eta_1)$ and $(\xi,\eta)$ belongs to the horizontal cones defined at $(x,y)$, i.e. $\mathcal{C}^{h,P}_{(x,y)}$ and/or $\mathcal{C}^{h,B}_{(x,y)}$.

Then the \He map is strictly but not strongly expanding with respect to $\mu$, that is
\[
\mu\left((x_1,y_1),(\xi_1,\eta_1)\right)>\alpha(x,y)\cdot \mu\left((x,y),(\xi,\eta)\right), \ \ \mbox{where}\ \ \alpha(x,y)>1,
\]
and $\alpha(x,y)$ is a constant in cases $(a)$, $(c )$ and $(d)$, and $\alpha(x,y)=\mathcal{E}(x,y)$, the expansion factor of the pull-back metric $\mu_{B}$, in case $(b)$. The expansion factor $\alpha(x,y)\rightarrow1$ if and only if $(x,y)$ tends to $W^{s}_{loc}(\fq_{a})$, the local stable manifold of the semi-parabolic fixed point.
\end{thm}
\proof There are four cases to consider:
\begin{itemize}
\item[(a)] Suppose
 \begin{eqnarray*}
 \mu\left((x,y),(\xi,\eta)\right)&=&\mu_P\left((x,y),(\xi,\eta)\right)\ \ \ \mbox{and}\\
  \mu\left((x_1,y_1),(\xi_1,\eta_1)\right)&=&\mu_P\left((x_1,y_1),(\xi_1,\eta_1)\right).
  \end{eqnarray*}
Since the Poincar\'e metric is smaller than the pull-back metric, it means that the points $(x,y)$ and $(x_{1},y_{1})$ are not very close to $q_{0}\times\D_{r}$. In particular by Remark \ref{remark: metricBPC} they must lie outside $\D_{\rho''}(q_{0})\times \D_{r}$.
By Proposition \ref{prop:cones-Ho-product-metric},
\[
\mu\left((x_1,y_1),(\xi_1,\eta_1)\right)>k\cdot
\mu\left((x,y),(\xi,\eta)\right).
\]

  \item[(b)]Suppose 
 \begin{eqnarray*}
 \mu((x,y),(\xi,\eta))&=&M\mu_B((x,y),(\xi,\eta))\ \ \ \mbox{and}\\
  \mu((x_1,y_1),(\xi_1,\eta_1))&=&M\mu_B((x_1,y_1),(\xi_1,\eta_1)).
  \end{eqnarray*}
 By Proposition \ref{lemma:HenonSectors} we  get 
$
\mu\left((x_1,y_1),(\xi_1,\eta_1)\right)>\mathcal{E}(x,y)\cdot
\mu\left((x,y),(\xi,\eta)\right)$. The expansion factor $\mathcal{E}(x,y)=1+\frac{\epsilon_{1}}{2}|\widetilde{x}|^{q}$, where $(\widetilde{x},\widetilde{y})=\phi_{a}(x,y)$.
  \item[(c)] Suppose 
   \begin{eqnarray*}
  \mu\left((x,y),(\xi,\eta)\right)&=&M\mu_B\left((x,y),(\xi,\eta)\right)\ \ \ \mbox{and}\\
  \mu\left((x_1,y_1),(\xi_1,\eta_1)\right)&=&\mu_P\left((x_1,y_1),(\xi_1,\eta_1)\right).
  \end{eqnarray*}  
By Remark \ref{remark: metricBPC} above, the point $(x_1,y_1)$ cannot be too close to  $q_{0}\times \D_{r}$ and it must stay outside the small tube $B''$.
By Proposition \ref{prop:cones-Ho-product-metric}, 
we have
\begin{eqnarray*}
\mu_P\left((x_1,y_1),(\xi_1,\eta_1)\right)&>&k\cdot
\mu_P\left((x,y),(\xi,\eta)\right)\\
&\geq& k \cdot M\mu_B\left((x,y),(\xi,\eta)\right)=k\cdot
\mu\left((x,y),(\xi,\eta)\right).
\end{eqnarray*}
  \item[(d)] Suppose 
   \begin{eqnarray*}
   \mu\left((x,y),(\xi,\eta)\right)&=&\mu_P\left((x,y),(\xi,\eta)\right)\ \ \ \mbox{and}\\
  \mu\left((x_1,y_1),(\xi_1,\eta_1)\right)&=&M\mu_B\left((x_1,y_1),(\xi_1,\eta_1)\right).
  \end{eqnarray*}
  In this case there are two subcases to consider:
 \begin{itemize}
 \item[(i)]
  If $(x,y)\in B'$, then by the choice of the constant $M$ we have
\begin{eqnarray*}
\mu_P\left((x,y),(\xi,\eta)\right)&<&\frac{2\cdot
\mu_P\left((x,y),(\xi,\eta)\right)}{\mu_B\left((x_1,y_1),(\xi_1,\eta_1)\right)}\cdot \frac{1}{2}
\mu_B\left((x_1,y_1),(\xi_1,\eta_1)\right)\\
&<&\frac{1}{2}\cdot M\mu_B\left((x_1,y_1),(\xi_1,\eta_1)\right)
\end{eqnarray*}
hence $\mu\left((x_1,y_1),(\xi_1,\eta_1)\right)>2\cdot
\mu\left((x,y),(\xi,\eta)\right)$.

\item[(ii)]
If $(x,y)\in B$, and the Poincar\'e metric is smaller than the
pull-back metric, then $(x,y)$ must be outside the small tube
$B''$ which encloses the local
stable manifold $W^s_{loc}(\fq_a)$. If we denote by
$
k':=\inf\limits_{\substack{(x,y)\in V-B''\\ |a|<\delta}}\mathcal{E}(x,y)
$ the
infimum of the expansion rate $\mathcal{E}(x,y)$ outside $B''$, then $k'>1$. By Proposition \ref{lemma:HenonSectors},
we know that
$
\mu_B\left((x_1,y_1),(\xi_1,\eta_1)\right)> \mathcal{E}(x,y)\cdot
\mu_B\left((x,y),(\xi,\eta)\right)
$.
Therefore
\begin{eqnarray*}
M\mu_B\left((x_1,y_1),(\xi_1,\eta_1)\right)&>& \mathcal{E}(x,y)\cdot
M\mu_B\left((x,y),(\xi,\eta)\right)\\
&\geq& k'\cdot \mu_P\left((x,y),(\xi,\eta)\right)\!,
\end{eqnarray*}
hence $\mu\left((x_1,y_1),(\xi_1,\eta_1)\right)>k'\cdot
\mu\left((x,y),(\xi,\eta)\right)$.
\qed
\end{itemize}
\end{itemize}

\section{Vertical and horizontal cones}\label{sec:cones}
\begin{defn}\label{def:vcones} In Section \ref{subsec: U-normalizing} we gave the definition \ref{def:hcone-in-petal} of a horizontal cone at a point $(x,y)$ from the set
$\D_{\rho}\times \D_r$, namely:
\[
\mathcal{C}^{h}_{(x,y)}=\left\{(\xi,\eta)\in
T_{(x,y)}\D_{\rho}\times\D_r,\ |\xi|>|\eta|\right\}.
\]
We will now define the vertical cone at a point $(x,y)$ from the set
$\D_{\rho}\times \D_r$ to be
\[
\mathcal{C}^{v}_{(x,y)}=\left\{(\xi,\eta)\in
T_{(x,y)}\D_{\rho}\times\D_r,\ |\xi|<|x|^{2q}|\eta|\right\}.
\]
\end{defn}

In Proposition \ref{prop:hcone-in-petal}, we showed that horizontal cones are invariant under $D\widetilde{H}$. Moreover, by Equation \ref{prop: expansion-hcones-in-petal} of Proposition \ref{lemma:HenonSectors}, the \He map expands the length of vectors from the horizontal cone $\mathcal{C}^{h}_{(x,y)}$ by a factor of $ 1+\frac{\epsilon_1}{2}|x|^q$. Now we will show that the vertical cones are invariant under
$D\widetilde{H}^{-1}$ and that $D\widetilde{H}^{-1}$ is expanding in the vertical direction.

\begin{prop}[\textbf{Vertical cones}]\label{prop:conesE}
Consider $(x,y)$ and $(x_1,y_1)$ in the repelling sectors of
$\D_{\rho}\times \D_r$ such that $\widetilde{H}(x,y)=(x_1,y_1)$.
Then
\[
D\widetilde{H}^{-1}_{(x_1,y_1)}\left(\mathcal{C}^{v}_{(x_1,y_1)}\right)\subset
Int\  \mathcal{C}^{v}_{(x,y)}
\]
and $ \big{\|} D\widetilde{H}^{-1}_{(x_{1},y_{1})}(\xi',\eta')\big{\|} \geq \frac{1}{|a|^{2}+1/2} \|(\xi',\eta')\| $ for $(\xi',\eta')\in \mathcal{C}^{v}_{(x_{1},y_{1})}$.
\end{prop}
\proof Let $(\xi',\eta')\in \mathcal{C}^{v}_{(x_1,y_1)}$ and
$(\xi,\eta)=D\widetilde{H}^{-1}_{(x_{1},y_{1})}(\xi',\eta')$. We need to
show that $(\xi,\eta)\in \mathcal{C}^{v}_{(x,y)}$ so we compute as before
\begin{eqnarray*}
\xi' &=&\lambda\left(1+(q+1)x^q+\partial_{x} g_a(x,y)\right)\xi+\lambda\partial_{y} g_a(x,y)\eta\\
\eta' &=&\partial_{x} (x h_{a})(x,y)\xi+\left(\mu+\partial_{y} (xh_{a})(x,y)\right)\eta\
\end{eqnarray*}
and estimate
\begin{eqnarray*}
|\xi'| &>& \left(|1+(q+1)x^q|-m|x|^{2q}\right)|\xi|-
M_a|x|^{2q+2}|\eta|\\
|\eta'|&<& N_a|\xi|+\left(|\mu|+N_a\right)|\eta|.
\end{eqnarray*}
Since $(\xi',\eta')$ belongs to the vertical cone at $(x_1,y_1)$,
we also know that
\begin{eqnarray*}
|\xi'|<|x_1|^{2q}|\eta'|<|x|^{2q}|1+x^q+ g_a(x,y)/
x|^{2q}|\eta'|<|x|^{2q}M_{1}^{2q}|\eta'|,
\end{eqnarray*}
where $M_{1}$ is the supremum of $|1+x^q+ g_a(x,y)/ x|$ on the repelling
sectors $W^{-}$ of the tubular neighborhood $B$, that is
\begin{equation}\label{eq:M1}
M_{1}:=\sup_{(x,y)\in W^{-},\ |a|<\delta}\big{|}1+x^q+g_a(x,y)/ x\big{|}.
\end{equation}
Clearly $M_{1}>0$. In fact we could take a constant $M_{1}>1$ because $Re(x^{q})>\epsilon_{1}|x|^{q}$ in the repelling sectors $W^{-}$.
By combining these inequalities we get
\[
\left(|1+(q+1)x^q|-m|x|^{2q}\right)|\xi|- M_a|x|^{2q+2}|\eta| <
M_{1}^{2q}N_a|x|^{2q}|\xi|+M_{1}^{2q}(|\mu|+N_a)|x|^{2q}|\eta|.
\]
After regrouping the terms, we obtain
\[
|\xi|<\frac{A_{2}}{A_{1}}|x|^{2q}|\eta|
\]
where $A_{1}$ and $A_{2}$ are defined as follows
\begin{eqnarray*}
A_{1} &:=& |1+(q+1)x^q|-(m+M_{1}^{2q}N_a)|x|^{2q} \\
A_{2} &:=& M_{1}^{2q}(|\mu|+N_a)+M_a|x|^2.
\end{eqnarray*}
Since $x$ is chosen from the repelling sectors we have $|1+(q+1)x^q|-m|x|^{2q}>1+\epsilon_{1}|x|^{q}$. The bounds $N_{a}$, $M_{a}$ and the eigenvalue $\mu$ all depend on $a$, and they tend to $0$ as $a\rightarrow 0$, so for $|a|$
small we can assume that $A_{1}>\frac{2}{3}$ and $A_{2}<\frac{1}{3}$. Hence
$(\xi,\eta)\in Int\ \mathcal{C}^v_{(x,y)}$. 

We will now show that inside the vertical cones,
$D\widetilde{H}^{-1}$ is expanding with respect to the Euclidean
metric. We have
\begin{eqnarray*}
|\eta'|&<& N_{a}|\xi|+\left(|\mu|+N_{a}\right)|\eta| < N_{a}\frac{A_{2}}{A_{1}}|x|^{2q}|\eta|+\left(|\mu|+N_{a}\right)|\eta| \\
&<&\left(\frac{1}{2}N_{a}|x|^{2q}+|\mu|+N_{a}\right)|\eta|<\left(|\mu|+\frac{3}{2}N_{a}\right)|\eta|,
\end{eqnarray*}
since $|x|<\rho<1$.  
We obtain $|\eta|>
\frac{1}{|\mu|+\frac{3}{2}N_{a}}|\eta'|$.
For $a$ sufficiently small we can assume that $N_{a}<1/3$, therefore
$|\eta|> \frac{1}{|\mu|+1/2}| \eta'|=\frac{1}{|a|^{2}+1/2}| \eta'|$.

Since $(\xi,\eta)$ and $(\xi',\eta')$
belong to vertical cones, we have $\| (\xi,\eta)\|= \max(|\xi|,|\eta|) = |\eta|$ and $\| (\xi',\eta')\| = \max(|\xi|',|\eta'|) = |\eta'|$, hence $D\widetilde{H}^{-1}$
expands in the vertical cones with a factor strictly greater than 1. 
\qed

The vertical cones that we have introduced in Definition \ref{def:vcones}, are taken with
respect to the Euclidean metric, in the normalized coordinates
around the local stable manifold of the semi-parabolic fixed point.

\begin{defn}\label{def:cones-pB} Let $\mathcal{C}^{v,B}_{(x,y)}:=D\phi_a^{-1}\left(\mathcal{C}^{v}_{\phi_{a}(x,y)}\right)$ and $\mathcal{C}^{h,B}_{(x,y)}:=D\phi_a^{-1}\left(\mathcal{C}^{h}_{\phi_{a}(x,y)}\right)$ denote the pull-back of the vertical cone $\mathcal{C}^{v}_{\phi_{a}(x,y)}$ and respectively of the horizontal cone $\mathcal{C}^{h}_{\phi_{a}(x,y)}$ defined in \ref{def:vcones}, from the normalized coordinates $\D_{\rho}\times \D_r$ into $B$, by the change of coordinate function $\phi_a$. 
\end{defn}


On the set $V$, outside of a small neighborhood
$\D_{\rho''}(q_{0})\times \D_r$ of the local stable manifold
$W^s_{loc}(\fq_{a})$ of the semi-parabolic fixed point
$\fq_{a}$, we will use a product metric $\mu_{U}\times \mu_E$ of the Poincar\'{e} metric $\mu_{U}$ of the set $U$ and the Euclidean metric $\mu_E$ on the vertical disk $\D_r$.  
We will also define vertical and horizontal cones with respect to the product metric and show invariance under $DH^{-1}$, respectively under $DH$.

Let us notice first that $p:U'\rightarrow U$ is a covering map, hence a local isometry from the set $U'$ endowed with the Poincar\'e metric of $U'$ to the set $U$ endowed with the Poincar\'e metric of $U$.  Since $U'$ is contained in $U$, the inclusion map is contracting with respect to the Poincar\'e metrics of $U'$ and $U$. Therefore, if we consider $x\in U'$ and $\xi$ a complex tangent vector, then 
\begin{equation}\label{eq:nonstrict-contraction}
\mu_{U}(p(x),p'(x)\xi)= \mu_{U'}(x,\xi)> \mu_{U}(x,\xi),
\end{equation} 
hence the polynomial $p$ is expanding with respect to the Poincar\'e metric of $U$. 
The boundaries of the sets $U$ and $U'$ touch at the parabolic fixed point $q_{0}$, so there is no constant of uniform expansion. However, the set $U'-\D_{\rho''}(q_{0})$ is compactly contained in $U$, therefore there exists a constant $k_{0}>1$ such that 
\begin{equation}\label{explanation-polynom}
\mu_{U}\left(p(x),p'(x)\xi\right)> k_{0}\cdot \mu_{U}(x,\xi),\ \ \mbox{for all}\ \ x\in U'-\D_{\rho''}(q_{0}).
\end{equation}
The constant $k_{0}$ from inequality \ref{explanation-polynom} is just $\liminf\limits_{x\in U'-\D_{\rho''}(q_{0})}\frac{\rho_{U}(p(x))|p'(x)|}{\rho_{U}(x)}$, where $\rho_{U}$ is the density function of the Poincar\'e metric $\mu_{U}$, that is $\rho_{U}$ is a positive continuous function such that
\[
\mu_{U}(x,\xi)=\rho_{U}(x)|\xi|,\ \ \mbox{for all}\ x\in U\ \mbox{and}\ \xi\in T_{x}U.
\] 
By inequality \ref{eq:nonstrict-contraction}, we know that $k_{0}\geq 1$. Suppose by contradiction that $k_{0}=1$. There exists a sequence of points $ \{x_{n}\}_{n\geq 1}$ in $U'-\D_{\rho''}(q_{0})$ such that $ \frac{\rho_{U}(p(x_{n}))|p'(x_{n})|}{\rho_{U}(x_{n})}\rightarrow 1$. 
 The sequence $x_{n}$ is bounded, therefore there exists a convergent subsequence $x_{n_{k}}\rightarrow x^{*}$, where $x^{*}$ belongs to the closure of $U'-\D_{\rho''}(q_{0})$. However, the limit point $x^{*}$ cannot belong to the part of the boundary of $U'-\D_{\rho''}(q_{0})$ given by $\partial U'-\D_{\rho''}(q_{0})$, because the latter is compactly contained in $U$, so $\rho_{U}(x^{*})$ is finite. However $p(\partial U'-\D_{\rho''}(q_{0}))\subset \partial U$, so the Poincar\'e metric of $\rho_{U}(p(x^{*}))$ would be infinite. Note also that $|p'(x^{*})|$ is different from $0$, as the closure of $U'$ does not contain the critical point $0$ of the polynomial $p$. This contradicts the fact that $\lim\limits_{k\rightarrow \infty} \frac{\rho_{U}(p(x_{n_{k}}))|p'(x_{n_{k}})|}{\rho_{U}(x_{n_{k}})}= 1$. The only possibility left is that $x^{*}$ belongs to $U'$, but then, by  inequality \ref{eq:nonstrict-contraction} we know that $\rho_{U}(p(x^{*}))|p'(x^{*})|>\rho_{U}(x^{*})$. In conclusion $k_{0}$ must also be strictly greater than $1$.

Lastly, we also make the observation that on the set $U'-\D_{\rho''}(q_{0})$, the Poicar\'{e} metric
$\mu_{U}$ is bounded above and below by the Euclidean metric, that is, there exist two positive constants
$m_1$ and $m_2$ such that $m_1 < \rho_{U}(x)<m_2$ for any $x\in U'-\D_{\rho''}(q_{0})$.

\begin{defn}\label{def:conesP}
Let $\con<1$ to be chosen later. Define the vertical cone at a point $(x,y)$
from the set $U'\times \D_r - \D_{\rho''}(q_{0})\times \D_r$ to
be
\[
\mathcal{C}^{v,P}_{(x,y)}=\left\{(\xi,\eta)\in
T_{(x,y)}U'\times\D_r,\ \ \mu_{U}(x,\xi)<\con |\eta|\right\}.
\]
Define the horizontal cone at a point $(x,y)$ from the set
$U'\times \D_r-\D_{\rho''}(q_{0})\times \D_r$ as
\[
\mathcal{C}^{h,P}_{(x,y)}=\left\{(\xi,\eta)\in
T_{(x,y)}U'\times\D_r,\ \ \mu_{U}(x,\xi)>|\eta|\right\}.
\]
\end{defn}

\begin{prop}[\textbf{Vertical cones}]\label{prop:cones-Ve-product-metric}
Consider $(x,y)$ and $(x',y')$ in $U'\times
\D_r-\D_{\rho''}(q_{0})\times \D_r$ such that $H(x',y')=(x,y)$.
Then
\[
DH^{-1}_{(x,y)}\left(\mathcal{C}^{v,P}_{(x,y)}\right)\subset Int\
\mathcal{C}^{v,P}_{(x',y')}
\]
and $\big{\|} DH^{-1}_{(x,y)}(\xi,\eta)\big{\|} \geq \frac{1}{|a|} \|(\xi,\eta)\| $ for $(\xi,\eta)\in
\mathcal{C}^{v}_{(x,y)}$.
\end{prop}
\proof
Let $(\xi,\eta)\in \mathcal{C}^{v}_{(x,y)}$ and denote by
$(\xi',\eta') = DH^{-1}_{(x,y)}(\xi,\eta)$.
Since
\begin{equation*}
\hvec{x'}{y'}=\hvec{\frac{y}{a}}{ \frac{
x-p(y/a)-a^2w}{a}}\ \ \ \mbox{and}\ \  \
DH_{(x,y)}^{-1}=\left[
                        \begin{array}{rr}
                          0 & \frac{1}{a} \\
			\vspace{-0.4cm}\\
                          \frac{1}{a} & -\frac{2y}{a^3}
                        \end{array}
                    \right],
\end{equation*}
we can compute $ \xi'=\frac{1}{a}\eta $ and $ \eta' =
\frac{1}{a}\left(\xi-\frac{2x'}{a} \eta\right)$. The vector
$(\xi,\eta)$ belongs to the vertical cone, so $
\mu_{U}(x,\xi)=\rho_{U}(x)|\xi|<\con  |\eta|$. This
implies
\begin{equation}\label{eq: cone-v}
 |\xi|<\frac{\con}{m_1} |\eta|.
\end{equation}
We can evaluate
\begin{equation}\label{eq: lh-1}
 \mu_{U}(x',\xi')=\rho_{U}(x') \frac{|\eta|}{|a|}\leq \frac{m_2}{|a|} |\eta|.
\end{equation}
Next, we compute using inequality \ref{eq: cone-v}
\begin{equation}\label{eq: rh-1}
 |\eta'|>\frac{1}{|a|} \left|\xi-\frac{2x'}{a} \eta\right|>\frac{1}{|a|}  \left(\frac{|2x'|}{|a|} -\frac{\con}{m_1}\right) |\eta|.
\end{equation}
The point $x'$ belongs to the set $U'$. The set $U'$ does not contain a neighborhood of the critical
point $0$ of the polynomial $p$. Hence there exists a lower bound
$r_1>0$ such that $r_1<|2x'|$. Choose $a$ small, so $
\frac{r_1}{|a|}-\frac{\con}{m_1} > \max\left\{\frac{2 m_2}{\con}, 1\right\}$. By combining Equations \ref{eq: lh-1} and \ref{eq: rh-1} we get
$ \mu_{U}(x',\xi') < \frac{\con}{2} |\eta|$.
Therefore we have shown the invariance of vertical cones $
DH^{-1}_{(x,y)}\left(\mathcal{C}^{v,P}_{(x,y)}\right)\subset Int\ \mathcal{C}^{v,P}_{(x',y')}$. 

Inequality \ref{eq: rh-1} shows that $ |\eta'|>\frac{1}{|a|} |\eta|$, so
$DH^{-1}$ is expanding in the vertical cones.
Since $\| (\xi,\eta) \| = \max(\mu_{U}(x,\xi),|\eta|) = |\eta|$ and $\| (\xi',\eta') \| = \max(\mu_{U}(x',\xi'),|\eta'|) = |\eta'|$, we obtain $ \| (\xi',\eta')\|> \frac{1}{|a|} \|(\xi,\eta)\|$, as claimed. 
\qed

\begin{remarka}
The scalar $0<\con<1$ in the definition of the vertical cone will
typically be chosen less than $\left(\frac{\rho}{2}\right)^{2q}$, so
that on a neighborhood of the boundary of $B$, the vertical cones
$\mathcal{C}^{v,P}_{(x,y)}$ from Definition \ref{def:conesP}
are contained in the interior of the pull-back cones $\mathcal{C}^{v,B}_{(x,y)}$ from Definition \ref{def:cones-pB}. 
In this way we can assure that $DH^{-1}_{(x,y)}\left(\mathcal{C}^{v,P}_{(x,y)}\right)\subset Int\
\mathcal{C}^{v,B}_{(x',y')}$.
\end{remarka}

To fully show the invariance of the two types of vertical cones under $DH^{-1}$, we have one more case to cover. 

\begin{prop} Let $(x',y')\in B'$ and $(x,y)\in B$ such that $H(x',y')=(x,y)$. Then
\[ 
DH^{-1}_{(x,y)}\left(\mathcal{C}^{v,B}_{(x,y)}\right)\subset Int\
\mathcal{C}^{v,P}_{(x',y')}.
\]
\end{prop}
\proof Consider $(\xi, \eta)\in \mathcal{C}^{v,B}_{(x,y)}$ and $(\xi', \eta')=DH^{-1}_{(x,y)}(\xi,\eta)$. We have to show that $\mu_{U}(x',\xi')<\con |\eta'|$. Let $(\tilde{x}, \tilde{y})=\phi_{a}(x,y)$ and $(\tilde{\xi}, \tilde{\eta})=D\phi_{a}|_{(x,y)}(\xi,\eta)$. By Definition \ref{def:cones-pB}, the vector $(\xi, \eta)$ belongs to the vertical cone $\mathcal{C}^{v,B}_{(x,y)}$ if and only if $|\tilde{\xi}|< |\tilde{x}|^{2q}|\tilde{\eta}|$. The change of coordinate $\phi_{a}(x,y)$ is holomorphic with respect to $a$ and it is $\bigO(a)$ close to $(\phi(x),y)$, therefore there exists $\kappa_{\phi}>0$ such that when $a$ is small we have $|\xi|<\kappa_{\phi} |\eta|$.

By using the same computations as in Proposition \ref{prop:cones-Ve-product-metric} we obtain $ \eta' =\frac{1}{a}\left(\xi-\frac{2x'}{a} \eta\right)$ and $ \xi'=\frac{1}{a}\eta $ and we have estimates analogous to relations \ref{eq: lh-1} and \ref{eq: rh-1}:
\begin{eqnarray}
\mu_{U}(x',\xi') &=& \rho_{U}(x')|\xi'|< \frac{m_{2}}{|a|}|\eta| \\
 |\eta'| &>& \frac{1}{|a|}  \left(\frac{|2x'|}{|a|} |\eta|-|\xi|\right) > \frac{1}{|a|}  \left(\frac{r_{1}}{|a|}-\kappa_{\phi}\right) |\eta|. \label{eq: lungime}
\end{eqnarray}
Therefore, to show that $ \mu_{U}(x',\xi')<\con |\eta'|$ all we need to know is that $ m_{2}<\con\left(\frac{r_{1}}{|a|}-\kappa_{\phi}\right)$, which is obviously true when $a$ is small. 
\qed

\begin{prop}[\textbf{Horizontal cones}]\label{prop:cones-Ho-product-metric}
Let $(x,y)$ and $(x',y')$ in
$(U'-\D_{\rho''}(q_{0}))\times \D_r$ such that
$H(x,y)=(x',y')$. Then
\[
DH_{(x,y)}\left(\mathcal{C}^{h,P}_{(x,y)}\right)\subset Int\
\mathcal{C}^{h,P}_{(x',y')}
\]
and $\big{\|}DH_{(x,y)}(\xi,\eta)\big{\|} \geq k\cdot
\|(\xi,\eta)\|$ for $(\xi,\eta)\in \mathcal{C}^{h,P}_{(x,y)}$, where $k>1$ is a constant that depends on the polynomial $p$ which has a parabolic fixed point at $q_{0}$.
\end{prop}
\proof Let $(\xi,\eta)\in \mathcal{C}^{h,P}_{(x,y)}$ and
$(\xi',\eta')=DH_{(x,y)}(\xi,\eta)$. We show that
$(\xi',\eta')\in Int\ \mathcal{C}^{h,P}_{(x',y')}$ and $\mu_{U}(x',\xi')>k\cdot \mu_{U}(x,\xi)$. Since
\[
\hvec{x'}{y'}=\hvec{p(x)+a^2w+ay}{ax}\ \ \ \mbox{and}\ \  \
DH_{(x,y)}=\left[
                        \begin{array}{rr}
                          2x & a \\
                          a & 0 \\
                        \end{array}
                      \right],
\]
we can compute $ \xi'=2x\xi+a\eta$ and $ \eta' = a\xi$. The
vector $(\xi,\eta)$ belongs to the horizontal cone at $(x,y)$, so
\begin{equation}\label{eq: cone-h}
|\eta|<\mu_{U}(x,\xi)=\rho_{U}(x)|\xi|<m_2|\xi|.
\end{equation}
We evaluate
\begin{equation}\label{eq: lh-2}
\mu_{U}(x',\xi')=\rho_{U}\left(p(x)+a^2w+ay\right)|2x\xi+a\eta|.
\end{equation}
The density function $\rho_{U}$ of the Poincar\'e metric is
bounded above and below on the set $U'':=
U'-\D_{\rho''}(q_{0})$ since we removed a disk around the fixed
point $q_{0}$, where the boundaries of $U'$ and $U$ touch.
Since $\rho_{U}$ is $C^{\infty}$-smooth on $U''$, its derivative
$\rho'_{U}$ is also bounded. There exists a constant $\cU>0$ such
that
\begin{eqnarray}\label{eq: density-smooth}
\frac{\left|\rho_{U}\left(p(x)+a^2w+ay\right)-\rho_{U}\left(p(x)\right)\right|}{|a|}
&\leq& |aw+y|\cdot \sup\limits_{U''} \rho'_{U}\cdot
\frac{\rho_{U}(p(x))}{\inf\limits_{U''} \rho_{U}} \nonumber
\\
&<& \cU \cdot \rho_{U}\left(p(x)\right).
\end{eqnarray}
\noindent By \ref{explanation-polynom}, the polynomial $p$ is expanding on $U'-\D_{\rho''}(q_{0})$ with respect
to the Poincar\'e metric $\mu_{U}$, that is there exists a constant $k_{0}>1$  such that
\begin{equation}\label{eq: c_t-polynom2}
\rho_{U}\left(p(x)\right)|p'(x)\xi|>k_{0}\cdot \rho_{U}(x)
|\xi|.
\end{equation}
The set $U'$ avoids a uniform neighborhood of the critical point $0$, so
there exists as before $r_1$ such that $0<r_1<|2x|$ for any $x\in U'$. We
now turn back to relation \ref{eq: lh-2}. Using \ref{eq:
density-smooth}, \ref{eq: c_t-polynom2} and \ref{eq: cone-h} one
gets
\begin{eqnarray}\label{eq: rh-2}
\mu_{U}(x',\xi') &>&(1-\cU |a|)\cdot
\rho_{U}\left(p(x)\right)|2x\xi|\cdot\frac{|2x\xi+a\eta|}{|2x\xi|}
\nonumber \\
&>& (1-\cU |a|)\cdot k_{0} \cdot \rho_{U}(x)|\xi| \cdot
\left(1-|a|\frac{|\eta|}{|2x||\xi|}\right) \nonumber \\
&>& k_{0} \cdot (1-\cU |a|)\cdot \left(1-|a|\frac{m_2}{r_1}\right)\cdot
\rho_{U}(x)|\xi|.
\end{eqnarray}
The constant $k_{0}$ is strictly bigger than $1$. The factors $1-c|a|$ and
$1-|a|\frac{m_2}{r_1}$  can be made
arbitrarily close to $1$ by reducing $|a|<\delta$. Let  $\delta$
be sufficiently small so that
\begin{equation}
k := k_{0} \cdot (1-\cU \delta)\cdot
\left(1-\delta\frac{m_2}{r_1}\right) >1.
\end{equation}
When $|a|<\delta$, relation \ref{eq: rh-2} gives
\begin{equation}
\mu_{U}(x',\xi') > k \cdot \rho_{U}(x)|\xi|=k\cdot
\mu_{U}(x,\xi)
\end{equation}
which proves that $DH$ expands in the horizontal cones. Also from
relation \ref{eq: rh-2} we infer that
\begin{equation}
|a|\cdot \mu_{U}(x',\xi') > k \cdot \rho_{U}(x)|a\xi|> k \cdot
m_1 \cdot |\eta'|
\end{equation}
which proves that $DH_{(x,y)}\left(\mathcal{C}^{h,P}_{(x,y)}\right)\subset
Int\ \mathcal{C}^{h,P}_{(x',y')}$, so the horizontal cones are invariant under $DH$ for $a$ small. \qed

If both types of horizontal cones are defined at some point $(x,y)\in V$, we cannot say in general that one is contained in the other, i.e. $\mathcal{C}^{h,P}_{(x,y)} \subset \mathcal{C}^{h,B}_{(x,y)}$ or $\mathcal{C}^{h,B}_{(x,y)} \subset \mathcal{C}^{h,P}_{(x,y)}$.  However, since the derivative of the \He map contracts the vertical component of tangent vectors by a factor of $a$,  we can assume that $DH_{(x,y)}\left(\mathcal{C}^{h,P}_{(x,y)}\right)\subset
Int\ \mathcal{C}^{h,B}_{(x',y')}$ and $DH_{(x,y)}\left(\mathcal{C}^{h,B}_{(x,y)}\right)\subset
Int\ \mathcal{C}^{h,P}_{(x',y')}$ whenever both types of cones are defined at $(x,y)$ and/or at $(x',y')=H(x,y)$.

\section{Distance between vertical-like curves}\label{sec:distance}

In this section we work entirely in the normalized coordinates from Theorem \ref{thm:UniformNF}. The notion of vertical-like curves translates as follows
\begin{defn}
We will call an analytic curve $\gamma\subset\D_{\rho}\times \D_{r}$
{\it vertical-like} if $\gamma$ is the graph of an analytic function
$\phi:\D_r\rightarrow \D_{\rho}$, and for all points $(x,y)$ on
$\gamma$, the tangent vectors $(\xi, \eta)$ to $\gamma$ at $(x,y)$
belong to the vertical cone $\mathcal{C}^{v}_{(x,y)}$ from Definition \ref{def:vcones}.
\end{defn}

Let us now consider two vertical-like curves in the same
repelling sector of $\D_{\rho}\times \D_{r}$, that are entirely
contained in the escaping set
$U^+$.
Denote these vertical curves
\[
f_1(z)=(\varphi_1(z),z)\ \mbox{and}\ f_2(z)=(\varphi_2(z),z).
\]

\noindent Let $g_1(\D_r)$ be the image under $\widetilde{H}^{-1}$ of
$f_1(\D_r)$, contained inside $\D_{\rho}\times\D_r$. More precisely,
\[
\widetilde{H}^{-1}(f_1(\D_r))\cap (\D_{\rho}\times\D_r),
\]
is a vertical-like fiber, that we can describe as the graph of an
analytic function
\[
g_1(z)=(\varphi'(z),z),\ \mbox{where}\ \varphi':\D_r\rightarrow
\D_{\rho}.
\]
Similarly, let $g_2(\D_r)$ be $\widetilde{H}^{-1}(f_2(\D_r))\cap
(\D_{\rho}\times\D_r)$ reparametrized by the second coordinate
$g_2(z)=(\varphi''(z),z)$. Notice that $g_1(\D_r)$ and $g_2(\D_r)$
are vertical-like curves (by Proposition \ref{prop:conesE}), both
contained in some other repelling sector of $\D_{\rho}\times \D_r$
and in $U^{+}$.

Much like in the hyperbolic setting, we would like to show that
$\widetilde{H}$ expands the horizontal distance between
vertical-like curves. We will measure the horizontal distance with
respect to the standard Euclidean metric on $\D_{\rho}\times \D_r$.
Define
\[
d(f_1,f_2)=\|\varphi_1-\varphi_2\|=\sup\limits_{z\in\D_r}|\varphi_1(z)-\varphi_2(z)|
\]
Notice that this distance between vertical-like curves
is just the distance between the parametrizing functions $\varphi_1$
and $\varphi_2$ with respect to the $\sup$-norm.

\begin{thm}\label{thm:distE} Let $d(g_1,g_2)$ and $d(f_1,f_2)$ be the horizontal
distance between the vertical-like curves $g_1$, $g_2$ and
respectively $f_1$, $f_2$. Then
\[
d(g_1,g_2)<Cd(f_1,f_2),
\]
where $C=C(f_{1},f_{2})<1$, so the normalized \He maps $\widetilde{H}$ expands strictly the distance between the vertical-like curves $g_1$ and
$g_2$.
\end{thm}
\proof Let $z\in \D_r$ be arbitrarily chosen and denote by
$x'=\varphi'(z)$, and $x''=\varphi''(z)$. The points $(x',z)$ and
$(x'',z)$ lie on the vertical-like curves $g_1$ and $g_2$. 

\begin{figure}[htb] 
\begin{center}
\includegraphics[scale=0.93, bb = 310 515 370 730]{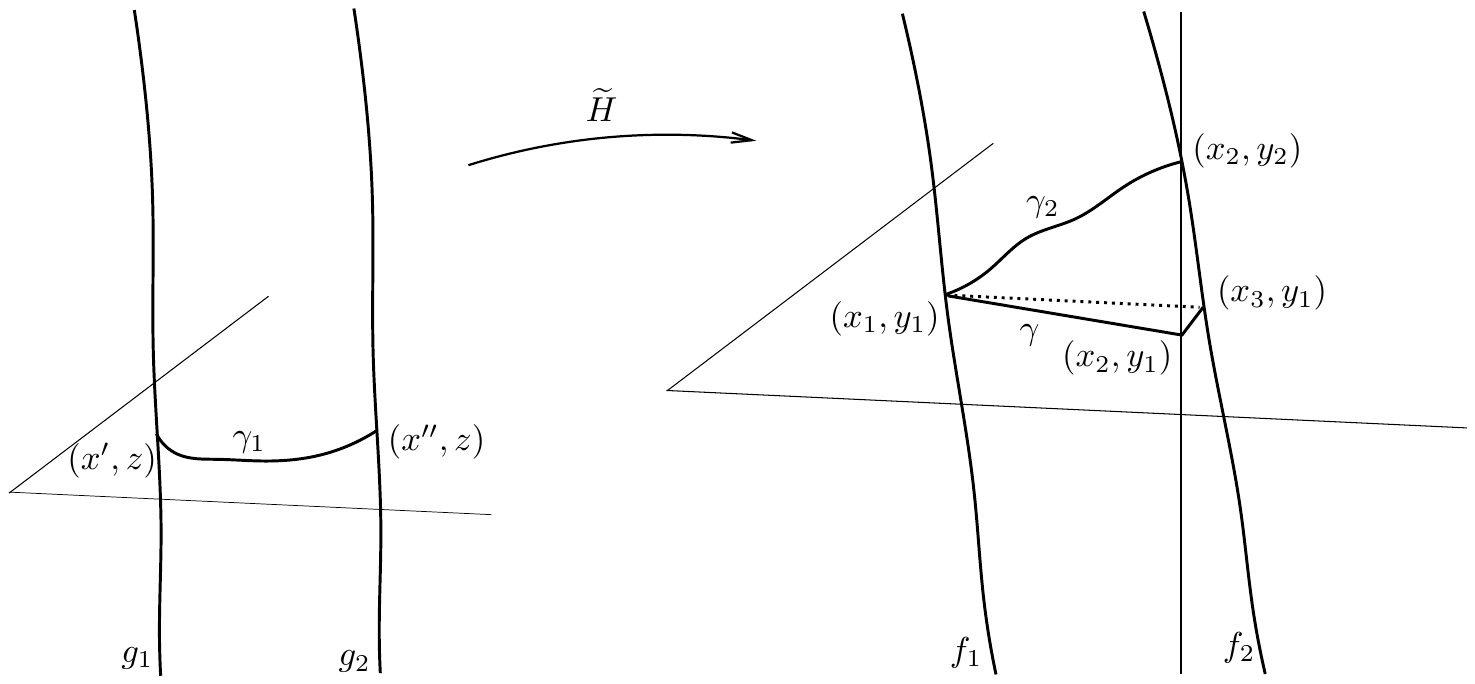}
\end{center}
\caption{Fibers $g_{1}, g_{2}$ and their image fibers $f_{1}, f_{2}$ under $\widetilde{H}$. }
\label{fig: FourFibers}
\end{figure}

Let
$(x_1,y_1)=\widetilde{H}(x',z)$ and $(x_2,y_2)=\widetilde{H}(x'',z)$
be the corresponding points on $f_1$ and $f_2$. Let
$(x_3,y_1)=(\varphi_2(y_2),y_2)$ be the point of intersection of the
curve $f_2$ with the horizontal plane $\C\times\{y_1\}$. Suppose
without loss of generality that $|x_2|\leq |x_1|$. Let
\[
\mathcal{I}(x_{1},x_{2}):=\int_{0}^{1}|tx_1+(1-t)x_2|^q dt.
\]

\begin{lemma}[\textbf{Step 1}]\label{lemma:x1x2} We have
\[
|x'-x''|<\left(1-\frac{\epsilon_{1}}{2M_{1}^{q}}\mathcal{I}(x_{1},x_{2})\right)|x_1-x_2|,
\]
where the constant $M_{1}$ is independent of $a$. The constant $\epsilon_{1}$ is given in Equation \ref{eq:eps1}.
\end{lemma}
\proof Choose a straight line in the $\C\times \{y_1\}$ plane,
\[
\gamma(t)=(x_{\gamma}(t),y_1),\ \mbox{ where }
x_{\gamma}(t)=tx_1+(1-t)x_2 \mbox{ and } t\in[0,1],
\]
connecting the points $(x_1,y_1)$
and $(x_2,y_1)$. There exists a horizontal curve
\[
\gamma_1(t):[0,1]\rightarrow \D_{\rho}\times\{z\},\ \
\gamma_1(t)=(x_{\gamma_1}(t),z),
\]
connecting the points $(x',z)=\gamma_1(0)$ and $(x'',z)=\gamma_2(0)$
and such that the projection of the curve $\gamma_2(t):=\widetilde{H}(\gamma_1(t))$ on the plane $\C\times\{y_1\}$
is exactly the straight line $\gamma(t)$. Formally, if we define
$pr:\D_{\rho}\times\D_r\rightarrow \D_{\rho}\times \{y_1\}$,
$pr(x,y)=(x,y_1)$, then $pr(\gamma_2(t))=\gamma(t)$.

\noindent By Lemma \ref{lemma:HenonSectors}, we know that $D\widetilde{H}$ expands
the horizontal length of vectors in $W^-$, so
\[
|\gamma'(t)|>C(x_{\gamma_1}(t))|\gamma_1'(t)|.
\]
We will compare the length of the curve $\gamma_1$ with the length
of $\gamma$. Note that $\gamma(t)$ is a just a horizontal line
segment, hence $|\gamma'(t)|=|x_1-x_2|,\mbox{ for all }t\in[0,1]$
and $l(\gamma)=|x_1-x_2|$.
\begin{eqnarray}\label{eq:len-gamma1}
l(\gamma_1)=\int_{0}^{1} |\gamma_1'(t)|dt <
\int_{0}^{1}\frac{1}{C(x_{\gamma_1}(t))}|\gamma'(t)|dt=
|x_1-x_2|\int_{0}^{1}\frac{1}{C(x_{\gamma_1}(t))}dt
\end{eqnarray}
Recall from Equation \ref{eq:eps2} that $C(x):=|1+(q+1)x^q|-m|x|^{2q}\geq 1+\epsilon_{1}
|x|^q$. 
We have $C(x)>1$ for all $(x,y)\in W^-$. Since $|x|<1$, we also have that 
\[
\frac{1}{C(x)}\leq\frac{1}{1+\epsilon_{1} |x|^q}\leq
1-\frac{\epsilon_{1}}{2}|x|^q.
\]
Recall also that for any $t\in[0,1]$
\[
x_{\gamma}(t)=\lambda x_{\gamma_1}(t)
\left(1+x_{\gamma_1}(t)^{q}+g_a(x_{\gamma_1}(t),z)/x_{\gamma_1}(t)\right)
\]
where $\big{|}1+x_{\gamma_1}(t)^{q}+g_a(x_{\gamma_1(t)},z)/x_{\gamma_1}(t)\big{|}<M_{1}$, as in Equation \ref{eq:M1}. By
combining these estimates we obtain
\begin{eqnarray*}
\int_{0}^{1}\frac{1}{C(x_{\gamma_1}(t))}dt <
1-\frac{\epsilon_{1}}{2}\int_{0}^{1}|x_{\gamma_1}(t)|^q dt &<&
1-\frac{\epsilon_{1}}{2M_{1}^{q}}\int_{0}^{1}|x_{\gamma}(t)|^q dt\\
&=& 1-\frac{\epsilon_{1}}{2M_{1}^{q}}\int_{0}^{1}|tx_1+(1-t)x_2|^q dt.
\end{eqnarray*}
Using that $|x'-x''|\leq l(\gamma_1)$ and Equation \ref{eq:len-gamma1} we get the desired inequality.
\qed

\begin{lemma}[\textbf{Technical estimate}]\label{lemma:technical} Let $q\geq 1$ be a natural number and $x_1,x_2\in \C$ be two complex
numbers, with $|x_2|\leq |x_1|$. Then
$
|x_1|^{q}\leq 2(q+1)\mathcal{I}(x_{1},x_{2}).
$
\end{lemma}
\proof If $x_1=0$ then $x_2=0$ and we have equality. Otherwise, set
$x=x_2/x_1$. Then $|x|\leq 1$ and we need to show that
\[
\frac{1}{2(q+1)}\leq
\int\limits_{0}^{1}\bigg{|}t\frac{x_2}{x_1}+(1-t)\bigg{|}^q
dt=\int\limits_{0}^{1}|tx+(1-t)|^q dt.
\]
For any $t\in[0,1]$ we have $|tx+(1-t)|\geq
|t|x|-(1-t)|=|t(1+|x|)-1|$. Let $u =t(1+|x|)-1$. Then $du =
(1+|x|)dt$ and
\begin{eqnarray*}
\int_0^1|t(1+|x|)-1|^q\, dt &=& \frac{1}{|x|+1}\int_{-1}^{|x|}
|u|^q\, du =
\frac{1}{|x|+1}\frac{|x|^{q+1}+1}{q+1}>\frac{1}{2(q+1)},
\end{eqnarray*}
since $0\leq|x|\leq1$.
\qed

Suppose for now that $a$ is small enough such that $N_a<\frac{1}{8q}$. This is similar to what we previously required for $N_a$. In Equation \ref{eq:Na-bound} we will impose another bound for $N_{a}$.

\begin{lemma}[\textbf{Step 2}]\label{lemma:Step2} 
$|x_2-x_3|< 4(q+1)N_{a} \mathcal{I}(x_{1},x_{2})|x_1-x_2|.$
\end{lemma}
\proof The geometric intuition behind the inequality is that the
curve $\gamma_2$ connecting $(x_1,y_1)$ and $(x_2,y_2)$ becomes
horizontal as $a\rightarrow0$, while the fibers $f_1$ and $f_2$
become vertical. The rigorous proof is outlined below.
From the proof of Proposition \ref{prop:hcone-in-petal} it follows that
$|y_1-y_2|<N_a|x_1-x_2|$, where $N_a\rightarrow 0$ as $a\rightarrow 0$.
The  curve
\[
t\rightarrow \left(\varphi_{2}(ty_1+(1-t)y_2), ty_1+(1-t)y_2\right), \
\ \ t\in [0,1]
\]
is vertical-like so in particular the horizontal distance is smaller
than the vertical distance and
\[
|\varphi_{2}(ty_1+(1-t)y_2) - \varphi_{2}(y_2)|<|ty_1 + (1-t)y_2 - y_2| =
t|y_1-y_2|,
\]
for $t\in [0,1]$. Using $\varphi_{2}(y_2)=x_2$ this gives
\[
|\varphi_{2}(ty_1+(1-t)y_2)|<|x_2|+t|y_1-y_2|<|x_2|+tN_a|x_1-x_2|.
\]
Hence
\begin{eqnarray*}
|x_2-x_3|&\leq&\int\limits_{0}^{1}\left|\frac{\partial}{\partial
t}\varphi_{2}(ty_1+(1-t)y_2)\right|dt\leq
\int\limits_{0}^{1}|y_1-y_2|\left|\varphi_{2}(ty_1+(1-t)y_2)\right|^{2q}dt\\
&\leq& |y_1-y_2| \left(|x_2|+N_a|x_1-x_2|\right)^{2q}\leq
N_{a}|x_1-x_2| \left(|x_2|+N_a|x_1-x_2|\right)^{2q}.
\end{eqnarray*}

Suppose without loss of generality that $|x_2|\leq |x_1|$ and $|x_{1}|<1$. From
the technical estimate Lemma \ref{lemma:technical} we get
\begin{eqnarray*}
|x_2-x_3|&<&N_{a}|x_1-x_2|(1+2N_a)^{2q}|x_1|^{2q}\\
&<& N_{a}|x_1-x_2|(1+2N_a)^{2q} \cdot 2(q+1)\mathcal{I}(x_{1},x_{2}).
\end{eqnarray*}
Since $N_a<\frac{1}{8q}$, we can use the following estimate
\[
2(q+1)(1+2N_{a})^{2q} < 2(q+1)\left(1+\frac{1}{4q}\right)^{2q}<4(q+1)
\]
to get  $|x_2-x_3|< 4(q+1)N_{a} \mathcal{I}(x_{1},x_{2})|x_1-x_2|$.
\qed

We now return to the proof of Theorem \ref{thm:distE}. 
We can use the triangle inequality in the $\D_{\rho}\times \{y_1\}$
disk to connect $|x_1-x_2|$ to the distance between the curves $f_1$
and $f_2$
\[
|x_1-x_2|-|x_2-x_3|\leq|x_1-x_3|\leq d(f_1,f_2).
\]
In Step 2 we showed that $|x_2-x_3|< 4(q+1)N_{a} \mathcal{I}(x_{1},x_{2})|x_1-x_2|$, so
\[
\left(1-4(q+1)N_{a} \mathcal{I}(x_{1},x_{2})\right)|x_1-x_2|<d(f_1,f_2).
\]
Using the bound on $|x_{1}-x_{2}|$ from Step 1, we get that 
\[
|x'-x''|<\frac{1-\frac{\epsilon_{1}}{2M_{1}^{q}} \mathcal{I}(x_{1},x_{2})}{1-4(q+1)N_a  \mathcal{I}(x_{1},x_{2})}\ d (f_1,f_2),
\]
where the quantity
\begin{equation}\label{eq:bigCL}
C= \frac{1-\frac{\epsilon_{1}}{2M_{1}^{q}} \mathcal{I}(x_{1},x_{2})}{1-4(q+1)N_a  \mathcal{I}(x_{1},x_{2})}<1,
\end{equation}
for $a$ small enough. Indeed, the constants $q, \epsilon_{1}$ and $M_{1}$ are independent of the
parameter $a$ whereas $N_a\rightarrow 0$ as $a \rightarrow 0$, so it
can be made small enough so that
\begin{equation}\label{eq:Na-bound}
N_a < \frac{\epsilon_{1}}{8(q+1)M_{1}^{q}}.
\end{equation}
The right hand side is a fixed constant, but this bound is not optimized. We get that
\[
d(g_{1},g_{2})\leq C d(f_{1},f_{2}),
\]
where $C<1$ depends on the distance between the curves and the $y-$axis.
This dependency is hidden in $\mathcal{I}(x_{1},x_{2})$. 

\vglue -.45cm
\qedof\textbf{of Theorem \ref{thm:distE}}

\section{The fixed point of a weakly contracting operator}\label{subsec:contraction}

In this section, we construct a function space $\mathcal{F}$ and a graph transform operator $F:\mathcal{F}\rightarrow \mathcal{F}$. We endow the space $\mathcal{F}$ with a metric induced by $\mu$ on the set $V$ and show that the operator $F$ is strictly (but not strongly) contracting. The key ingredients will be the invariance of vertical cones constructed in Section  \ref{sec:cones} under $DH^{-1}$ and the weak expansion of $DH$ in the horizontal cones. We use a generalization of the Banach fixed point theorem, due to Browder, to claim the existence of a unique fixed point $f^{*}$ of $F$.

\begin{lemma}\label{lemma:preim}
Let $(x,y)\in V\cap \overline{U^{+}}$ and $(x',y')=H^{-1}(x,y)$. If $|y'|<r$ then $(x',y')\in V$.
\end{lemma}
\proof The point $(x',y')$ is in $U^+\cup J^+$ hence it cannot lie in the sets that have been removed from $\D_r\times \D_r$ when constructing the set $V$ as they belong to the interior of $K^{+}$, as shown in Lemmas \ref{lemma:omega} and \ref{lemma:HenonSectors}. 

If $(x,y)$ is in $V$, then it belongs to $W^{-}_B$, $W^{-}_{B'}$ or $(U'-\D_{\rho'}(q_{0}))\times \D_r$. In the first case, if $(x,y)$ belongs to the repelling sectors $W^{-}_B$, then $(x',y')\in W^{-}_B \cup W^{-}_{B'} \subset V$. In the last two cases,  if $(x,y)$ belongs to $W^{-}_{B'}$ or $(U'-\D_{\rho'}(q_{0}))\times \D_r$, then for $a$ chosen small enough, we can assume that the disk of radius $2r|a|$ around $x$ is contained in $U$. 
The point
$
\hvec{x'}{y'}=\frac{1}{a}\hvec{y}{x-p(y/a)-a^2w}
$
belongs to $V$ because $|y'|<r$ (by hypothesis) and $y/a\in U'$. 
We can use the inequality $|x-p(y/a)-a^2w|<r|a|$ to show that $y/a\in U'$. Indeed, since  $a$ is chosen small enough
so that the disk of radius $r|a|+|w||a|^{2}<2r|a|$ around $x$ is still in $U$, it follows that $p(y/a)\in U$, hence $y/a\in U'$. Then $(x',y')\in U'\times \D_{r}$ and $(x',y')\in U^{+}\cup J^{+}$, hence $(x',y')$ belongs to $V$.
\qed

\begin{defn}
Let $L=\{(f(z),z),\ z\in \D_{r}\}\subset V$, be the graph of an analytic function $f:\D_{r}\rightarrow \D_{r}$. The analytic curve $L\subset V$ is \textit{vertical-like}, if the following conditions are met. Choose $(x,y)\in L$ and $(\xi,\eta)$  a tangent vector to $L$ at $(x,y)$. If $(x,y)\in B$, then $(\xi,\eta)$ belongs to the pull-back vertical cone $\mathcal{C}^{v,B}_{(x,y)}$ described in Definition \ref{def:cones-pB} using Definition \ref{def:vcones}.   If $(x,y)$ is outside $B''$ then $(\xi,\eta)$ belongs to the vertical cone $\mathcal{C}^{v,P}_{(x,y)}$ described in Definition \ref{def:conesP}. 
\end{defn}

\begin{lemma}\label{lemma: two-preimages}
Let $L$ be a vertical-like curve in $V\cap (U^{+}\cup J^{+})$. Then $H^{-1}(L)\cap V$ is the union of two vertical-like curves $L_{1}$ and $L_{2}$.
\end{lemma}
\proof
Since the curve $L$ is vertical-like, it is the graph of a holomorphic function $f:\D_{r}\rightarrow \D_{r}$, with $|f'(z)|<1$. Hence $L=\{(f(z),z), z\in \D_{r}\}$.
Then
\[
H^{-1}(L)=\left\{ \left(z/a,(f(z)-p(z/a)-a^{2}w)/a\right),\ z\in \D_{r}\right\}
\]
is an analytic curve whose horizontal folding points cannot belong to $V$. Suppose there is a folding point inside $V$. By construction of $V$, the first coordinate $z/a$ of the folding point must be bounded away from $0$ (independent of $a$). It follows that the equation $f'(z)-2z/a^{2}=0$ cannot have solutions inside $\D_{r}$, as $ (2/a)\cdot (z/a)$ gets arbitrarily large when $a$ is small enough, whereas $f'(z)$ remains bounded, because the curve is vertical-like.

Therefore the degree of the projection of $H^{-1}(L)\cap V$ on the second coordinate is constant. It is easy to see that the degree is $2$, by looking at the number of intersections of $H^{-1}(L)$ with the $x\mbox{-axis}$. The curve $L$ is vertical-like in $V$, hence it intersects $H(x\mbox{-axis})$ in exactly two points. Then $H^{-1}(L)$ intersects the $x\mbox{-axis}$ in two points.

Thus $H^{-1}(L)\cap \C\times \D_{r}$ is a union of two analytic curves $L_{1}$ and $L_{2}$. By Lemma \hyperref[lemma:preim]{\ref{lemma:preim}}, $L_{1}$ and $L_{2}$ are contained in $V$, hence $H^{-1}(L)\cap V$ is the union of two analytic curves $L_{1}$ and $L_{2}$.

Let us make one more remark about $L_{1}$ and $L_{2}$. Since $L$ is a vertical like disk in $U^{+}$, its projection on the first coordinate is almost constant and bounded away from $0$, the critical point of $p$. Let $\Delta$ be the image of the projection on the first coordinate. There are two holomorphic branches $p_{1}$ and $p_{2}$ of $p^{-1}$ defined on $\Delta$.  Let now $(f(z),z)$ be any point of $L$ such that $H^{-1}(f(z),z)\in V$. By Lemma \hyperref[lemma:preim]{\ref{lemma:preim}} it suffices to check that the condition $|(f(z)-p(z/a)-a^{2}w)/a|<r$ is met.  This condition means exactly that $z/a$ is $\bigO(a)$ close to either $p_{1}(f(z))$ or $p_{2}(f(z))$. The curves $L_{1}$ and $L_{2}$ correspond to different choices of the branch of $p^{-1}$.

By Lemmas \ref{prop:conesE} and \ref{prop:cones-Ve-product-metric}, the curves $L_{1}$ and $L_{2}$ are vertical-like, so they are graphs of functions on $\D_{r}$. The projection $pr_{2}: L_{1}\rightarrow \D_{r},\ pr_{2}(x,y)=y$ is a degree one covering map, hence by the Implicit Function Theorem,  $L_{1}$ is the graph of a holomorphic function $x=\phi(y)$. The map $\phi$ must also be injective, as the projection on the first coordinate $pr_{1}: L_{1}\rightarrow U'$, $pr_{1}(x,y)=x$ is injective.  
\qed

Choose $R>2$ as in the construction of the neighborhood $U$ of the Julia set $J_{p}$ from Equation \ref{eq:nbdU} and define the sequence of equipotentials $\gamma_{n}:\R/\Z\rightarrow \C$ of the parabolic polynomial $p$,
\begin{equation}\label{eq:gamma}
    \gamma_{n+1}(t)=\Phi_{p}\left(R^{1/2^{n+1}}e^{2\pi i t}\right) = p^{-1}(\gamma_{n}(2t)),
\end{equation}
where $\Phi_{p}$ is the inverse B\"ottcher isomorphism of $p$ \ref{eq:Phip}.
Note that $\gamma_{-1}(\R/\Z)\subset\partial U$ and $\gamma_{0}(\R/\Z)\subset\partial U'$.

\begin{defn}
Let $f_{0}:\s^{1}\times \D_{r}\rightarrow \partial V$ be the function $f_{0}(t,z)=(\gamma_{0}(t),z)$. The image of $f_{0}$ is the outer boundary of $V$ and it is contained in the escaping set $U^{+}$.
\end{defn}

We will construct a sequence of pull-backs of the map $f_{0}$ under the \He map.  This choice of $f_0$ will simplify the computations from Section \ref{subsec: conjugacy} where we establish the conjugacy of the \He map to a model map. For the definition of the function space we could start with any function $f_0:\s^1\times\D_r\rightarrow U^+\cap V$, $f_0(t,z)=(\phi_t(z),z)$, continuous with respect to $t$ and analytic with respect to $z$ such that $\phi_t(\D_r)$ is a vertical-like disk and $t\mapsto \phi_t(0)$ is homotopic to $\gamma_0(t)$. 

\begin{defn}\label{def:F}
Consider the space of functions:
\begin{equation*}
\mathcal{F} = \left\{ f_n:\s^{1}\times \D_{r} \rightarrow V\ |\ 
        f_0(t,z)=(\gamma_0(t),z),\ f_n(t,z)=F\circ f_{n-1}(t,z)\
        \mbox{for}\ n\geq 1
        \right\},
\end{equation*}
where the graph transform $F:\mathcal{F}\rightarrow\mathcal{F}$ is 
defined as
\[
 F(f)=\tilde{f},
\] where $ \tilde{f}\big{|}_{t\times\D_{r}}$ is the conformal map of the component
of $H^{-1}\left( f(2t\times\D_{r}) \right)\cap V$ ``homotopic to''
$f_{0}(t\times\D_{r})$, normalized via the Implicit
Function Theorem (the projection on the second coordinate).
\end{defn}
On the function space $\mathcal{F}$ we consider the metric
\begin{eqnarray*}
 d(f, g) = \sup_{t\in \s^{1}}\sup_{z\in \D_{r}} d_{\mu}\left(f(t,z),\ g(t,z)\right).
\end{eqnarray*}
where $d_{\mu}\left(f(t,z),\ g(t,z)\right)$ is the infimum of the
length of horizontal rectifiable paths $\tau:[0,1]\rightarrow V$
with $\tau(0)=f(t,z)$ and $\tau(1)=g(t,z)$. The length is
measured with respect to $\mu$, which is defined in Section \ref{sec:metric}. Note that $d_{\mu}$ is a metric.

\medskip
We will begin by describing the function space $\mathcal{F}$ and showing how the function $f_{1}$ is constructed. Proposition \ref{prop: F-welldef} explains the construction of the graph-transform $F$ and describes the main properties of the functions from the space $\mathcal{F}$. 

\begin{figure}[htb]
\begin{center}
\includegraphics[scale=0.93, bb = 200 450 400 720]{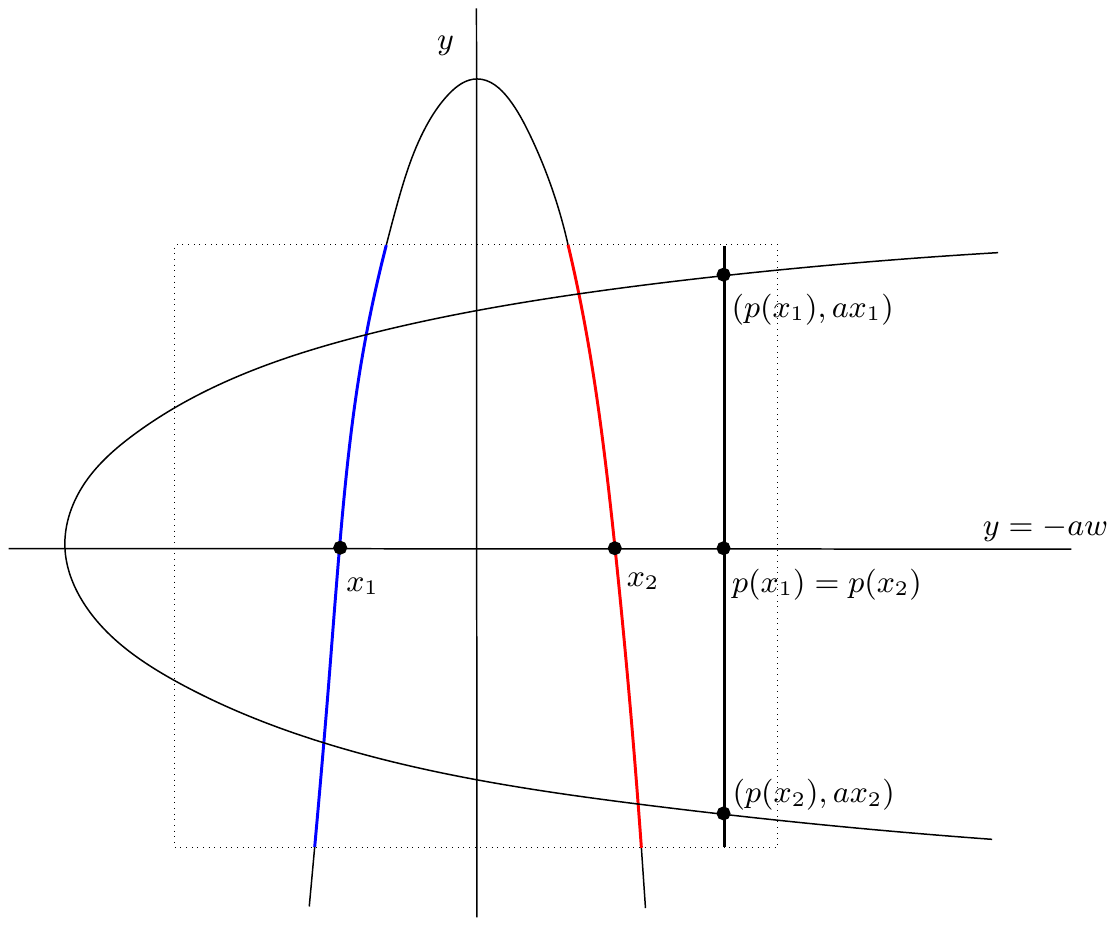}
\end{center}
\caption{The preimage of a fiber of $f_{0}$ in the neighborhood $V$. }
\label{pic:pbs_0}
\end{figure}

For any fixed $t\in \s^{1}$,  the set $f_{0}(2t\times\D_{r})$ is a vertical disk.
By Lemma \ref{lemma: two-preimages}, the preimage $H^{-1}(f_{0}(2t\times \D_{r}))\cap V$ is a disjoint union of two vertical-like disks (as shown in Figure \ref{pic:pbs_0}) that we would like to label as $t$ and $t+1/2$.

The preimage $H^{-1}(\gamma_{0}(2t),z)=(z/a, (\gamma_{0}(2t)-p(z/a)-a^{2}w)/a)$ belongs to $V$ if and only if the second component belongs to $\D_{r}$. The inequality  $|\gamma_{0}(2t)-p(z/a)-a^{2}w|<r|a|$ implies that the first coordinate $z/a$ is $\bigO(a)$ close to one of the two preimages of $\gamma_{0}(2t)$ under the polynomial $p$.

Denote by $C_{t}$ and $C_{t+1/2}$ the components of $H^{-1}\circ f_{0}(2t\times \D_{r})\cap V$ that cross the horizontal axis $y=-aw$ at $\left(\gamma_{1}(t),-aw\right)$, and respectively at $\left(\gamma_{1}(t+1/2),-aw\right)$. 

Notice that $pr_{2}:C_{t}\rightarrow \D_{r},\  pr_{2}(x,z)=z$ is a degree one covering map, hence, by the Implicit Function Theorem, $C_{t}$ is the graph of a holomorphic function $x=\phi_{t}(z)$. Let us define a new function $f_{1}:\s^{1}\times \D_{r} \rightarrow V$ as
$f_{1}(t,z):=(\phi_{t}(z),z)$. 

\begin{remarka} $f_{1}$ is homotopic to $f_{0}$ by construction, since $\gamma_{1}(t)$ and $p^{-1}(\gamma_{0}(2t))$ are homotopic. Moreover, since $a$ is small, $f_{1}(\s^{1}\times \D_{r})$ and $f_{0}(\s^{1}\times \D_{r})$ are disjoint.
\end{remarka}

\begin{prop}\label{prop: F-welldef}
    The map $F:\mathcal{F}\rightarrow\mathcal{F}$ is well defined. Choose any function $f\in \mathcal{F}$. For any $t\in\s^{1}$, $f(t\times \D_{r})$ is a vertical-like disk parametrized by the second coordinate. 
There exists $\varphi_t:\D_{r}\rightarrow \D_{r}$ analytic with respect to $z$ and $a$ and continuous with respect to $t$ such that $f(t,z)=(\varphi_{t}(z),z)$.
\end{prop}
\proof
We have $F\circ f_{0}(t,z)=f_{1}(t,z)$. Assume by induction that the functions $f_{n-1}, f_{n}:\s^{1}\times \D_{r} \rightarrow V $, $f_{n}(t,z)=(\varphi^{n}_{t}(z),z)$, $f_{n-1}(t,z)=(\varphi^{n-1}_{t}(z),z)$ have been constructed for $n\geq 1$ and let us show how to define $f_{n+1}$.

For each $t\in \s^{1}$, $H^{-1}(f_{n}(2t\times\D_{r}))\cap V$ is a union of two vertical like disks in $U^{+}$, $C_{t}$ and $C_{t+1/2}$. A choice of labeling is involved and we will first explain how this is done. Intuitively we would like to label by $t$ the disk which is closer to $f_{n}(t\times \D_{r})$ and by $t+1/2$ the disk which is closer to $f_{n}((t+1/2)\times \D_{r})$.  Let $(x_{1}, -aw)$ and $(x_{2}, -aw)$ be the two intersection points of $H^{-1}(f_{n}(2t\times\D_{r}))$ with the axis $y=-aw$. The points $(p(x_{1}), ax_{1})$ and $(p(x_{2}), ax_{2})$ belong to $f_{n}(2t\times \D_{r})$ and we would like to label them. The disk $f_{n-1}(2t\times\D_{r})$ contains two labeled points $b^{n}_{t}=H(f_{n}(t, -aw))$ and $b^{n}_{t+1/2}=H(f_{n}(t+1/2, -aw))$. Let $p_{1}$ and $p_{2}$ be two holomorphic branches of $p$ such that $p_{1}\circ p(\varphi_{t}^{n}(-aw))=\varphi_{t}^{n}(-aw)$ and $p_{2}\circ p (\varphi_{t+1/2}^{n}(-aw))=\varphi_{t+1/2}^{n}(-aw)$.   If $p_{1}\circ p(x_{1})=x_{1}$ then we label the point $(p(x_{1}), ax_{1})$ as $b^{n+1}_{t}$ and the component of $H^{-1}(f_{n}(2t\times\D_{r}))\cap V$ that intersects the axis $y=-aw$ at $(x_{1}, -aw)$ as $C_{t}$. Otherwise, if $p_{2}\circ p(x_{1})=x_{1}$, we label it as $C_{t+1/2}$. 

As before, the projection on the second coordinate $pr_{2}:C_{t}\rightarrow \D_{r},\  pr_{2}(x,z)=z$ is a degree one covering map, hence, by the Implicit Function Theorem, $C_{t}$ is the graph of a holomorphic function $x=\varphi_{t}^{n+1}(z)$. Let then $f_{n+1}$ be defined as \[f_{n+1}(t,z):=(\varphi^{n+1}_{t}(z),z).\] 
It is also easy to see that $\varphi_{t}^{n+1}$ is injective, by the definition of $H^{-1}$.  The function $f_{n+1}$ is holomorphic with respect to $z$ and $a$ and continuous with respect to $t$. 
\qed

\begin{prop}\label{prop:strict}
    The operator $F:\mathcal{F}\rightarrow \mathcal{F}$ is a strict contraction.
    \[
    d(F(f),F(g))<d(f,g), \mbox{ for any } f,g\in \mathcal{F}.
    \]
\end{prop}

The proof is an immediate consequence of the following proposition.
\begin{prop}\label{prop:inegal2t} Let $f,g\in \mathcal{F}$ and $t\in \s^{1}$. Then
\[
d\left(F\circ f(t\times \D_r),F\circ g(t\times
\D_r)\right)<C(f,g,t)\cdot d\left(f(2t\times \D_r),g(2t\times \D_r)\right),
\]
where $C(f,g,t)$ is a contraction factor which depends on the fibers $f(t\times\D_{r})$ and $g(t\times \D_{r})$ and $0\leq C(f,g,t)<1$.
\end{prop}
\proof
The delicate case is when the curves enter $B''$ and come close to $W^{s}_{loc}(\fq_{a})$. By Theorem \ref{thm: mu-expansion}, case (b), the derivative of the  \He map still expands in the horizontal direction but the expansion factor goes to $1$ as we approach the stable manifold. This case has already been carefully analyzed in Section \ref{sec:distance}.

The case where the curves are outside of a small neighborhood of the local stable manifold $W^{s}_{loc}(\fq_{a})$ can be treated as in the hyperbolic setting, because by Theorem \ref{thm: mu-expansion} the derivative of the \He map expands in the horizontal direction with a fixed expansion factor, independent of $a$.

Suppose that the fibers $f(2t\times \D_{r}), g(2t\times \D_{r}),  F\circ f(t\times \D_{r})$, $F\circ g(t\times \D_{r})$  belong to $V-B=(U'-\D_{\rho'}(q_{0}))\times \D_{r}$. We show that there exists a constant $C<1$ such that
\begin{equation}\label{eq:localineq}
\sup\limits_{z\in \D_{r}} d_{\mu}\left(F\circ f(t,z), F\circ g(t,z) \right)\leq C \sup\limits_{z\in \D_{r}}d_{\mu}\left(f(2t,z), g(2t,z) \right).
\end{equation}

Recall that $f(2t\times \D_{r}), g(2t\times \D_{r}), F\circ f(t\times \D_{r})$ and $F\circ g(t\times \D_{r})$ are  vertical-like disks  parametrized by the second coordinate. There exists conformal maps $\varphi_{1},\varphi_{2}:\D_{r} \rightarrow U'$ such that
$F\circ f(t,z)=(\varphi_{1}(z),z)$ and $F\circ g(t,z)=(\varphi_{2}(z),z)$. Let $z$ be any point in $\D_{r}$. Denote by $x'= \varphi_{1}(z)$, $x''= \varphi_{2}(z)$. Then
$(x_1,y_1)=H(x',z)=(p(x')+az+a^{2}w, ax')$ and $(x_2,y_2)=H(x'',z)=(p(x'')+az+a^{2}w, ax'')$. Finally, denote by $(x_3,y_1)$ the intersection point of the vertical-like curve $g(2t\times \D_{r})$ with the horizontal line $\C\times\{y_1\}$. The configuration is the same as in Figure  \ref{fig: FourFibers} (fibers named $f$ and $g$ here correspond to $f_1$ and $g_1$ on the picture).

 Let $\gamma$ be any horizontal curve in $V$ between the point $(x_1,y_1)$ on the curve $f(2t\times \D_{r})$ and the point $(x_2,y_1)$. There exists a horizontal curve $\gamma_1$ in $V$ between the point $(x',z)$ on the curve $F\circ f(t\times \D_{r})$ and $(x'',z)$ on the curve $F\circ g(t\times \D_{r})$ such that $\gamma_2=H(\gamma_1)$ is a curve inside the parabola $H(\C\times {z})$ linking $(x_1,y_1)$ to $(x_2,y_2)$ such that its projection on the plane $\C\times{y_1}$ is exactly the curve $\gamma$. The curves $\gamma$, $\gamma_{1}$ and $\gamma_{2}$ are just local variables here. 

The \He map expands the length of horizontal vectors, therefore by Definition \ref{def:Poincaremetric} and Proposition \ref{prop:cones-Ho-product-metric} there exists $k>1$ such that
\[
\mu_P(\gamma_2(s), \gamma_2'(s))>k\mu_P(\gamma_1(s),\gamma_1'(s)), \mbox{ for any } s\in[0,1].
\] 
Notice also that the tangent vector $\gamma_2'(s)$ belongs to the horizontal cone at $\gamma_2(s)$, so in particular, $\mu_P$ is just the Poincar\'e metric of the projection on the first coordinate, so \[\mu_P(\gamma_1(s),\gamma_1'(s))=\mu_P(\gamma(s),\gamma'(s))=\mu_U(pr_1(\gamma(s)),pr_1(\gamma'(s))).
\]
After passing to the infimum after the length of all horizontal curves $\gamma$
we conclude that 
\[
d_{\mu}((x_1,y_{1}),(x_2,y_{1}))>kd_{\mu}((x',z),(x'',z)),
\]
where the distance is measured with respect to the Poincar\'e metric of U.

Let $\psi: \D_{r} \rightarrow g(2t\times \D_{r})$ be the conformal isomorphism which parametrizes the fiber. Then $\psi(ax')=(x_3,y_1)$ and $\psi(ax'')=(x_2, y_2)$. The fiber is vertical-like, therefore
\[
d_{\mu}((x_2,y_{1}), (x_3,y_{1}))\leq \con |ax'-ax''|
\]
The points $x'$ and $x''$ are in $U'-\D_{\rho'}(q_{0})$ which is compactly contained in $U$, so we have
\[
    m_{1}|x'-x''|\leq d_{\mu}((x',z),(x'',z))\leq m_{2}|x'-x''|.
\]
In conclusion we have
$
d_{\mu}((x_2,y_{1)}, (x_3,y_{1})) \leq |a|\frac{\con}{m_1} d_{\mu}((x',z),(x'',z)).
$
By the triangle inequality it follows that
\begin{eqnarray*}
d_{\mu}((x_1,y_{1}),(x_2,y_{1}))-d_{\mu}((x_2,y_{1}),(x_3,y_{1}))&\leq& d_{\mu}((x_1,y_{1}),(x_3,y_{1}))\\
& \leq &\sup\limits_{z\in \D_{r}}d_{\mu}(f(2t,z),g(2t,z)).
\end{eqnarray*}
So by combining the previous inequalities we get
\[
d_{\mu}((x',z),(x'',z))\leq \frac{1}{(k- |a|\frac{\con}{m_1}) } \sup\limits_{z\in \D_{r}}d_{\mu}(f(2t,z),g(2t,z))
\]
The expansion factor $k$ is strictly bigger than $1$, so when $a$ is small, $C:=(k- |a|\frac{\con}{m_1})^{-1}$ is a contraction factor strictly less than $1$, which gives inequality \ref{eq:localineq}.
\qed

\begin{thm}[\textbf{Contracting map}]
There exists a monotonically increasing and right continuous function $h:[0,\infty)\rightarrow[0,\infty)$ such that $h(s)<s$ for each $s>0$ and
\[
d\left(F(f),F(g)\right)\leq h\left(d(f,g)\right),
\]
for any $f,g\in\mathcal{F}$.
\end{thm}
\proof Let $h:[0,\infty)\rightarrow[0,\infty)$ be
\[
h(s):=\sup\limits_{\substack{f,g\in \mathcal{F}, t\in S^1\\
d\left(f(2t\times\D_r),g(2t\times\D_r)\right)\leq s}} d(F\circ
f(t\times \D_r),F\circ g(t\times \D_r)).
\]
It is easy to see that $h$ is increasing and that $h(0)=0$. Moreover, by definition
\[
d\left(F(f),F(g)\right)\leq h\left(d(f,g)\right),
\]
for any $f,g\in\mathcal{F}$.

By Proposition \ref{prop:inegal2t} we know that
\begin{equation}\label{eq:C-local}
d\left(F\circ f(t\times \D_r),F\circ g(t\times
\D_r)\right)<C(f,g,t)d\left(f(2t\times \D_r),g(2t\times \D_r)\right),
\end{equation}
where $0\leq C(f,g,t)<1$ is the contraction factor. It follows that $h(s)\leq s$ for all $s\geq 0$. This right-hand limit $h(s+) := \lim_{\delta\searrow 0}h(s+\delta)$ exists everywhere since the function $h$ is monotonically increasing. We want to show that $h(s+)<s$ for all $s>0$.

Suppose that $h(s+)=s$ for some $s>0$. Let $(\delta_{n})_{n\geq 1}$ be a strictly decreasing sequence of positive numbers converging to $0$. For each $n$ there exists fibers $f_{n}$, $g_{n}$ and a $t_{n}\in \s^{1}$ such that
\begin{equation}\label{eq:sup-h}
d\left(F\circ f_{n}(t_{n}\times \D_r),F\circ g_{n}(t_{n}\times \D_r)\right)>h(s+\delta_{n})-\delta_{n}
\end{equation}
and where $d\left(f_{n}(2t_{n}\times \D_r),g_{n}(2t_{n}\times \D_r)\right)\leq s+\delta_{n}$.
This follows from the definition of $h(s+\delta_{n})$ as a supremum. In view of relation \ref{eq:C-local} we get that
\begin{eqnarray*}
h(s+\delta_{n})-\delta_{n}&<&d\left(F\circ f_{n}(t_{n}\times \D_r),F\circ g_{n}(t_{n}\times \D_r)\right)\\
&<& C_{n}d\left(f_{n}(2t_{n}\times \D_r),g_{n}(2t_{n}\times \D_r)\right)\leq C_{n} (s+\delta_{n})<s+\delta_{n},
\end{eqnarray*}
where $C_{n}:=C(f_{n},g_{n},t_{n})$ is a number as in Equation \ref{eq:C-local} above, with $0\leq C_{n}<1$ for every $n\geq 1$. Dividing both sides by $s+\delta_{n}$ and passing to the limit as $n\rightarrow \infty$ yields
\[
\frac{h(s+)}{s}=1\leq \lim_{n\rightarrow \infty}C_{n}\leq 1.
\]
Thus $\lim_{n\rightarrow \infty}C_{n}$ exists and is equal to 1. However, this can only happen if for all $n\geq n_{0}$ the fibers $f_{n}$ and $g_{n}$ belong to the normalizing tubular neighborhood of the semi-parabolic fixed point and the distance between the fibers is measured in the Euclidean metric (in fact the pull-back of the Euclidean metric under the normalizing map). Otherwise, the contraction factor $C(f_{n},g_{n},t_{n})$ is bounded by a uniform constant $K<1$.

The contraction factor $C_{n}$ is constructed explicitly in Section \ref{sec:distance}. It is of the form
\[ C_{n}= \frac{1-\alpha \mathcal{I}(x_{1,n},x_{2,n})}{1-\beta  \mathcal{I}(x_{1,n},x_{2,n})},
\]
where $\alpha,\beta$ are fixed constants with $0<\beta<\alpha$. The numbers $x_{1,n}$ and $x_{2,n}$ are the $x$-coordinates of two points that belong to the fibers $f_{n}(2t_{n}\times \D_{r})$, respectively $g_{n}(2t_{n}\times \D_{r})$. Recall that $ \mathcal{I}(x_{1,n},x_{2,n})=\int_{0}^{1}|tx_{1,n}+(1-t)x_{2,n}|^{q}dt$. If $C_{n}\rightarrow 1$ then $ \mathcal{I}(x_{1,n},x_{2,n})\rightarrow0$. In Lemma \ref{lemma:technical} we showed that $ \mathcal{I}(x_{1,n},x_{2,n})\geq \frac{1}{2(q+1)}\max(|x_{1,n}|^{q},|x_{2,n}|^{q})$, so $x_{1,n}\rightarrow 0$ and $x_{2,n}\rightarrow 0$.  But then $|x_{1,n}-x_{2,n}|\rightarrow 0$. It follows from Lemma \ref{lemma:x1x2} and the choice of $x_{1,n}$ and $x_{2,n}$ that $d\left(F\circ f_{n}(t_{n}\times \D_r),F\circ g_{n}(t_{n}\times \D_r)\right)\rightarrow 0$ as $n\rightarrow \infty$.

Passing to the limit in Equation \ref{eq:sup-h} yields $0\geq h(s+)=s$, thus $s=0$. Contradiction! Therefore $h(s+)<s$ for all $s>0$. The function $\tilde{h}:s\mapsto h(s+)$ is continuous from the right and verifies all properties of the function $h$. With a small abuse of notation we will consider this as the function $h$ from the hypothesis.
\qed

\begin{thm}[Browder \cite{Br}]\label{thm:Browder} Let $(X,d)$ be a complete metric space and suppose $f:X\rightarrow X$ satisfies
\[
	d(f(x),f(y))<h(d(x,y))\ \ \ \mbox{for all}\ \ x,y\in X,
\]
where $h:[0,\infty)\rightarrow [0,\infty)$ is increasing and continuous from the right such that $h(s)<s$ for all $s>0$.  Then $f$ has a unique fixed point $x^{*}$ and $f^{n}(x)\rightarrow x^{*}$ for each $x\in X$.
\end{thm}
\proof We will follow the proof from \cite{KS}. For a fixed $s>0$, the sequence $(h^{n}(s))_{n\geq 0}$ is monotone decreasing (not necessarily strictly) and bounded below, so it has a limit as $n\rightarrow \infty$. Since $h$ is continuous from the right, the sequence converges to a fixed point of $h$. But $0$ is the only fixed point of $h$, so $h^{n}(s)\rightarrow 0$ for each $s>0$.

Let $x_{0}\in X$ be fixed and consider $x_{n}=f^{n}(x_{0})$, $n=1,2, \ldots$. We can show inductively that $d(x_{n},x_{n+1})< h^{n}(d(x_{0},x_{1}))$ for all $n\geq 0$. Passing to the limit, we get that
\[
0\leq \lim_{n\rightarrow \infty}d(x_{n},x_{n+1})\leq \lim_{n\rightarrow \infty}h^{n}(d(x_{0},x_{1}))  =0.
\]
Thus $\lim_{n\rightarrow \infty}d(x_{n},x_{n+1})=0$. We now show that $(x_{n})_{n\geq 1}$ is Cauchy. Let $\epsilon>0$. Since $\epsilon-h(\epsilon)>0$, we can choose $n$ large enough so that
$d(x_{n},x_{n+1})<\epsilon-h(\epsilon)$.
Consider the ball $B(x_{n},\epsilon)=\{x\in X\ |\ d(x_{n},x)<\epsilon \}$ of radius $\epsilon$ around $x_{n}$. Let $z\in B(x_{n},\epsilon)$. Then
\begin{eqnarray*}
d(x_{n},f(z))&\leq& d(x_{n},f(x_{n}))+d(f(x_{n}),f(z))\\
 &\leq& d(x_{n},x_{n+1})+h(d(x_{n},z))\\
 &\leq& (\epsilon-h(\epsilon)) + h(\epsilon) = \epsilon.
\end{eqnarray*}
In the last step, we have used the fact that $h$ is increasing, so $d(x_{n},z)<\epsilon$ implies $h(d(x_{n},z)) \leq h(\epsilon)$.  Therefore $f:B(x_{n},\epsilon)\rightarrow B(x_{n},\epsilon)$. It follows that $d(x_{n},x_{n+m})<\epsilon$ for all $m\geq 0$. Thus our sequence is Cauchy, hence convergent since  $X$ is complete. Let $\lim_{n\rightarrow \infty}f^{n}(x)=x^{*}\in X$. Then $f(x^{*})=x^{*}$ since $f$ is continuous. Uniqueness of $x^{*}$ follows from the contractive condition.
\qed

A weaker version of this theorem was used in \cite{DH} to prove local connectivity of the Julia set of a parabolic polynomial. 
Note that for $h(s)=Ks$ with $0<K<1$, the theorem reduces to the classical Banach fixed point theorem.

Let $\mathcal{\overline{F}}$ be the completion of the space $\mathcal{F}$ in the $d$-metric defined above.

\begin{prop}\label{prop:f*}
    The map $F:\mathcal{\overline{F}}\rightarrow\mathcal{\overline{F}}$ has a unique fixed point $f^{*}$.
\end{prop}
\proof
The operator $F$ is contracting in the metric defined on $\mathcal{\overline{F}}$. The existence and uniqueness of a fixed point follows from the fixed point Theorem \ref{thm:Browder}. \qed

The fixed point $f^{*}$ is a continuous surjection $f^{*}:\s^{1}\times \D_{r}\rightarrow J^{+}\cap\overline{V}$. 
As defined in Section \ref{subsec: V-parabolic} by $\overline{V}$ we denote the set $V$ together with the local stable manifold $W^{s}_{loc}(\fq_{a})$ and its preimage $H^{-1}(W^{s}_{loc}(\fq_{a})) \cap B'$, which are both in the boundary of $V$.

\begin{prop} $Im(f^{*}) = J^{+}\cap \overline{V}$.
\end{prop}
\proof By Lemma  \ref{lemma: neighJ} we have
$
J^{+}\cap \overline{V}=\bigcap_{n\geq 0}H^{-\circ n}(\overline{V}\cap \overline{U^{+}} )
$.
By construction, $f_{0}(t,z)=(\gamma_{0}(t),z)$, so $f_{0}(\s^{1}\times \D_{r})$ is the outer boundary of $\overline{V}$ and is entirely contained in $U^{+}$. Moreover, $f^{*}$ was obtained as a limit of the functions $f_{n}:\s^{1}\times \D_{r}\rightarrow \overline{V}$, where $f_{n}(\s^{1}\times \D_{r})=H^{-1}(f_{n-1}(\s^{1}\times\D_{r}))\cap V$, so $f_{n}(\s^{1}\times\D_{r})$ is the outer boundary of the set $\bigcap_{0\leq k\leq n}H^{-\circ k}(\overline{V}\cap \overline{U^{+}})$. Hence $Im(f^{*})=\bigcap_{n\geq 0}H^{-\circ n}(\overline{V}\cap \overline{U^{+}})$.
\qed

\begin{prop} The fixed point $f^{*}:\s^{1}\times \D_{r}\rightarrow J^{+}\cap\overline{V}$ has the form
\[
f^{*}(t,z)=(\varphi_{t}(z),z),
\]
where $\varphi_{t}(z)$ is continuous with respect to $t$, holomorphic with respect to $z$ and $a$.
\end{prop}
\proof 
The fixed point $f^{*}(t,z)= (\varphi_{t}(z),z)$ 
is obtained as a uniform limit of the sequence $(f_n)_{n\geq 0}$, where $f_{0}(t,z)=(\gamma_{0}(t),z)$
and $f_n(t,z)=F^{\circ n}(f_{0})(t,z)=(\varphi^{n}_{t}(z),z)$, when $n \geq 1$. Each function $f_n$ is continuous with respect to $t$ and holomorphic with respect to $z$ and $a$ and therefore $f^*$ is also continuous with respect to $t$ and holomorphic in $z$ and $a$. Notice also that for any $t\in \s^1$,  $\varphi^{n}_{t}(z)$ is injective when $n\geq 1$, so $\varphi_t(z)$ will either be injective or constant by Hurwitz's theorem. 
\qed

\section{The conjugacy}\label{subsec: conjugacy}

In this section we analyze the properties of the fixed point $f^{*}$ in more detail and construct the conjugacies to a unique model map. 

Consider $f^{*}(t,z) = (\varphi_{t}(z),z)$, where $\varphi_{t}(z)$ is continuous with respect to $t\in \s^{1}$ and analytic with respect to $z\in \D_{r}$. Let $\sigma:\s^{1}\times \D_{r}\rightarrow \s^{1}\times \D_{r}$
\begin{equation}\label{eq:sigma}
\sigma(t,z) = \left( 2t, a\varphi_{t}(z) \right).
\end{equation}
For sufficiently small $|a|>0$ the map $\sigma$ is well-defined. We will see that is also open and injective.
Suppose the semi-parabolic \He map is written as in Equation \ref{eq:Hw}. The following theorem is an immediate consequence of our construction.

\begin{thm}\label{thm:conjugacy2} Let $p(x)=x^{2}+c_{0}$ be a polynomial with a parabolic fixed point of multiplier $\lambda=e^{2 \pi i p/q}$. There exists $\delta>0$  such that for all parameters $(c,a)\in \mathcal{P}_{\lambda}$ with $0<|a|<\delta$ the diagram
\[
\diag{\s^{1}\times \D_{r}}{J^{+}\cap \overline{V}}{\s^{1}\times \D_{r}}{J^{+}\cap \overline{V}} {f^{*}}{\sigma}{H_{c,a}}{f^{*}}
\]
commutes.
\end{thm}
\proof From the definition of $f^{*}$, we have that $H\circ f^{*}(t\times \D_{r})$ is compactly contained in $f^{*}(2t\times \D_{r})$. Thus we can write
\[
H\circ f^{*}(t,z)=\hvec{p(\varphi_{t}(z))+a^{2}w+az}{a\varphi_{t}(z)} = \hvec{\varphi_{2t}(a\varphi_{t}(z))}{a\varphi_{t}(z)} =  f^{*}\circ \sigma (t,z).
\]
The last equality follows from $ f^{*}\circ \sigma(t,z) = f^{*}(2t,a\varphi_{t}(z)) = (\varphi_{2t}(a\varphi_{t}(z)),a\varphi_{t}(z))$. Therefore $f^{*}$ semiconjugates $H$ on $J^{+}\cap \overline{V}$ to $\sigma$ on $\s^{1}\times \D_{r}$, as claimed.
\qed

\begin{lemma}\label{lemma:expansion} We have the following expansion for $\varphi_{t}(z)$
\begin{equation*}
\varphi_{t}(z) = \gamma(t) -\frac{1}{p'(\gamma(t))}az + \bigO(a^{2}).
\end{equation*}
\end{lemma}
\proof Consider the sequence $f_{n}(t,z)=F^{\circ n}(f_{0})(t,z)=(\varphi_{t}^{n}(z),z)$, for all $n\geq 1$, and  $f_{0}(t,z) = (\gamma_{0}(t),z)$. By construction we have that $H\circ f_{n+1}(t\times \D_{r})$ is compactly contained in $f_{n}(2t\times \D_{r})$, hence
\[
H\circ f_{n+1}(t,z)=\hvec{p(\varphi_{t}^{n+1}(z))+a^{2}w+az}{a\varphi_{t}^{n+1}(z)} = \hvec{\varphi_{2t}^{n}(a\varphi_{t}^{n+1}(z))}{a\varphi_{t}^{n+1}(z)}
\]
and in particular
\begin{equation}\label{eq:recurrence}
p(\varphi_{t}^{n+1}(z))+a^{2}w+az=\varphi_{2t}^{n}(a\varphi_{t}^{n+1}(z)).
\end{equation}

Consider the sequence of equipotentials $\gamma_{n}(t)$ as defined in equation \ref{eq:gamma}. Since the Julia set $J_{p}$ is connected, $p'(\gamma_{n}(t))$ does not vanish. Moreover, if $p$ is parabolic, $p'(\gamma(t))$ does not vanish either, where $\gamma$ is the Charat\'eodory loop of the parabolic polynomial $p$. We have the following two relations
\begin{eqnarray*}
\gamma_{n+1}(t) &=& p^{-1}(\gamma_{n}(2t))\\
\left(p^{-1}\right)' (\gamma_{n}(2t)) &=& \frac{1}{p'(\gamma_{n+1}(t))}.
\end{eqnarray*}

Note that for $n=0$, $p(\varphi_{t}^{1}(z))+a^{2}w+az = \gamma_{0}(2t)$ so for $a$ sufficiently small
the following expansion holds
\begin{eqnarray*}
\varphi_{t}^{1}(z) &=& p^{-1}\left(\gamma_{0}(2t) -az-a^{2}w\right) = p^{-1}(\gamma_{0}(2t))- \left(p^{-1}\right)' (\gamma_{0}(2t))az + \bigO(a^{2})\\ &=& \gamma_{1}(t) - \frac{az}{p'(\gamma_{1}(t))} +\bigO(a^{2}).
\end{eqnarray*}

We show by induction that for $n\geq 1$
\begin{equation*}
\varphi_{t}^{n}(z) = \gamma_{n}(t) -\frac{az}{p'(\gamma_{n}(t))} + \bigO(a^{2}).
\end{equation*}
Indeed, rearranging equation \ref{eq:recurrence} yields
\begin{eqnarray*}
\varphi_{t}^{n+1}(z) &=& p^{-1}\left(\varphi_{2t}^{n}(a\varphi_{t}^{n+1}(z)) -az-a^{2}w\right) \\
&=& p^{-1}\left(\gamma_{n}(2t) -\frac{a^{2}\varphi_{t}^{n+1}(z)}{p'(\gamma_{n}(2t))} -az + \bigO(a^{2}) \right) = p^{-1}\left(\gamma_{n}(2t) -az + \bigO(a^{2}) \right)\\
&=& p^{-1}(\gamma_{n}(2t))  - \left(p^{-1}\right)'(\gamma_{n}(2t))az + \bigO(a^{2})= \gamma_{n+1}(t) - \frac{az}{p'(\gamma_{n+1}(t))} + \bigO(a^{2}).
\end{eqnarray*}
Letting $n\rightarrow \infty$ we get the desired expansion for $\varphi_{t}(z)$.
\qed

In fact, since the polynomial $p$ is quadratic, $p'(\gamma(t))$ is just  $2\gamma(t)$ in the expansion of $\varphi_{t}(z)$.

\begin{prop}\label{prop:Prop0}
Let $p$ be hyperbolic or parabolic. For sufficiently small $|a|>0$ the map $\sigma$ is open and injective. Also $\sigma(\s^{1}\times \D_{r})\subset \s^{1}\times \D_{|a|r'}$, with $r'<r$.
\end{prop}
\proof If $p$ be hyperbolic or parabolic then there are no critical points in $J_{p}$ and there exists $\epsilon>0$ such that if $\xi_{1}\neq\xi_{2}\in J_{p}$ such that $p(\xi_{1}) =p(\xi_{2})$ then $|\xi_{1}-\xi_{2}|>\epsilon$. Thus when $p$ is hyperbolic or parabolic $|\gamma(t)-\gamma(t+1/2)|>\epsilon$ for $t\in\s^{1}$. From Lemma \ref{lemma:expansion} there exists $M>0$ such that $|\varphi_{t}(z)-\gamma(t)|<|a|M$ for all $t\in \s^{1}$ and $z\in \D_{r}$. Then for $|a|<\frac{\epsilon}{2M}$ the map $\sigma$ is injective. It is also open because locally it is a homeomorphism. The Julia set $J_{p}$ is inside a disk of radius 2 \cite{Bu} so $|\gamma(t)|<2$ and $|\phi_{t}(z)|<2+|a|M$. Since $r>3$, we can easily find $r'<r$ such that the image of $\sigma$ is inside $\s^{1}\times \D_{|a|r'}$. 
\qed
\begin{prop}\label{prop:Prop1} Consider $f^{*}(t,z) = (\varphi_{t}(z),z)$ and suppose
that $f^{*}(t_1,z_1)=f^{*}(t_2,z_2)$ for some $t_1,t_2\in \s^1$ and
$z_1,z_2\in \D_r$. Then $\varphi_{t_1}(z)=\varphi_{t_2}(z)$ for all
$z\in \D_r$.
\end{prop}
\proof If $f^{*}(t_1,z_1)=f^{*}(t_2,z_2)$ then $(\varphi_{t_1}(z_1),z_1)=(\varphi_{t_2}(z_2),z_2)$, hence $z_1=z_2$ and $\varphi_{t_1}(z_1)=\varphi_{t_2}(z_1)$.
Denote by  $s:\D_r\rightarrow \C$ the holomorphic function $s(z)=\varphi_{t_1}(z)-\varphi_{t_2}(z)$
and assume that $s(z)$ has an isolated zero  of order $m$ at $z_1$.

The functions $\varphi_{t_1}(z)$ and respectively $\varphi_{t_2}(z)$
were obtained as the limit of the uniformly convergent sequence of
holomorphic functions $\varphi_{t_1}^{n}(z)$ and respectively
$\varphi_{t_2}^{n}(z)$. By Hurwitz's theorem, there exists $\rho>0$
such that for sufficiently large $n>n_0$, the function
$\varphi_{t_1}^{n}(z)-\varphi_{t_2}^{n}(z)$ has exactly $m$ zeros in
the disk $|z-z_1|<\rho$. This is a contradiction, since by
construction $\varphi_{t_1}^{n}(z)\neq\varphi_{t_2}^{n}(z)$ for any $n\geq 0$ and $z\in \D_{r}$. 
Hence $z_1$ cannot be an isolated zero of the function $s$ on $\D_r$.
It follows that $s$ vanishes identically on $\D_r$ and so
$\varphi_{t_1}(z)=\varphi_{t_2}(z)$ for all $z\in\D_r$.
\qed

The fixed point $f^{*}(t,z)=(\varphi_{t}(z),z)$ depends on the
parameter $a$. We will use the notation
$f_a^{*}(t,z)=(\varphi_{t}(z,a),z)$ whenever we want to stress out
the dependence on $a$. Let $\delta>0$ be chosen as in Theorem \ref{thm:conjugacy2}.

\begin{prop}\label{prop:Prop2} Fix $z\in \D_r$ and $a'\in \D_{\delta}$ and assume that
$\varphi_{t_1}(z,a')=\varphi_{t_2}(z,a')$ for some $t_1,t_2\in
\s^1$. Then $\varphi_{t_1}(z,a)=\varphi_{t_2}(z,a)$ for any $a$ with
$|a|<\delta$.
\end{prop}
\proof Let $s:\D_{\delta}\rightarrow \D_r$ be the holomorphic
function $s(a)=\varphi_{t_1}(z,a)-\varphi_{t_2}(z,a)$. Denote by
$s_n$ the holomorphic functions
$s_n(a)=\varphi_{t_1}^n(z,a)-\varphi_{t_2}^n(z,a)$.

For any $n\geq 0$, and any $a$ with $|a|<\delta$, we
have $\varphi_{t_1}^n(z,a)\neq\varphi_{t_2}^n(z,a)$
by construction. Therefore $s_n(a)\neq 0$ for any $n\geq0$ and any $a$ with $|a|<\delta$.

The sequence $s_n$ converges uniformly to $s$ on $\D_{\delta}$.
By Hurwitz, $s$ has either no zeros on $\D_{\delta}$ or vanishes
identically on $\D_{\delta}$. Since we know that $s(a_1)=0$ it
follows that $s$ vanishes identically, thus $\varphi_{t_1}(z,a)=\varphi_{t_2}(z,a)$ for any $a$ with $|a|<\delta$. 
\qed

\begin{prop}\label{prop:Prop3} Consider $t_1\neq t_2\in \s^1$. The following statements are equivalent
\begin{itemize}
  \item [a)] $f_a^{*}(t_1,z)=f_a^{*}(t_2,z)$ for some $a$ with $|a|<\delta$ and
some $z\in \D_r$
  \item[b)]
$f_a^{*}(t_1,z)=f_a^{*}(t_2,z)$ for any $z\in \D_r$ and for any $a$
with $|a|<\delta$
\item [c)] $\gamma(t_1)=\gamma(t_2)$.
\end{itemize}
\end{prop}
\proof
By Propositions \ref{prop:Prop1} and \ref{prop:Prop2} we know that if $f_a^{*}(t_1,z)=f_a^{*}(t_2,z)$ for some $a$ with $|a|<\delta$ and
some $z\in \D_r$ then $f_a^{*}(t_1,z)=f_a^{*}(t_2,z)$ for any
$a\in \D_{\delta}$ and for any $z\in \D_r$. In particular when $a=0$ we must have
$f_0^{*}(t_1,z)=f_0^{*}(t_2,z)$. This is equivalent to
$(\gamma(t_1),z)=(\gamma(t_2),z)$, hence $\gamma(t_1)=\gamma(t_2)$.
\qed

\medskip

On $\s^{1}$ we have a natural equivalence relation $\sim_{p}$ given by the Thurston lamination of the polynomial $p$ as follows: $t_{1}\sim_{p}t_{2}$ whenever $\gamma(t_{1})=\gamma(t_{2})$. Then the set $\s^{1}/_{\sim_{p}}$ is homeomorphic to $J_{p}$ and the polynomial $p$ acting on its Julia set $J_{p}$ is topologically conjugate to the angle doubling map on $\s^{1}/_{\sim_{p}}$ as in \cite{Th} and \cite{Th1}.

This allows us to determine the equivalence classes of $f^{*}$. We define an equivalence relation $\sim$ on $\s^{1}\times \D_{r}$ so that $(t_{1},z)\sim(t_{2},z)$ whenever $\gamma(t_{1})=\gamma(t_{2})$. By Lemma \ref{lemma:expansion} $\varphi_{t}(z)$ can be written as
\[
\varphi_{t}(z)=\gamma(t)-\frac{az}{p'(\gamma(t))}+a^{2}\beta(t,z,a).
\]
In view of Proposition \ref{prop:Prop3} above, $\beta(t_{1}, z,a) = \beta(t_{2}, z,a)$ whenever $\gamma(t_{1}) = \gamma(t_{2})$.   Clearly $\sim$ is closed.  We would like to identify the quotient space $\s^{1}\times \D_{r}/_{\sim}$ with $J_{p}\times \D_{r}$ and the map $\sigma$ on $\s^{1}\times \D_{r}$ defined in Equation \ref{eq:sigma} with a similar map $\sigma_{p}$ acting on $J_{p}\times \D_{r}$.

Consider a map $\sigma_{p}: J_{p}\times \D_{r}\rightarrow J_{p}\times \D_{r}$ of the form
\begin{equation}\label{eq:sigma-p}
\sigma_{p}(\zeta,z) = \left(p(\zeta),a\left(\zeta-\frac{az}{p'(\zeta)}+a^{2}\beta\left(\gamma^{-1}(\zeta), z,a\right)\right)\right).
\end{equation}
It is well defined, in view of the discussion above.

The map $g:\s^{1}\times \D_{r}/_{\sim}\rightarrow J_{p}\times \D_{r}$, $g(t,z)=(\gamma(t),z)$ is a homeomorphism which makes the diagram
\begin{equation}\label{eq:g}
\diag{\s^{1}\times \D_{r}/_{\sim}}{J_{p}\times \D_{r}}{\s^1\times \D_{r}/_{\sim}}{J_{p}\times \D_{r}}
{g}{\sigma}{\sigma_{p}}{g}
\end{equation}
commute. The conjugacy follows directly from the fact that $p(\gamma(t))=\gamma(2t)$.

The map $\sigma_{p}$ on $J_{p}\times \D_{r}$ has the form
\begin{equation*}\label{eq:sigma-p1}
\sigma_{p}(\zeta,z) =\left(p(\zeta),a\zeta-\frac{a^{2}z}{p'(\zeta)}+\bigO(a^{3})\right).
\end{equation*}
and can be further conjugated to a solenoidal map
\[
\psi(\zeta, z)=\left(p(\zeta),a\zeta -\frac{a^{2}z}{p'(\zeta)}\right).
\]

For $|a|>0$ small enough $\sigma_{p}$ and $\psi$ are well-defined, open, and injective. Both maps depend on $a$ and we will use the notation $\psi_{a}$ to mark the dependence of $\psi$ on $a$, but we will use $\psi$ when there is no confusion. We will show in Lemma \ref{lemma:h-conj2} that for $0<|a|<\delta$ all $\psi_{a}$ are conjugate to each other. Fix $\epsilon$ so that $0<\epsilon<\delta$. Then $\psi_{a}$ and $\psi_{\epsilon}$ are conjugate and $\psi_{\epsilon}$ does not depend on $a$.

\begin{lemma}\label{lemma:h-conj1}
There is a homeomorphism $h:J_{p}\times \D_{r}\rightarrow J_{p}\times \D_{r}$ conjugating $\sigma_{p}$ to $\psi$.
\end{lemma}
\proof We first show that there exists a homeomorphism
\[h: J_{p}\times \D_{r} - \sigma_{p}(J_{p}\times \D_{r}) \rightarrow J_{p}\times \D_{r}-\psi(J_{p}\times \D_{r})
\]
which is the identity on the outer boundary $J_{p}\times \partial\D_{r}$ and given by the formula
\[ h(\zeta,z) = \psi\circ\sigma_{p}^{-1}(\zeta,z)
\]
on the inner boundary $\sigma_{p}(J_{p}\times \partial\D_{r})$. Define the space $\mathcal{H}$ of fiber homeomorphisms
\[
 J_{p}\times \D_{r} - \sigma_{p}(J_{p}\times \D_{r}) \rightarrow J_{p}\times \D_{r}-\psi(J_{p}\times \D_{r})
\]
that agree with $h$ on the boundary as a fiber bundle over $J_{p}$. Let $\zeta\in J_{p}$ and let $\mathcal{H}_{\zeta}$ be the fiber above $\zeta$ in $\mathcal{H}$. We know that $|p'(\zeta)|=2|\zeta|$ is bounded above and below since $J_{p}$ does not contain the critical point of $p$. The fiber above $\zeta$ in the range of the homeomorphism $h$ is a disk of radius $r$ with two disjoint disks of radius $\frac{r|a|^{2}}{2|\zeta|}$ removed, that is
\[
	\D_{r} - \bigcup_{\xi\in p^{-1}(\zeta)}\D_{\frac{r|a|^{2}}{2|\zeta|}}(a\xi).
\]
There are $d$ such disks removed if the polynomial has degree $d$. Similarly, the fiber above $\zeta$ in the domain is the disk $\D_{r}$ with two disjoint simply connected domains removed. These are topological disks of center $a\xi+\bigO(a^{3})$ and radius at most $\frac{r|a|^{2}}{2|\zeta|}+\bigO(|a|^{3})$, for all  $\xi\in p^{-1}(\zeta)$.

In $\mathcal{H}_{\zeta}$ we consider only those fiber homeomorphisms $h'$ which agree with $h$ on the boundary and which move all points by at most $\bigO(|a|^{3})$. Since the term $\bigO(|a|^{3})$ is much smaller compared to $\frac{r|a|^{2}}{2|\zeta|}$ when $a$ is small, there are no Dehn twists created as $\zeta$ moves on $J_{p}$. Therefore all such homeomorphisms are homotopic and this defines a preferred class of homeomorphisms. Note that $\mathcal{H}_{\zeta}$ is not empty. Furthermore, $\mathcal{H}_{\zeta}$ is contractible. This argument is in the same spirit as Lemma 6.8 in \cite{HOV2} and follows from a theorem of Hamstrom \cite{Ham} (which states that if $S$ is a compact surface with nonempty boundary -- in our case a disk with two disjoint disks removed -- then the  components of the group of homeomorphisms which are the identity on the boundary are contractible).

$\mathcal{H}$ is a locally trivial fiber bundle over $J_{p}$, with contractible fibers. A fiber bundle with contractible fibers over a paracompact base has a continuous section. Hence there exists a map
$s:J_{p}\rightarrow \mathcal{H}$, $s(\zeta)=h_{\zeta}$, which associates to each $\zeta$ a homeomophism $h_{\zeta}$, so that the choice is continuous with respect to $\zeta$. Set $h$ to be $s$. We now extend $h$ on the inner levels by the dynamics, so we are able to construct a homeomorphism
\begin{equation*}
h: J_{p}\times \D_{r} -\bigcap_{n\geq 0} \sigma_{p}^{\circ n}(J_{p}\times \D_{r}) \rightarrow J_{p}\times \D_{r}-\bigcap_{n\geq 0}\psi^{\circ n}(J_{p}\times \D_{r})
\end{equation*}
which conjugates $\sigma_{p}$ to $\psi$. Furthermore, we  extend to the Cantor set (in each fiber) by continuity.
\qed

\begin{lemma}\label{lemma:h-conj2}
There exists a homeomorphism $h_{a,\epsilon}:J_{p}\times \D_{r}\rightarrow J_{p}\times \D_{r}$ conjugating $\psi_{a}$ to $\psi_{\epsilon}$.
\end{lemma}
\proof We need to consider the space of homeomorphisms $\mathcal{H}$ and construct a preferred class of homeomorphisms. The proof is the same as the proof of Lemma \ref{lemma:h-conj1} above. The same idea was also used in Lemma 5.5 from \cite{R}.
\qed

Consider the linear change of variables $(\zeta,z)\mapsto (\zeta,az)$. For $|a|>0$ this  conjugates $\psi_{a}:J_{p}\times \D_{r} \rightarrow J_{p}\times \D_{r}$ to a map $\psi_{a}':J_{p}\times \D_{r'} \rightarrow J_{p}\times \D_{r'}$, where $r'=r/|a|$ and
\begin{equation}\label{eq:fp-model}
\psi_{a}'(\zeta, z)=\left(p(\zeta),\zeta -\frac{a^{2}z}{p'(\zeta)}\right).
\end{equation}
Similarly $\psi_{\epsilon}$ is conjugate to $\psi_{\epsilon}'$. Note that all these maps depend on the polynomial $p$.  When $p$ is hyperbolic, Lemma \ref{lemma:h-conj2} is Proposition 6.13 from \cite{HOV2}, and as we have seen, the situation is not very different when $p$ is parabolic.
The map $\psi_{\epsilon}'$  is the same model map that was used in  \cite{HOV2} in understanding \He maps that are small perturbations of hyperbolic polynomials. 

We now have all the ingredients to complete the proof of the theorems described in the introduction.

\medskip
\noindent \textbf{Proof of Theorem  \ref{thm:Parabolics-Jp}.}
The proof follows directly from Theorem \ref{thm:conjugacy2} and from Lemma \ref{lemma:h-conj1}. 
\qed

Based on the construction of the set $V$ from Section \ref{subsec: V-parabolic}, the Julia set $J$ of the \He map is $J = \bigcap_{n\geq 0} H^{\circ n}(J^+\cap \overline{V})$. Let $\Sigma := \bigcap_{n\geq 0}\sigma^{\circ n}(\s^{1}\times \D_{r})$. Then $\Sigma$ is a (dyadic) solenoid for $0<|a|<\delta$ and in view of Theorem \ref{thm:conjugacy2}, Proposition \ref{prop:Prop3} and the  above discussion, we can present $J$ as a quotiented solenoid, $J\simeq\Sigma/_{\sim}$. Therefore $J$ is connected. More directly, we can regard $J$ as
\begin{equation*}
	J \simeq  \bigcap_{n\geq 0}\psi_{\epsilon}^{\circ n}(J_{p}\times \D_{r}).
\end{equation*}

\noindent \textbf{Proof of Theorem  \ref{cor:J}.}
The proof follows directly from the model of $J$ described above and Corollary \ref{cor:J*} below. 
\qed

\begin{thm}\label{thm:conjugacy3} Let $p(x)=x^{2}+c_{0}$ be a polynomial with a parabolic fixed point of multiplier $\lambda=e^{2 \pi i p/q}$. There exists $\delta>0$  such that for all parameters $(c,a)\in \mathcal{P}_{\lambda}$ with $0<|a|<\delta$ there exists a homeomorphism $g^{*}$ which is continuous with respect to $\zeta$ and analytic in $a$ and $z$ and which makes the diagram
\[
\diag{J_{p}\times \D_{r}}{J^{+}\cap \overline{V}}{J_{p}\times \D_{r}}{J^{+}\cap \overline{V}} {g^{*}}{\sigma_{p}}{H_{c,a}}{g^{*}}
\]
commute.
\end{thm}
\proof The homeomorphism $g^{*}$ is a composition between $f^{*}$ and the inverse of the homeomorphism $g$ defined above, in Equation \ref{eq:g}. The map $\sigma_{p}$ is given in Equation \ref{eq:sigma-p}.
\qed

\begin{cor}\label{cor:J*} The Julia set $J$ equals $J^{*}$, the closure of the saddle periodic points.
\end{cor}
\proof
The Julia set $J$ is homeomorphic to a quotiented solenoid. Since the periodic points are dense in the solenoid, we get that $J$ is the closure of the periodic points of the \He map. Let  $x_{a}\in J$ be a periodic point of period $k$ of the \He map $H_{a}$, different from the semi-parabolic fixed point $\fq_{a}$. The periodicity of $x_{a}$ induces a periodicity on the disks that foliate $J_{p}\times \D_{r}$, namely there exists a periodic point $\zeta\in J_{p}$, $p^{\circ k}(\zeta)=\zeta$  of the parabolic polynomial $p$ such that $x_{a}\in g^{*}(\zeta\times
\D_{r})$ and $\sigma_{p}^{\circ k}(\zeta\times \D_{r})$ is compactly contained inside $\zeta\times
\D_{r}$. Note that $\zeta\neq q_{0}$, where $q_{0}$ is the parabolic fixed point of $p$. The conjugacy map $g^{*}(\zeta,z)$ is holomorphic with respect to $z$, so the stable multipliers of the \He map coincide with the stable multipliers of the map
\[
 \sigma_{p}(\zeta,z)=\left(p(\zeta), a\zeta-\frac{a^{2}z}{p'(\zeta)}+\bigO(a^{3})\right).
\]
Let $\lambda^{s/u}$ be the eigenvalues of $DH^{\circ k}_{x_{a}}$. Then $\lambda^{s}=\bigO(a^{2k})$ and $\lambda^{u} = (p^{\circ n})'(\zeta)+\bigO(a)$. The function $g^{*}$ is holomorphic with respect to $a$, so the disks that foliate $J^{+}\cap \overline{V}$ move holomorphically with $a$. The point $x_{a}$ moves holomorphically with $a$ and we have $x_{a}\rightarrow\zeta$ as $a\rightarrow0$.

The polynomial Julia set $J_{p}$ is the closure of the repelling periodic points \cite{M}. By the Fatou-Shishikura inequality \cite{S}, a polynomial of degree $d\geq 2$ has at most $d-1$ non-repelling cycles. Since $p$ is quadratic and has a parabolic fixed point $q_{0}$, all other periodic cycles are repelling. Therefore $|(p^{\circ n})'(\zeta)|>1$. Clearly, when $a$ is small, $|\lambda^{u}|>1$ and $|\lambda^{s}|<1$, so the periodic point $x_{a}$ is a saddle point of the \He map.

Let $\delta$ be as in Theorem \ref{thm:Parabolics-Jp}. We show that the periodic point $x_{a}$ is saddle. It is easy to see that $|\lambda^{s}|<1$, so we only need to show that $|\lambda^{u}|>1$.
Assume that $|\lambda^{u}|=1$ for some parameter $a_{0}$ with $0<|a_{0}|<\delta$. Then we can perturb $a_{0}$ so that $|\lambda^{u}|$ becomes strictly smaller than $1$. Otherwise $1/|\lambda^{u}|$ would have a local maximum at $a_{0}$, which is not possible. Thus we can find a parameter $a$ close to $a_{0}$ for which $x_{a}$ is a sink, and as such it must belong to the interior of $K^{+}$ and not to $J^{+}$; contradiction.
It follows that all periodic points are saddles, except the semi-parabolic fixed point, hence $J=J^{*}$.
\qed

\section{Conclusions}
Let $\mathcal{P}_{\lambda}^{n}$  be the set of parameters $(c,a)\in\C^{2}$ for which the \He map $H_{c,a}$ has a cycle of period $n$ with one multiplier $\lambda$ a root of unity. In other words, the $n^{th}$ iterate $H_{c,a}^{\circ n}$ of the \He map  has a fixed point $\fq$ such that its derivative $DH_{c,a}^{\circ n}(\fq)$ has eigenvalues $\lambda$ and $\mu=(-1)^n a^{2n}/\lambda$. 

When $n=1$ the curve $\mathcal{P}_{\lambda}^{1}$ is just the curve $\mathcal{P}_{\lambda}$ from Equation \ref{eq:P-lambda} which has a nice global characterization. We also have a nice description for $n=2$. However, when $n\geq 3$, it is hard to give an explicit formula for the curve $\mathcal{P}_{\lambda}^{n}$. In any case, $\mathcal{P}_{\lambda}^{n}$ is an algebraic set that intersects the parametric line $a=0$ in a discrete set of points. Suppose that $(c_0,0)$ is such a point of intersection, then $c=c(a)$ is locally a function of the parameter $a$ when $(c,a) \in \mathcal{P}_{\lambda}^{n}$ and $|a|<\delta$ is sufficiently small.
The results stated in Section \ref{sec:intro} hold for small perturbations inside the curve $\mathcal{P}_{\lambda}^{n}$ of the quadratic polynomial $x\mapsto x^2+c_0$ with a parabolic $n$-cycle of multiplier $\lambda$.

The technique presented in this paper and the results from Section \ref{sec:intro} can also be easily generalized to \He maps that are small perturbations (inside appropriate algebraic sets analogous to the curves $\mathcal{P}_{\lambda}^n$) of a polynomial $p$ of degree $d\geq 2$ whose critical points are attracted either to attractive or parabolic fixed points (or cycles). 
The local model space for $J^{+}$ is $J_{p}\times\D_{r}$ and the model map acting on it is
$
\psi(\xi,z)=\left(p(\xi),a\xi-\frac{a^{2}z}{p'(\xi)}\right)
$.
The map $\psi$ is again well defined because the Julia set $J_{p}$ contains no critical points.


\vspace{0.75cm}
\begin{tabularx}{\textwidth}{X  p{0.05\textwidth} X}
  {\sc Remus Radu} & & {\sc Raluca Tanase} \\
 Institute for Mathematical Sciences & &  Institute for Mathematical Sciences\\
 Stony Brook University & & Stony Brook University \\
 Stony Brook, NY 11794-3660& & Stony Brook, NY 11794-3660 \\
 e-mail: rradu@math.sunysb.edu& & e-mail: rtanase@math.sunysb.edu\\
\end{tabularx}

\end{document}